\documentclass[11pt]{article}
\usepackage{latexsym,amsmath,amssymb,epsf,amsthm,times,amscd}
\usepackage[T1]{fontenc} 
\usepackage[utf8]{inputenc} 
\usepackage[usenames]{color}
\usepackage{epsf}
\usepackage{graphicx}
\usepackage{graphics}
\usepackage{subcaption}
\usepackage{hyperref}

\newfont{\bb}{msbm10 at 11pt}
\newfont{\bbp}{msbm10 at 9pt}

\def\r{\hbox{\bb R}}

\def\b{\hbox{\bb B}}

\def\s{\hbox{\bb S}}
\def\sp{\hbox{\bbp S}}

\newcommand{\Max}{\textup{Max}}
\newcommand{\U}{\mathcal{U}}

\newcommand{\abs}[1]{\left\vert #1 \right\vert}
\newcommand{\set}[1]{\left\{#1\right\}}
\newcommand{\meta}[2]{\langle #1,#2 \rangle }

\newcommand{\pscalar}[2]{\left\langle #1, #2 \right\rangle}

\numberwithin{equation} {section}

\begin{document}
	\theoremstyle{plain}\newtheorem{lemma}{Lemma}[section]
	\theoremstyle{plain}\newtheorem{definition}{Definition}[section]
	\theoremstyle{plain}\newtheorem{theorem}{Theorem}[section]
	\theoremstyle{plain}\newtheorem{proposition}{Proposition}[section]
	\theoremstyle{plain}\newtheorem{remark}{Remark}[section]
	\theoremstyle{plain}\newtheorem{corollary}{Corollary}[section]
	\begin{center}
		\rule{15cm}{1.5pt} \vspace{.6cm}
		
		{\Large \bf Extremal domains in $\mathbb{S}^2$:\\[3mm]Geometric and Analytic methods} \vspace{0.4cm}
		
		\author[Jos\'{e} M. Espinar$\mbox{}^\dag$ and Diego A. Mar\'{i}n$\mbox{}^\ddag$
		
		\vspace{0.3cm} \rule{15cm}{1.5pt}
	\end{center}
	
	\vspace{.3cm}
	
	\noindent $\mbox{}^\dag$ Department of Geometry and Topology and Institute of Mathematics (IMAG), University of Granada, 18071, Granada, Spain; e-mail:
	jespinar@ugr.es\vspace{0.2cm}
	
	\noindent $\mbox{}^\ddag$ Department of Geometry and Topology and Institute of Mathematics (IMAG), University of Granada, 18071, Granada, Spain; e-mail: damarin@ugr.es
	
	\vspace{.3cm}
	
\begin{abstract}
In this article, we study domains $\Omega \subset \mathbb{S}^2$ that support positive solutions of the overdetermined problem 
\[
\Delta u + f(u,|\nabla u|)=0 \quad \text{in } \Omega,
\]
subject to the boundary conditions $u=0$ on $\partial\Omega$ and $|\nabla u|$ being locally constant along $\partial\Omega$. 
We refer to such domains as $f$--extremal domains.

In the first part of the paper, we extend the moving plane method in $\mathbb{S}^2$ and show that if an $f$--extremal domain $\Omega$ contains a simple curve of maximum points of $u$, then both $\Omega$ and $u$ are either rotationally symmetric or antipodally symmetric. 
Using the Alexandrov reflection method, we establish an analogous symmetry result for properly embedded constant mean curvature (CMC) surfaces with capillary boundaries that contain a simple curve of minimum distance to the origin (a neck).

In the second part, we strengthen these conclusions for specific nonlinearities, including the eigenvalue problem ($f(x)=\lambda x$), the Serrin problem ($f(x)=\lambda x + c$), harmonic domains ($f(x)=c$), and nonlinearities of the form $f(x)=\lambda x + x^\beta$, for constants $\lambda \geq 2$, $c \geq 0$, and $\beta \in (0,1)$. 
In these cases, we prove that the domain $\Omega$ must be rotationally symmetric.

Throughout the paper, we restrict our attention to the analytic setting in order to simplify the exposition and highlight the main ideas.
\end{abstract}

{\bf Key Words:} Overdetermined elliptic problems, moving plane method, capillary surfaces.

{\bf MSC2020:}  35Nxx, 35Bxx, 53Axx.

\section{Introduction}\label{Introduction}

The study of Overdetermined Elliptic Problems (OEPs) has garnered significant attention in the fields of differential geometry and partial differential equations (PDEs). These problems typically seek to determine whether a given domain \(\Omega\) within a Riemannian manifold \((\mathcal{M}, g)\) can support a solution to a particular elliptic PDE while simultaneously satisfying both Dirichlet and Neumann boundary conditions. More precisely, let \(\Omega \subset \mathcal{M}\) be a domain where the following OEP holds:
\begin{eqnarray} \label{OEP}
    \left\{
    \begin{array}{llll}
        \Delta{u} + f(u, |\nabla u|) = 0  &\mathrm{in}~\Omega,\\
        u > 0  &\mathrm{in}~\Omega, \\
        u = 0  &\mathrm{on}~\partial\Omega,\\
        \langle\nabla{u},\eta\rangle = \alpha_i  &\mathrm{on}~\Gamma_i \subset \partial \Omega,
    \end{array}
    \right.
\end{eqnarray}
where \(\alpha_i \leq 0\) is a constant on each connected boundary component \(\Gamma_i \subset \partial \Omega\), and \(f\) is a Lipschitz function depending on \(u\) and \(|\nabla u|\). Here, \(\Delta\) denotes the Laplacian associated with the Riemannian metric \(g\)  and $\eta$ the outer conormal to $\Omega$. The imposition of both Dirichlet and Neumann boundary conditions in \eqref{OEP} often imposes strong geometric constraints on the domain \(\Omega\) and the solution \(u\), giving rise to notable rigidity and symmetry phenomena. Domains that admit such solutions are referred to as \(f\)-\textit{extremal domains} and a solution to \eqref{OEP} is denoted by \((\Omega, u)\) (cf. \cite{Ros}). 

The study of OEPs has a rich and extensive history, beginning with foundational results in classical potential theory and extending into contemporary geometry and analysis. A seminal result in this area was established by J. Serrin \cite{Se}, who extended Alexandrov's reflection method, initially developed for the study of constant mean curvature (CMC) surfaces \cite{Ale}, to the context of PDEs. Serrin demonstrated that if a solution to the classical torsion problem in the Euclidean space $\mathbb{R}^n$ for the Laplace equation (i.e., \(f \equiv 1\)) satisfies both Dirichlet and Neumann boundary conditions, the domain must necessarily be a ball, and the solution must exhibit radial symmetry. This remarkable result, known as Serrin's Theorem, revealed the profound connections between OEPs and geometric methods. Shortly thereafter, H. F. Weinberger \cite{We} provided an alternative proof of Serrin's Theorem, employing purely analytical techniques via $P-$functions. Weinberger’s method underscored the power of maximum principles and comparison methods, highlighting an analytic pathway to understanding the rigidity of solutions to overdetermined problems. These complementary approaches—geometric via Alexandrov’s reflection and analytic via $P-$functions—have inspired extensive research in exploring the rigidity and symmetry properties of $f-$extremal domains (see \cite{CWZ}, and the references therein, for an excellent state-of-the-art overview.).

In recent decades, the study of OEPs has expanded significantly, incorporating a variety of geometric tools. One important approach to such problems is the moving plane method, which has been extended to space forms like the closed hemisphere \(\textup{cl}(\mathbb{S}^n_+)\) and hyperbolic space \(\mathbb{H}^n\), as explored by Kumaresan and Prajapat \cite{KP}. Additionally, the theory of CMC hypersurfaces has had a profound influence on the analysis of \(f\)-extremal domains (cf. \cite{Ber, EFM, EMao, Rei2, RRS} and references therein). Notably, Hélein et al. \cite{Hel} introduced a Weierstrass-type representation for \(0\)-extremal domains (often referred to as \textit{exceptional domains}) in the Euclidean plane, drawing a connection between the geometry of such domains and minimal surfaces. This connection was further extended by Traizet \cite{Trai}, underscoring the deep ties between exceptional domains and minimal surfaces. In the spherical setting, Espinar and Mazet \cite{EM} classified simply connected \(f\)-extremal domains in \(\mathbb{S}^2\) for a broad class of functions \(f\). Their work parallels Nitsche’s celebrated theorem \cite{JNit} on capillary CMC disks in the Euclidean unit ball. Capillary surfaces naturally arise in fluid mechanics as the interface between two fluids or between a fluid and air, subject to surface tension forces. In the context of geometric analysis, these surfaces provide a rich framework for studying minimal and constant mean curvature (CMC) surfaces, especially in the setting of overdetermined boundary problems. Recently, Cavalcante and Nunes \cite{Cav} examined the stability of extremal domains in \(\mathbb{S}^n\), utilizing techniques from the study of capillary CMC hypersurfaces in \(\mathbb{B}^n\), as established by Ros and Vergasta \cite{RosVer}. These developments further solidify the connection between overdetermined problems and capillary surfaces.

While geometric methods have played a central role in this area, analytic techniques—such as the P-function method introduced by Weinberger—have also been essential in advancing the field, especially in situations where geometric methods encounter limitations. This approach has been extended to general ambient spaces (cf. \cite{And,CV} and references therein), and a wider class of differential equations (cf. \cite{FarVal, Gao} and references therein). Another powerful analytical tool is the Implicit Function Theorem, first employed by Sicbaldi \cite{Sic} to construct bifurcating families of extremal domains. This method has since been widely applied to generate examples of \(f\)-extremal domains in various manifolds (cf \cite{FM,FMW2, Kam,Ruiz} and references therein). Further developments involve comparison techniques via pseudo-radial functions, leading to classification results in various contexts. (cf. \cite{ABBM,ABM,Chr} and references therein).

Consequently, the existence and classification of \(f\)-extremal domains remain active areas of research, with both geometric and analytic techniques being deployed to reveal new properties. These developments demonstrate how the interplay between geometry and analysis continues to uncover the intricate relationship between domain shape and the behavior of solutions to overdetermined elliptic problems. 

In this paper, our primary objective is to further explore the geometry of solutions to \eqref{OEP} in \(\mathbb{S}^2\) by combining the moving plane method with analytic methods, such as detailed gradient estimates, to provide a comprehensive study of $f$-extremal domains in \(\mathbb{S}^2\) with disconnected boundaries

\subsection{Statement of the results and organization of the paper}

We will investigate analytic solutions \((\Omega, u)\) to the overdetermined elliptic problem \eqref{OEP} in order to avoid technical complications. Specifically, we restrict our attention to \textit{finite domains} \(\Omega \subset \mathbb{S}^2\), meaning that \(\mathbb{S}^2 \setminus {\rm cl}(\Omega) = \bigcup_{i=1}^k D_i\), where \(D_i\) are disjoint open topological disks with regular (analytic) boundaries \(\Gamma_i\), and \(f \) is an analytic function. Therefore, regularity theory (cf.~\cite{GT}) implies that $u$ is analytic up to the boundary. We denote 
\[
\Max(u) := \{ p \in \Omega : u(p) = \max_{\Omega} u \},
\]
the set of maximum points of the function \(u\). Since \(u\) is real analytic, \(\Max(u)\) is a compact
real analytic (in particular, subanalytic) subset of \(\Omega\). Moreover, as \(u\) is not constant,
\(\Max(u)\) has empty interior. Hence, by the \L{}ojasiewicz structure theorem \cite{Kr}, \(\Max(u)\) admits
a finite stratification by real analytic submanifolds of dimensions \(0\) and \(1\). In particular,
\(\Max(u)\) is the union of finitely many isolated points and finitely many real analytic
one--dimensional strata (real analytic arcs/curves), which may meet at finitely many singular
points (vertices). In particular, a simple (embedded) closed curve $\gamma\subset \Max(u)$ is called an \emph{isolated component} if there exists an open neighborhood $\U\subset \Omega$ of $\gamma$ such that $\Max(u)\cap \U=\gamma$. (The same notion applies to curves contained in ${\rm Min}(u)$). Note that $\gamma$ must then be an embedded regular curve; otherwise, any neighborhood would intersect the lower--dimensional strata of $\textup{Max}(u)$.

Let \(\mathcal{I}: \mathbb{S}^2 \to \mathbb{S}^2\) be an isometry of the sphere. We introduce the following notions:

\begin{itemize}

\item 
$(\Omega,u)$ is {\it rotationally symmetric} if there exist a point $p_0\in\mathbb{S}^2$ and constants 
$0\leq a<b\leq\pi$ such that $\Omega=d_{p_0}^{-1}([a,b])$ and 
$u=U\circ d_{p_0}$ for some continuous function $U\in\mathcal{C}([a,b])$, 
where $d_p:\mathbb{S}^2\to\r^+$ denotes the distance function to $p\in\mathbb{S}^2$, 
that is, $d_p(q)=\mathrm{dist}_{\mathbb{S}^2}(q,p)$ for all $q\in\mathbb{S}^2$.

\item 
$(\Omega,u)$ is {\it antipodally symmetric} if $\mathcal{A}(\Omega)=\Omega$ and 
$u(p)=u(\mathcal{A}(p))$ for all $p\in\Omega$, where $\mathcal{A}:\mathbb{S}^2\to\mathbb{S}^2$ 
is the antipodal map $\mathcal{A}(p)=-p$.

\item 
Two pairs $(\Omega,u)$ and $(\tilde\Omega,\tilde u)$ are {\it congruent} if there exists an isometry 
$\mathcal{I}\in{\rm Iso}(\mathbb{S}^2)$ such that $\mathcal{I}(\Omega)=\tilde\Omega$ and 
$u(\mathcal{I}(p))=\tilde u(p)$ for all $p\in\Omega$. 
In this case, we write $(\Omega,u)\equiv(\tilde\Omega,\tilde u)$.

\end{itemize}

Our first main result is the following structural symmetry theorem.

\begin{quote}
{\bf Theorem A:} {\it Let $(\Omega, u)$, $\Omega \subset \s^2$ a finite domain, be an analytic solution to \eqref{OEP}. Assume that $\Max(u)$ contains a simple closed curve $\gamma$ that is an isolated component of $\Max(u)$. Then, $(\Omega ,u)$ is either antipodally or rotationally symmetric. In the latter case, $\Omega$ is either a disk or an annulus.}
\end{quote}

Theorem~A shows that the structure of the set of maximum points $\Max(u)$ strongly constrains the geometry of the associated $f-$extremal domain $\Omega$, without any assumption on the sign of the function $f$. Its proof is purely geometric and differs substantially from the arguments developed in \cite{EMa} for the special case $f(x)=2x$.

It combines an adapted moving plane method with a lifting argument inspired by classical works on constant mean curvature surfaces (cf.~\cite{Me}). Unlike the traditional moving plane method on the $2-$sphere, which typically requires that the complement of the domain contains an entire geodesic \cite{KP} (although a half-geodesic suffices to apply the method), our approach only requires that the complement of the domain contains two antipodal points. Finally, we emphasize that the techniques developed here are flexible enough to extend the symmetry result to domains with multiple boundary components and mixed Dirichlet–Neumann boundary conditions. This phenomenon is further illustrated by a result \`a la Sirakov \cite{Sir} (see Lemma~\ref{LemSirakov}). In particular, Lemma~\ref{LemSirakov} extends the Kumaresan--Prajapat Theorem \cite{KP} under the sole assumption that a connected component of the complement of $\textup{cl}(\Omega)$ contains two antipodal points.

A particularly relevant consequence of Theorem~A arises when the nonlinearity
$f$ satisfies a natural positivity condition.

\begin{quote}\label{CorThA}
{\bf Corollary A.} {\it Let $(\Omega, u)$, $\Omega \subset \s^2$ a finite domain, be an analytic solution to \eqref{OEP} and assume that $f(x,0)>0$ for all $x>0$. Then, either $\textup{Max}(u)$ contains only finitely many isolated points or the solution $(\Omega ,u)$ is either antipodally or rotationally symmetric. In the latter case, $\Omega$ is an annulus.}
\end{quote}

Corollary~A applies directly to several classical and widely studied models without additional topological assumptions, including:
the eigenvalue problem ($f(x)=\lambda x$), the Serrin problem ($f(x)=\lambda x+1$), harmonic domains in $\mathbb{S}^2$ ($f(x)=1$), and the Allen--Cahn equation ($f(x)=x(1-x^2)$). 

\begin{figure}[h!]
	\centering
	\begin{subfigure}[b]{0.4\linewidth}
		\includegraphics[width=\linewidth]{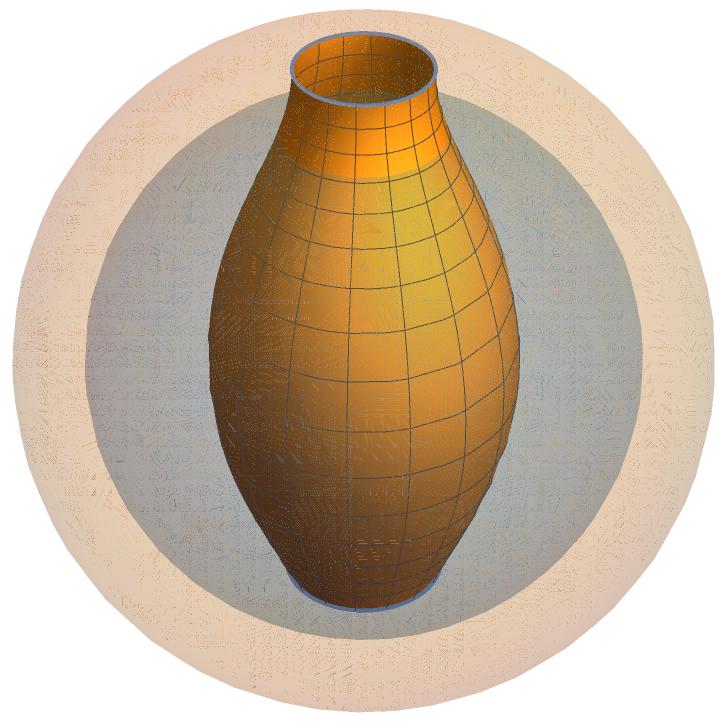}
		\label{fig:UnduloidCapillaryBoundaries}
	\end{subfigure}
	\begin{subfigure}[b]{0.4\linewidth}
		\includegraphics[width=\linewidth]{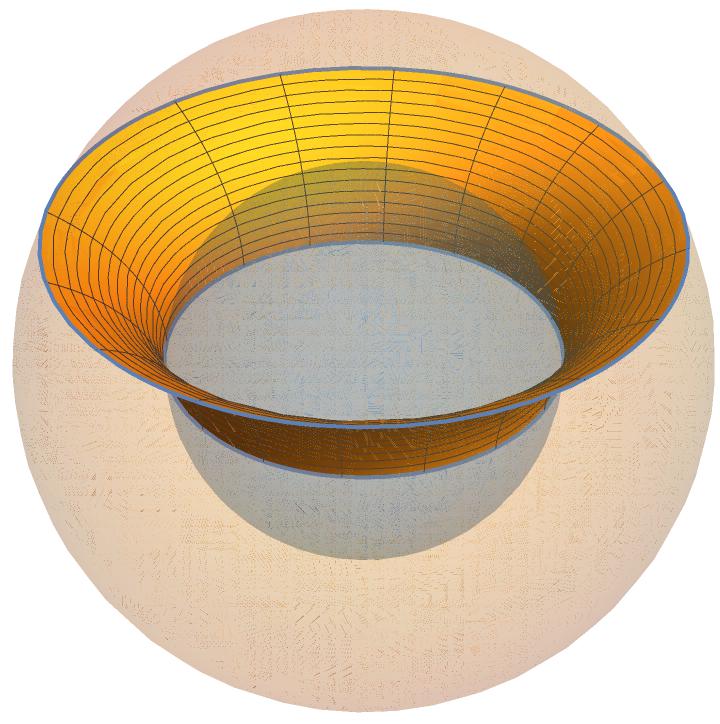}
		\label{fig:NodoidCapillaryBoundaries}
	\end{subfigure}
	\caption{We present two examples of rotationally symmetric CMC annuli with capillary boundaries. {\bf Left:} There exists an unduloid intersecting the unit sphere and a smaller sphere at different angles. {\bf Right:} A nodoid intersects the unit sphere orthogonally and is tangent to a smaller sphere.}
	\label{fig:CMCWithCapillaryBoundaries}
\end{figure}

We next show that the geometric techniques developed above can be extended to characterize constant mean curvature surfaces with capillary boundaries that contain a neck. A properly embedded (analytic) surface \(\Sigma \subset \mathbb{R}^3 \setminus \set{\bf 0}\) is said to have \textit{capillary boundaries} if each boundary component intersects from the inside a sphere centered at the origin at a constant angle, possibly with different radii and angles (cf. Definition \ref{surfacesWithFreeBoundaries} and Figure \ref{fig:CMCWithCapillaryBoundaries}). A {\it neck} $\gamma\subset\Sigma \subset \mathbb{R}^3 \setminus \set{\bf 0}$ is defined as a closed (analytic) curve along which the distance function to the origin in $\mathbb{R}^3$, 
${\bf d}_\Sigma(p):=|p|$ for $p\in\Sigma$, attains its absolute minimum, that is,
\[
\gamma\subset{\rm Min}({\bf d}_\Sigma):=\{p\in\Sigma : {\bf d}_\Sigma(p)=\min_\Sigma {\bf d}_\Sigma\}.
\]

With this terminology, in Section~\ref{SectGA}, we state the following geometric counterpart of Theorem~A, which constitutes the core rigidity statement of the paper in the setting of constant mean curvature surfaces.

\begin{quote}
{\bf Theorem B:} {\it Let $\Sigma \subset \mathbb{R}^3 \setminus \set{\bf 0} $ be a properly embedded (analytic) CMC annulus with capillary boundaries, $\partial \Sigma= \zeta ^- \cup \zeta ^+$, such that $\zeta ^- \subset \mathbb{S}^2 (r_-)$ and $\zeta ^+ \subset \mathbb{S}^2 ( r_+)$ for some $0<  r_- \leq r_+$. Set $0< r_0 := {\rm min}_{\Sigma} \, {\bf d}_\Sigma < r_-$ and suppose that 
\begin{enumerate}
\item[(i)] $\Sigma$ contains a neck $\gamma\subset\Sigma\cap\s^2(r_0)$ which is an isolated component of ${\rm Min}({\bf d}_\Sigma)$.

\item[(ii)] Let $\Sigma ^\pm \subset \Sigma \setminus \gamma$ be the connected component such that $\partial \Sigma ^\pm = \zeta ^\pm \cup \gamma$, $i=1,2$. The cone over $\gamma$, 
$C(\gamma):=\set{ \lambda \cdot p \in \mathbb{R}^3 \, : \, \lambda > 0 \text{ and } p \in \gamma}$, separates $\Sigma ^-$ and $\Sigma ^+$. 

\item[(iii)]  ${\rm int}(\Sigma ^\pm) \subset \textup{cl}(\b ^3 (r_\pm)) \setminus \b ^3 (r_0)$.	
\end{enumerate}
Then, $\Sigma$ is either antipodally or rotationally symmetric. In the latter case, $\Sigma$ is a part of either a catenoid, if $H=0$, or a Delaunay surface, if $H\neq 0$. Moreover, when $\Sigma$ is antipodally symmetric, then $r_-= r_+$.}
\end{quote}

In Section~\ref{SectGA}, we discuss the hypotheses of Theorem~B in their proper geometric context in detail, clarify the role and necessity of assumptions {\it (i)--(iii)}, and explain how Theorem~B extends to surfaces with more complicated topology. We also show how \cite[Theorem~B]{EMa} can be recovered as a particular case of Theorem~B. It is worth stressing that the techniques developed in \cite[Theorem~B]{EMa} are not sufficient to establish Theorem~B in the generality considered here.

The geometric techniques developed in this paper for classifying extremal domains apply to a wide range of functions \( f \). However, this method does not suffice to avoid antopodally symmetric solutions. To achieve a stricter classification, additional assumptions on the function \(f\) in \eqref{OEP} are necessary.
\begin{quote}
{\bf Theorem C:} {\it Let $(\Omega, u)$, $\Omega \subset \s ^2$ a finite domain, be an analytic solution to \eqref{OEP}, where we assume that either $f\equiv 1$ or $f\in\mathcal{C}^{\omega}$ satisfies $f(x)\geq x f'(x)$ and $f'(x)\geq 2$ for all $x>0$.

Then, either $\textup{Max}(u)$ contains only finitely many isolated points or the solution $(\Omega ,u)$ is rotationally symmetric. In the latter case, $\Omega$ is an annulus.
}
\end{quote}

Theorem~C applies to several classical and relevant situations, including the
eigenvalue problem ($f(x)=\lambda x$), the Serrin problem
($f(x)=\lambda x+1$), harmonic domains in $\mathbb{S}^2$ ($f(x)=1$), as well as
nonlinearities of the form $f(x)=\lambda x+x^\beta$ for $\lambda\ge2$ and
$\beta\in(0,1)$.

The proof of this result relies on an analytic method that generalizes the
approach developed in \cite{EMa} for the special case $f(x)=2x$, which itself builds on ideas introduced in~\cite{ABM}. In \cite[Section~1.5]{ABM}, the authors explicitly ask whether their method extends to more general semilinearities $f$ in the Euclidean setting.
Theorem~C shows that this is indeed the case on the $2$--sphere.
In a forthcoming work, we will show that the same strategy applies in Riemannian manifolds satisfying a lower bound on the Ricci curvature as well.

The main idea consists in introducing a family of rotationally symmetric {\it model solutions} to \eqref{OEP} (see Section~\ref{SectModel}) and developing a comparison argument that relates the geometry of the level sets of a general solution of \eqref{OEP} to those of a model solution.  To this end, we introduce the $\overline{\tau}$--function and use it to derive sharp gradient estimates for solutions of overdetermined problems in $\mathbb{S}^2$. 
These estimates play a crucial role in controlling the behavior of solutions near the boundary and near the set $\Max(u)$, and are the key ingredient in the proof of rotational symmetry.

The comparison framework developed in Section~\ref{SectAM} also clarifies the role and necessity of conditions imposed on the nonlinearity $f$, which naturally emerge in the second part of the paper, where the proof of Theorem~C is carried out.

\subsection{Concluding Remarks}

We have adopted the terminology Serrin problem for $f(x) = \lambda x +1$ (cf. \cite{CV}) and harmonic domains for $f(x) =1$ (cf. \cite{Savo}) to distinguish both cases in $\mathbb{S}^2$. However, the terminology $f(x) =1$ for the Serrin problem also appears in the literature (cf. \cite{FM}).

One can generalize Theorem A by replacing $\Delta u + f(u, |\nabla u|)$ in~\eqref{OEP} with a fully nonlinear, uniformly elliptic PDE (cf. \cite{BirDem} and references therein) of the form $F(\nabla^2 u, \nabla u, u) = 0$, provided that $F$ is invariant under isometries of $\mathbb{S}^2$ and the difference of any two solutions satisfies the maximum principle (including the Strong Maximum Principle, Hopf's Lemma, and Serrin's Corner Lemma).

Theorem B can be extended to a broader class of surfaces, such as those with positive constant extrinsic curvature \(K \). More generally, we can consider elliptic special Weingarten (ESW) surfaces in \(\mathbb{R}^3\), which satisfy a relation of the form \(H = \mathcal{F}(H^2 - K)\), where \(\mathcal{F}\) is a smooth function defined on an interval \({\cal J} \subseteq [0,\infty)\) with \(0 \in {\cal J}\), such that \(4t\mathcal{F}'(t)^2 < 1\) for all \(t \in {\cal J}\). A key insight is that Alexandrov's reflection principle and Nitsche type Theorem remain valid for ESW-surfaces (cf. \cite{AEG,EspFer}), thereby allowing symmetry results to be established even in these generalized contexts.

The sharpness of the characterizations for $f$-extremal domains (Theorem A) and capillary CMC surfaces (Theorem B) is not entirely clear. For Theorem A, we believe it is possible to construct a non-rotationally symmetric but antipodally symmetric $f$-extremal domain containing a curve of maximum points for some non-linearity $f$. In \cite{FMW}, Fall, Minlend and Weth construct a family of non-rotational annular harmonic domains ($f(x)=1$), however, this family only contains a finite number of maxima. Regarding Theorem B, Fernández, Hauswirth and Mira \cite{FHM} constructed a family of embedded capillary minimal annuli (extended to free boundary CMC annuli \cite{CMF}) foliated by spherical curvature lines. In this family, all surfaces are non-rotational, and some are antipodally symmetric. However, $\Sigma$ cannot contain a neck, as this would imply tangency with a sphere, which is excluded in their construction (see \cite[Remark 2.4]{CMF}).

Theorem~A uses analyticity only through an analytic continuation argument. By contrast, Lemma~\ref{LemSirakov} (à la Sirakov) remains valid under the weaker hypotheses that $\Omega$ has $\mathcal{C}^2$ boundary, $u\in\mathcal{C}^2(\overline{\Omega})$, and $f$ is Lipschitz. Likewise, the regularity assumptions in Theorem~C and, more generally, throughout Section~\ref{SectAM}, can be substantially relaxed.

\section{Geometric Method: Extended moving plane method}\label{SectGM}

Let $\mathbb{R}^3$ denote the usual Euclidean space, where $(x,y,z)$ represent Cartesian coordinates and $\pscalar{\cdot}{\cdot}$ denotes the Euclidean scalar product. In this section, we extend the classical moving plane method to prove that any analytic solution \((\Omega, u)\) to \eqref{OEP}, \(\Omega \subset \mathbb{S}^2\) a finite domain, is either rotationally symmetric or antipodally symmetric, provided that the level set \(\textup{Max}(u) := \set{p \in \Omega \, : \, \, u (p) = M}\) contains at least one simple closed curve that is an isolated component of $\textup{Max}(u)$. Set \(M:= u_{\textup{max}}\) the maximum value of \(u\). We first observe that:

\begin{lemma}\label{lemmaTheoremA1}
Let $(\Omega,u)$ be an analytic solution to \eqref{OEP}, $\Omega \subset \s^2$ a finite domain, such that $\textup{Max}(u)$ contains at least a closed curve $\gamma \subset \Omega$. If $\gamma$ is antipodally symmetric, then $(\Omega,u)$ is antipodally symmetric.
\end{lemma}
\begin{proof}
Let $\mathcal A : \s ^2 \to \s ^2$ be the antipodal map defined by $\mathcal A (p) =-p$, and consider $\tilde u := u \circ \mathcal A$ and $\tilde \Omega = \mathcal A (\Omega)$. Then, $( \tilde \Omega , \tilde u)$ is an analytic solution to \eqref{OEP}. 
	
Since $\gamma \subset \Omega \cap \tilde \Omega $ is invariant under the antipodal map, and $u (p) = M = u_{\textup{max}}$ for all $p \in \gamma$, it follows that $u=\tilde u$ along $\gamma$. Moreover, because $u$ attains its maximum along $\gamma$, we have $|\nabla u | =0$ along $\gamma$ and, similarly, $|\nabla \tilde u | =0$ along $\gamma$. Thus, $\nabla u = \nabla \tilde u$ along $\gamma$. 
	
	Therefore, $u$ and $\tilde u$ are analytic solutions to the same analytic PDE \eqref{OEP}, and their values and first derivatives coincide along $\gamma$. By the uniqueness guaranteed by the Cauchy–Kovalevskaya Theorem (cf. \cite{Kr}), it follows that $u = \tilde u$ in a neighborhood of $\gamma$. This extends to the entire domain by analyticity, that is, $(\Omega ,u) $ is antipodally symmetric.
\end{proof}

To illustrate the extended moving plane method that we develop, we derive the following result \`a la Sirakov:

\begin{lemma}\label{LemSirakov}
Let $\Omega \subset \s^2$ be a finite domain, $\s ^2 \setminus {\rm cl}(\Omega) = \bigcup _{i=1}^k D_i $ where each $D_i$ is an open topological disk. Write $\partial \Omega = \Gamma _0 \cup \Gamma _M$ with $\partial D_i=\Gamma_i \subset \Gamma _0$ for $i=1,\ldots,m$ and $\partial D_i=\gamma_i \subset \Gamma _M$ for $i=m+1,\ldots,k$ being analytic curves, and let $(\Omega ,u)$ be an analytic solution to
\begin{eqnarray*} 
    \left\{
    \begin{array}{llll}
        \Delta{u} + f(u, |\nabla u|) = 0  & \mathrm{in}~\Omega \subset \mathbb{S}^2, \\
        0 < u \leq M  & \mathrm{in}~\Omega, \\
        u = 0  & \mathrm{on}~\Gamma _0 \subset \partial \Omega ,\\
        u = M  & \mathrm{on}~\Gamma _M \subset \partial \Omega ,\\
        \langle \nabla u, \eta \rangle = \alpha_i  & \mathrm{on}~\Gamma_{i} \subset \Gamma _0,\\
        \langle \nabla u, \eta \rangle = 0  & \mathrm{on}~\gamma_i \subset \Gamma _M,
    \end{array}
    \right.
\end{eqnarray*}
where $\alpha_i \leq 0$ are constants for $i \in \set{1, \dots, m}$, and $\eta $ is the outer conormal along $\partial \Omega$. Assume moreover that each $\gamma _i \subset \Gamma _M$ is an isolated component of $\Max (u)$. Suppose that we have one of the following two cases
\begin{itemize}
	\item[(i)] There exists a point $x \in \s^2$ such that $x$ and $-x$ both belong to $D_i$, for some $i \in \{m+1,\ldots,k\}$. In this case, $\Gamma_0$ could be empty.
	\item[(ii)] There exists a point $x \in \s^2$ such that $x$ and $-x$ both belong to $D_i$, for some $i \in \{1,\ldots,m\}$ and $|\alpha _j | \leq |\alpha _i|$ for all $j \in \{1, \ldots , m \}$ and $j \neq i$. In this case, $\Gamma_M$ could be empty.
\end{itemize}
Then, $(\Omega ,u)$ is rotationally symmetric.
\end{lemma}

\begin{proof}
	We shall prove the first case $(i)$ of the above lemma. Case $(ii)$ will follow with the same argument.
	
	Thus, without loss of generality, assume that $x,-x\in D_{m+1}\setminus \{{\bf n},{\bf s}\}$, here ${\bf n}$ and ${\bf s}$ are the north and south poles. Let $\beta _x \subset \s^2$ be the unique great circle passing through $x$ and ${\bf n}$ and denote by $\beta_x^+$ (resp. $\beta_x^-$) the closed geodesic semi-arc in $\mathbb{S}^2 \subset \r ^3$ joining $x$ and $-x$, passing through the north pole ${\bf n}$ (resp. south pole ${\bf s}$). Let $P_x$ be the unique plane in $\mathbb{R}^3$ that contains $\beta _x$ and define $P_x^\pm$ as the closed half-plane in $\mathbb{R}^3$ containing $\beta_x^\pm$ whose boundary is the line $\mathcal{L}(x) = \{ t x : t \in \mathbb{R} \}$. For a rotation $R_{x,\theta}: \mathbb{R}^3 \to \mathbb{R}^3$ fixing the line $\mathcal{L}(x)$ by an angle $\theta$, denote the rotated half-planes as $P^\pm_{x,\theta} = R_{x,\theta}(P^\pm_x)$. Finally, denote by $\mathcal{R}_{x,\theta}: \mathbb{R}^3 \to \mathbb{R}^3$ the reflection across the plane $P_{x,\theta} = P^-_{x,\theta} \cup P^-_{x,\theta + \pi}$. 

	Consider the following radial graph associated to $u$:
\[
\Sigma:=\{(1+M-u(p))\,p\in\mathbb{R}^3 : p\in\Omega \subset \mathbb{S}^2\}.
\]
	The closure $\textup{cl} (\Sigma)$ is a properly embedded surface with boundary
$\partial\Sigma=\zeta^-\cup\zeta^+$ lying on two concentric spheres centered at the origin, where $\zeta^-=\Gamma_M\subset\mathbb{S}^2(1) \equiv \mathbb{S}^2$ and $\zeta^+=(1+M)\Gamma_0\subset\mathbb{S}^2(1+M)$. Here, $\mathbb{S}^2 (\rho ) \subset \mathbb{R}^3$ denotes, as usual, the standard sphere of radius $\rho $ centered at the origin. 

Let $p \in \Sigma $ and let $\{e_1,e_2\}$ be an orthonormal basis of $T_p \Sigma $. A direct computation shows that the unit normal vector to $\Sigma$ at $p$ is given by
$$ N(p) = \frac{(1+M - u(p))}{\sqrt{(1+M - u(p))^2 + |\nabla u|^2}} p + \frac{1}{\sqrt{(1+M - u(p))^2 + |\nabla u|^2}} \sum_{i=1}^2 (-1)^i \langle \nabla u, e_i \rangle e_i \wedge p.$$
In particular, 
$$\langle N(p), p \rangle = \frac{1 + M}{\sqrt{(1 + M)^2 + \alpha_{+,i}^2}} ,$$along each connected component $\zeta^+_i := (1+M) \Gamma_i \subset \zeta ^+ \subset \mathbb{S}^2(1+M)$, where $\alpha _{+,i} \leq 0$ is the constant Neumann condition along each $\Gamma _i  \subset \Gamma  _0$, and 
$$\langle N(p), p \rangle = 1,$$along each connected component $\zeta ^- _i := \gamma_i \subset \Gamma _M \subset \mathbb{S}^2(1)$, which implies that  the surface intersects the spheres $\mathbb{S}^2(1)$ and $ \mathbb{S}^2(1+M)$ at constant angles.

	Set, what we can call the {\it ceiling} and {\it floor} of $\Sigma$, as 
\[
\mathcal C:=\{(1+M)p\in\mathbb{R}^3 : p\in D_i,\ i=1,\ldots,m\} \text{ and } \mathcal F:=\bigcup_{i=m+1}^k D_i,
\]and let $\mathcal O$ denote the (possibly disconnected) open region bounded by
${\rm cl}(\Sigma) \cup \mathcal C \cup \mathcal F$. Observe that $\partial\mathcal O$ is a topological $2-$sphere, possibly with tangential self-intersections at points of $\Max(u)\cap \mathcal F$.

	We start with a simple case to establish the overall strategy, which will clarify the general case.
	
\begin{quote}
{\bf Claim A:}  {\it If $D_{m+1}$ contains the closed geodesic semi-arc $\beta_x^-$, then $(\Omega ,u)$ is rotationally symmetric.}
\end{quote}
	
\begin{figure}[h!]
\centering
\includegraphics[width=0.6\linewidth]{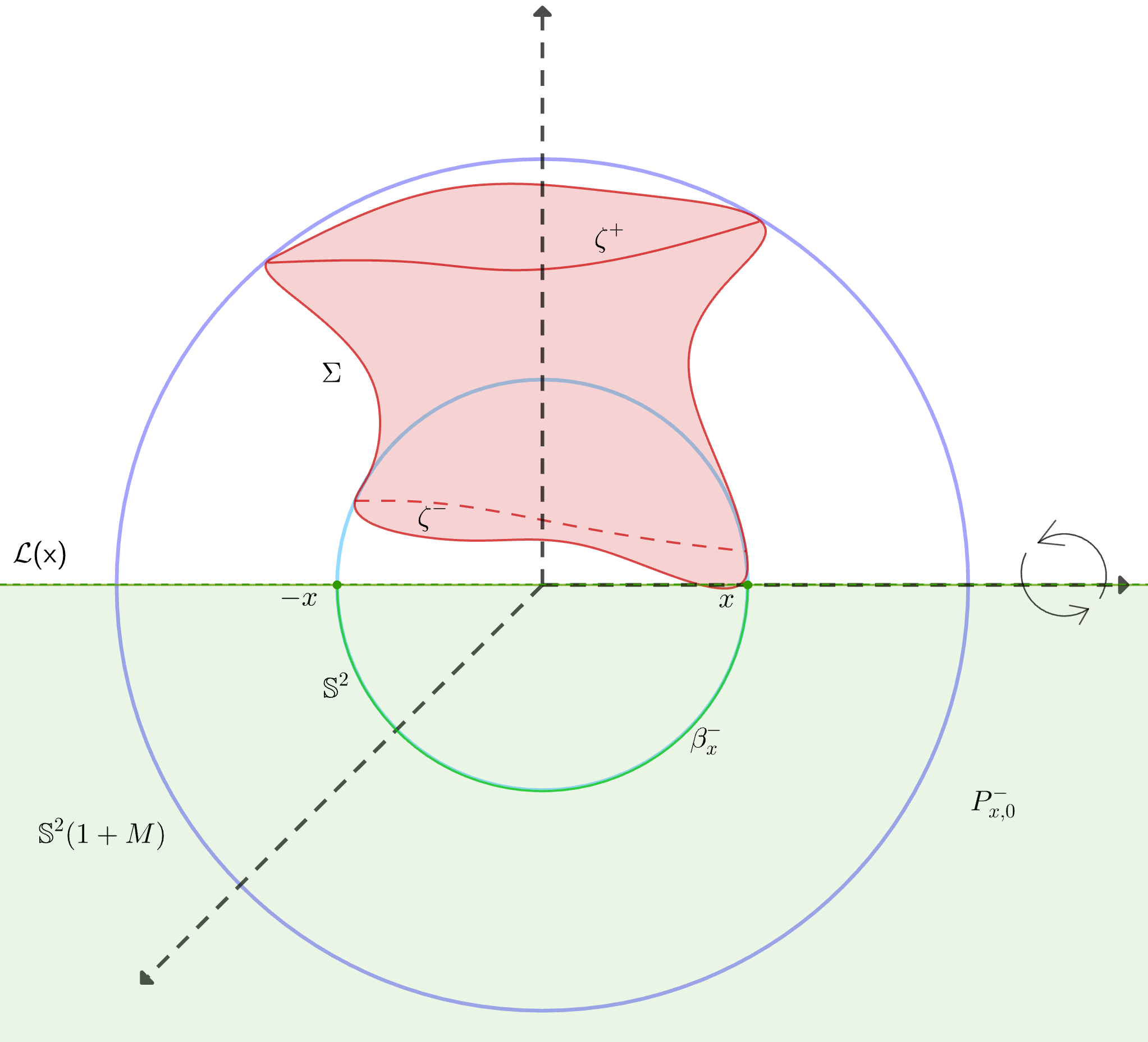}
\caption{In this figure, we present a radial graph $\Sigma$ over a topological annulus $ \Omega  \subset \mathbb{S}^2$. Observe that the curve $\zeta^-$ does not intersect the half-geodesic $\beta_x^-$, which implies that the half plane $P_{x,0}^-$ is disjoint from $\Sigma$. For this reason, we can start the moving plane method in the standard manner.}
\label{fig:RadialGraph}
\end{figure}
	
\begin{proof}[Proof of Claim A]
Construct the foliation $P_{x,\theta}^-$ for $\theta \in [0, 2\pi)$ and define the sector:
$$\mathcal{S}(\theta) = \bigcup_{s \in [0, \theta)} P_{x,s}^-, \text{ for all }\theta \in (0, 2\pi).$$

Since $P^-_{x,0} \cap {\rm cl}(\Sigma) = \emptyset$, by the hypothesis $\beta_x^- \subset D_{m+1}$, there exists $\theta_0>0$ such that $P^-_{x,\theta} \cap  {\rm cl}(\Sigma) = \emptyset$ for all $\theta \in [0,\theta_0)$ and $P^-_{x,\theta_0}$ touches $ {\rm cl}(\Sigma)$ by the first time. Now, apply Alexandrov's reflection method for $\Sigma $ (cf. Figure \ref{fig:RadialGraph}). Thus, for $\theta> \theta_0$ with $\theta-\theta_0$ small, reflecting $\Sigma \cap \mathcal{S}(\theta)$ with respect to $P_{x,\theta}^-$ yields a reflected surface $\widetilde{\Sigma}_\theta := \mathcal{R}^-_{x,\theta}\left( \Sigma \cap \mathcal{S}(\theta) \right)$ that remains within $\mathcal{O}$, since $\partial P_{x,\theta}^- = \mathcal L (x)$ for all $\theta$. Hence, there exists a critical angle $\bar{\theta} \in (0, 2\pi)$ such that:
		
$$\left(\Sigma \setminus \mathcal{S}(\theta)\right) \cap \widetilde{\Sigma}_\theta = \emptyset \text{ and } \widetilde{\Sigma}_\theta \setminus P^-_{x,\theta} \subset \mathcal{O} \text{ for all } \theta \in (\theta_0, \bar{\theta}),$$leading to one of the following cases:
				
\begin{enumerate}
	\item[(a)] $\widetilde{\Sigma}_{\bar{\theta}}$ and $\Sigma \setminus \mathcal{S}(\bar{\theta})$ intersect at an interior point;
	\item[(b)] $\widetilde{\Sigma}_{\bar{\theta}}$ and $\Sigma \setminus \mathcal{S}(\bar{\theta})$ intersect at a boundary point on the same boundary component;
	\item[(c)] $\widetilde{\Sigma}_{\bar{\theta}}$ is orthogonal to $P^-_{x,\bar{\theta}}$ at some point along $\partial \tilde{\Sigma}_{\bar{\theta}}$;
	\item[(d)] a boundary point of $\widetilde{\Sigma}_{\bar{\theta}}$ intersects a point of $(\textup{cl}(\Sigma) \setminus \mathcal{S}(\bar{\theta}) ) \cap \mathbb{S}^2(1)$;
\end{enumerate}

		Define $\widetilde{\Omega}_{\bar{\theta}} = \mathcal{R}^-_{x,\theta}\left( \Omega \cap \mathcal{S}(\bar{\theta}) \right)$ and consider the function  $w_{\bar{\theta}} : \widetilde{\Omega}_{\bar{\theta}} \to \mathbb{R}$ given by $w_{\bar{\theta}}(p) = u\left( \mathcal{R}^-_{x,\theta}(p) \right) - u(p)$. Then, $w_{\bar{\theta}}$ is non-negative in $\widetilde{\Omega}_{\bar{\theta}}$, and there exists a second order linear operator $\mathcal{Q} = \Delta + c(u, \nabla u)$, where $ c $ is a bounded function, such that  $\mathcal{Q} w_{\bar{\theta}} = 0$  in $\widetilde{\Omega}_{\bar{\theta}} $ (see \cite{GT}).
		
		Thus, cases (a), (b) and (c) above are ruled out, respectively, by the interior maximum principle (case (a)), the boundary maximum principle (case (b)), and Serrin's Lemma (case (c)). Specifically, in any of these cases, $P_{x,\bar\theta}$ is a plane of symmetry for $\Omega $ and $u$, i.e., $w_{\bar \theta} \equiv 0 $ in $\textup{cl}( \tilde \Omega^+ _{\bar \theta})$.

	The only situation requiring additional justification is case (d) and why we can apply the boundary maximum principle even when the first contact occurs at a vertex of $\Max(u)$. If case (d) holds, there exists $\bar p \in \left( \Omega \cap \Max (u) \cap \partial \widetilde{\Omega}_{\bar{\theta}} \right) \setminus P^- _{x , \bar \theta}$. Take $\epsilon > 0 $ small enough such that the open disk $D_{\bar p}(\epsilon)$ centered at $\bar p$ of radius $\epsilon$ is contained in $\Omega  $, this follows since each $\gamma _i \subset \Gamma _M$ is an isolated component. Now, shrinking $\epsilon $ if necessary, we can assume that $\mathcal D :=  \widetilde \Omega _{\bar \theta} \cap D_{\bar p}(\epsilon)$ is diffeomorphic to an open half-disk and $\bar p \in \partial \mathcal D \setminus \partial D_{\bar p}(\epsilon)$; recall that $\partial \mathcal D \setminus \partial D_{\bar p}(\epsilon)$ is an analytic arc (no singularities) since each $\gamma _i \subset \Gamma _M$ is embedded and analytic. Consider the function  $w_{\bar{\theta}} : \mathcal D \to \mathbb{R}$, given by $w_{\bar{\theta}}(p) = u\left( \mathcal{R}^-_{x,\theta}(p) \right) - u(p)$, as above, which is well-defined. Then, $\omega _{\bar \theta}$ is non-negative in $\mathcal D$ and $\nabla \omega _{\bar \theta} (\bar p) =0$, which implies that $\omega _{\bar \theta} (\bar p) \equiv 0$ in $\mathcal D$ by the boundary maximum principle. Therefore, $P_{x,\bar\theta}$ is a plane of symmetry for $\Omega $ and $u$ (cf.~Figure~\ref{fig:Immersed}).

\begin{figure}[h!]
\centering
	\includegraphics[width=0.5\linewidth]{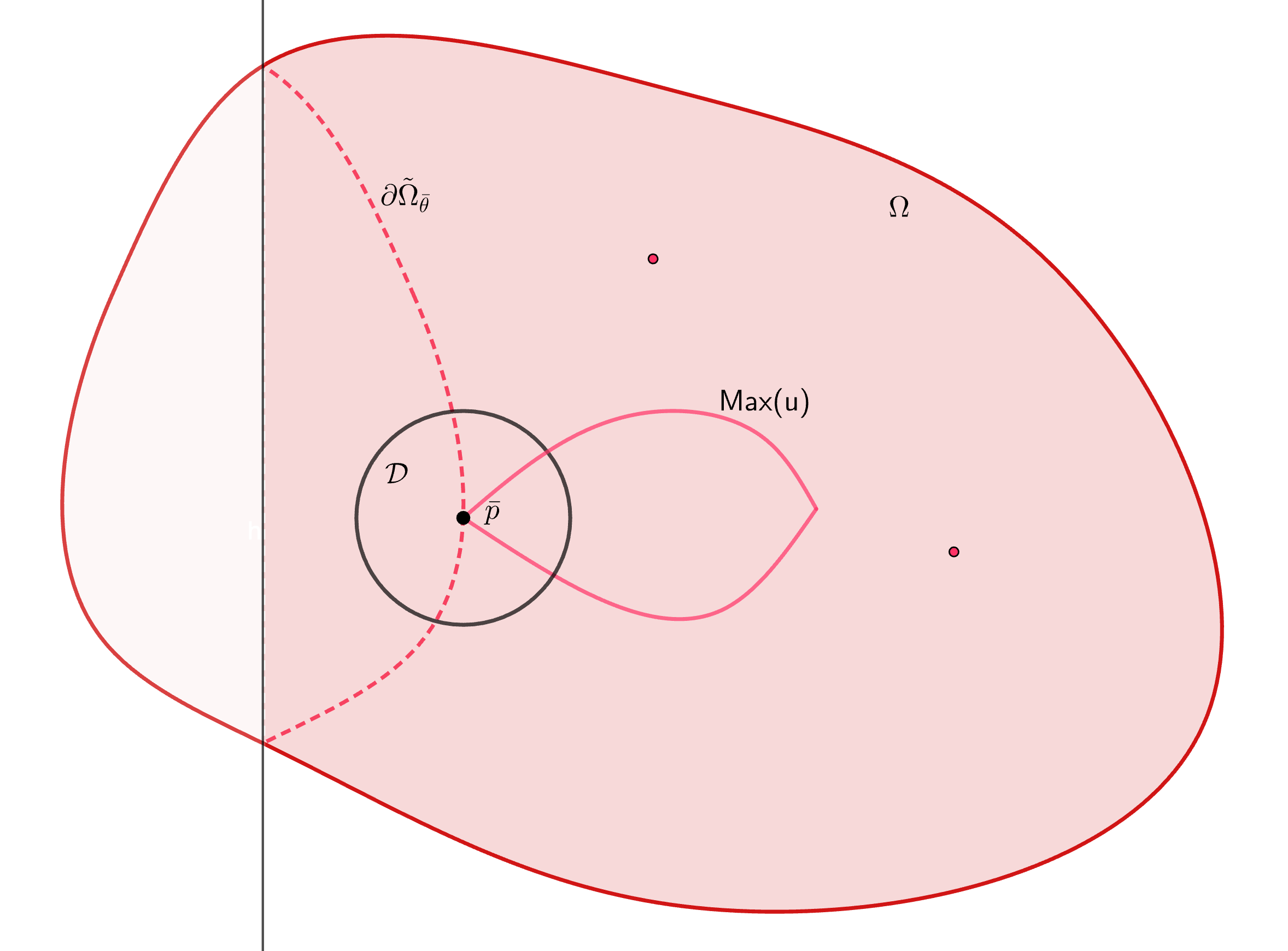}
\caption{Consider the stereographic projection of $\s ^2 \setminus \{ {\bf s}\}$, where ${\bf s}$ is the south pole, and we can suppose without loss of generality that ${\bf s} \in \s ^2 \setminus \textup{cl}(\Omega)$. Here we are assuming that $\partial \Omega = \Gamma_M$ is an isolated component of $\textup{Max}(u)$. In the setting described in case $(d)$, we can apply the boundary maximum principle inside the region $\mathcal{D}$ between $\partial \tilde{\Omega}_{\bar \theta}$ and $D_{\bar p} (\epsilon)$.}
	\label{fig:Immersed}
\end{figure}

		Finally, since $D_{m+1}$ is open, we can find a small neighborhood $U_x$ of $ x$ in $D_{m+1} \subset \s^2$ such that $\beta _y ^- \subset D_{m+1}$ for all $y \in U_x$. Hence, applying the above argument, we can find $\bar \theta (y) \in (0, 2\pi)$ so that $P _{y, \bar \theta (y)}$ is a plane of symmetry of $\Sigma $ for all $y \in U_x$, yielding that $\Sigma$ is rotationally symmetric. This completes the proof of Claim A.
	\end{proof}
	
	Next, consider the case in which $x , -x \in D_{m+1}$ but $\beta ^- _x \cap  \partial D_{m+1}\neq \emptyset$. Note that we can not begin with the moving plane method directly in this case (cf. Figure \ref{fig:Lifting}, left). However, by lifting the surface $\Sigma $ to the universal cover of $\r^3 \setminus \mathcal L  (x)$, we avoid the complications that arise from the presence of the half-geodesic $\beta _x ^-$ and facilitate the application of reflections and isometries, recalling that ${\rm cl}(\Sigma) \cap \mathcal L  (x) = \emptyset$.

\begin{figure}[h!]
\centering
\begin{subfigure}[b]{0.43\linewidth}
\includegraphics[width=\linewidth]{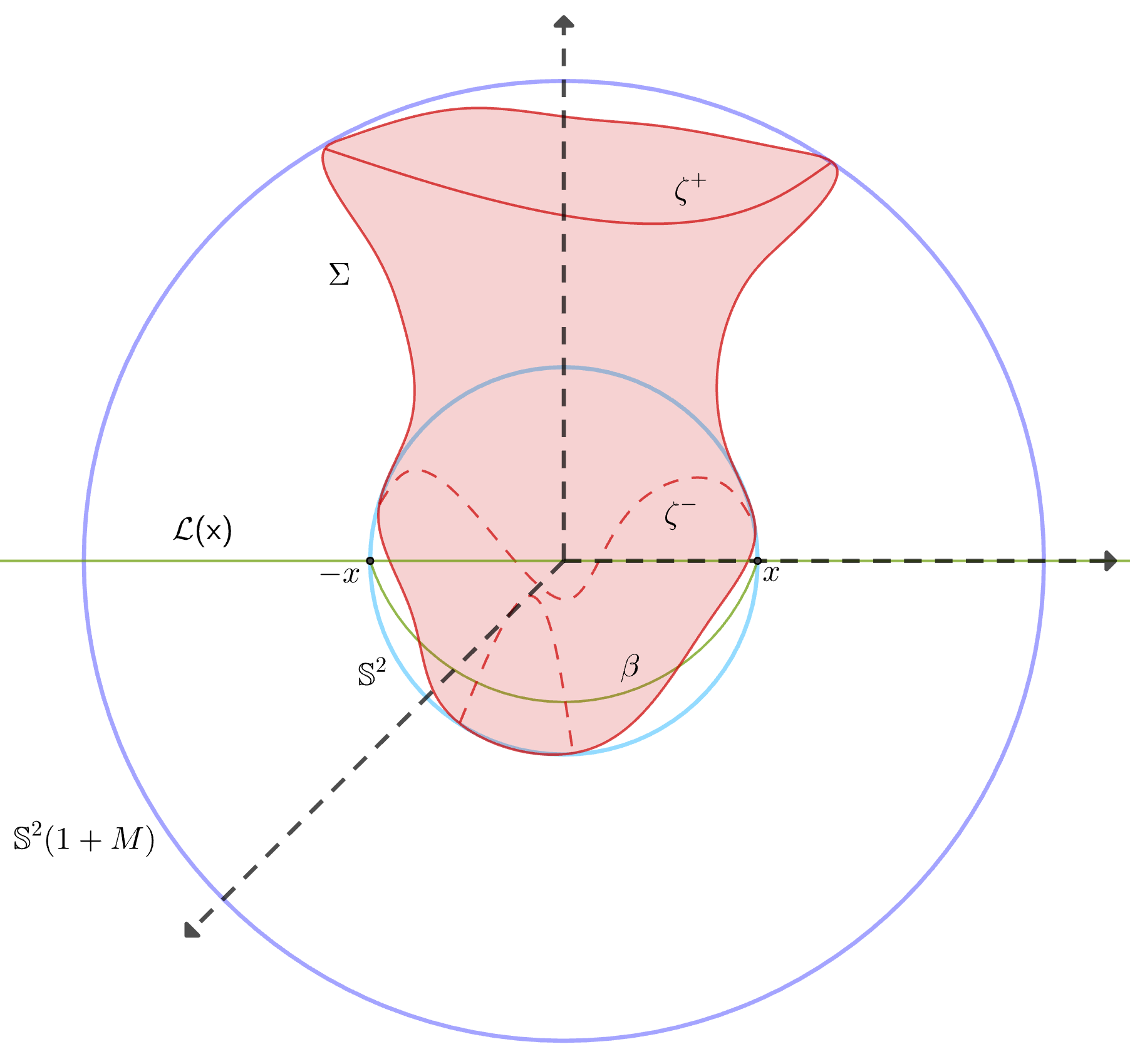}
\end{subfigure}
\begin{subfigure}[b]{0.53\linewidth}
\includegraphics[width=\linewidth]{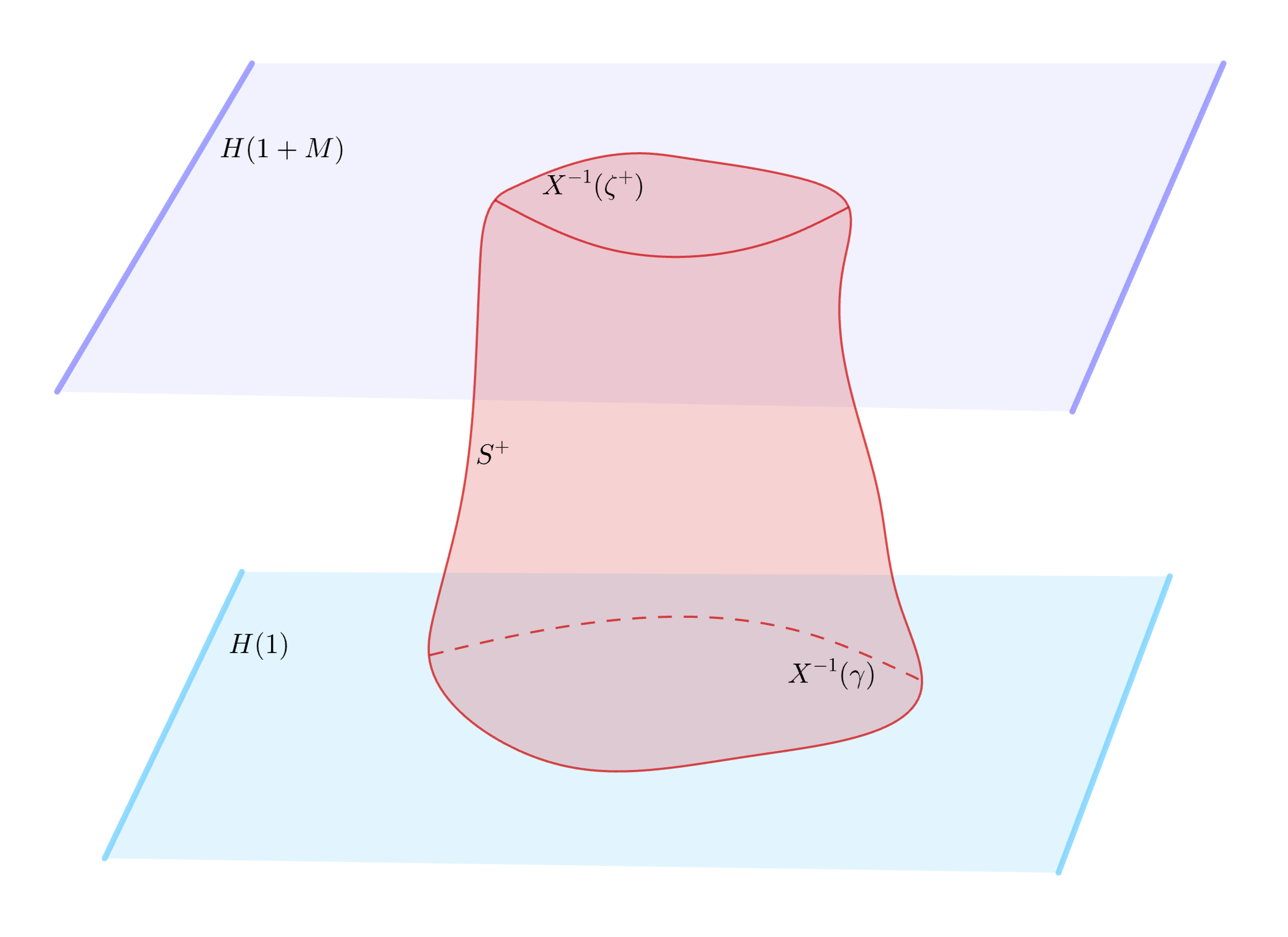}
\end{subfigure}
\caption{{\bf Left:} This figure illustrates the radial graph $\Sigma$ as described at the beginning of this section. Note that the line $\mathcal{L}(x)$ does not intersect $\Sigma$, but any geodesic semi-arc $\beta$ with endpoints in $x$ and $-x$ intersects the curve $\partial D_{m+1} = \zeta^-$ transversely. {\bf Right:} This figures shows the lifting $S=X^{-1}(\Sigma)$ into the parameter space $\mathcal U$. Since $x,-x \in D_{m+1}$, it follows that $\Sigma \subset \r^3 \setminus \mathcal{L}(x)$, where $\mathcal{L}(x)$ denotes the line passing through $x$ and $-x$. This allows a diffeomorphic lifting of $\Sigma$ into $\U$ using the parametrization defined in equation \eqref{lifting}. In the parameter space $\U$, each Euclidean sphere maps to a horizontal band of width $\pi$. }
\label{fig:Lifting} 
\end{figure}
	
\begin{quote} 
{\bf Claim B:} {\it If $D_{m+1}$ contains two antipodal points, $x, -x \in D_{m+1}$, then $(\Omega ,u)$ is rotationally symmetric.}
\end{quote}
\begin{proof}[Proof of Claim B]
Without loss of generality, set $x = (1,0,0)$. Parametrize $\mathbb{R}^3 \setminus \mathcal{L}(x)$ by $\mathcal{U} := (0, \pi) \times \mathbb{R} \times \mathbb{R}_+$, using the covering map $X : \mathcal{U} \to \mathbb{R}^3 \setminus \mathcal{L}(x)$, given by:
		\begin{equation}\label{lifting}
			X(t,s,\rho) = \left(\rho \cos t, \rho \sin t \cos s, -\rho \sin t \sin s \right).
		\end{equation}
		By lifting the Euclidean metric $g_E$ via $X$ on $\mathcal{U}$, the pullback metric becomes:
		$$g := \rho^2 dt^2 + \rho^2 \sin^2 t ds^2 + d\rho^2.$$
		Thus, $X : \left(\mathcal{U}, g \right) \to \left(\mathbb{R}^3 \setminus \mathcal{L}(x), g_E \right)$ is a local isometry. The reflection map in this setting, $\mathcal{R}_{\bar{s}}: (t, s, \rho) \mapsto (t, 2\bar{s} - s, \rho)$, is an isometry of $\left(\mathcal{U}, g \right)$.
		
		Next, consider the surface $\Sigma$ as defined earlier, and note that $\partial \mathcal{O}$ is a contractible topological $2-$sphere in $\mathbb{R}^3 \setminus \mathcal{L}(x)$. Therefore, $X^{-1}(\partial \mathcal{O})$ is also a topological sphere, allowing us to lift $\Sigma$ to $(\mathcal{U}, g)$, giving $S = X^{-1}(\Sigma)$. 
		
		The surface $S$ lies in $\{1 \leq \rho \leq 1 + M\} \subset (\mathcal{U}, g)$, with boundary $\partial S = X^{-1}(\zeta ^-) \cup X^{-1}(\zeta^+)$. Note that $ X^{-1}(\zeta ^-) \subset H(1)$ and $ X^{-1}(\zeta^+) \subset H(1+M)$, where $H(r) := \{\rho = r\}$. Under the map $X^{-1}$, each sphere of radius $\rho \in (0, +\infty)$, $\mathbb{S}^2(\rho) \setminus \{(\pm \rho, 0, 0)\} \subset \mathbb{R}^3 \setminus \mathcal{L}(x)$, is mapped to a horizontal band $H(\rho)$ at height $\rho$ in $(\mathcal{U}, g)$. Since $\Sigma$ forms a constant angle with a sphere along each boundary component, the lifting $S$ intersects $H(1)$ tangentially along $X^{-1}(\zeta ^-)$ and intersects $H(1+M)$ at a constant angle along $X^{-1}(\zeta^+)$ (cf. Figure \ref{fig:Lifting}, right).
		
		For fixed $s \in \mathbb{R}$, the vertical band $Q(s) := (0, \pi) \times \{s\} \times \mathbb{R}^+$ in $\mathcal{U}$ corresponds, under $X$, to the half-plane $P_{s \, (\text{mod} \, 2\pi), x} := X(Q(s))$, which is the half-plane with boundary $\mathcal{L}(x)$ passing through $(0, \cos s, -\sin s) \in \mathbb{R}^3$. 
		
		In this setting, we apply the moving plane method to $S$ with respect to the foliation of vertical bands $\{Q(s)\}_{s \in \mathbb{R}}$ in $(\mathcal{U}, g)$. Define $Q^-(s) = \bigcup_{-\infty}^s Q(s)$ and $\widetilde{S}_s = \mathcal{R}_s(S \cap Q^-(s))$ (cf. Figure \ref{fig:Liftin2}, left). 
		
		\begin{figure}[h!]
			\centering
			\begin{subfigure}[b]{0.47\linewidth}
				\includegraphics[width=\linewidth]{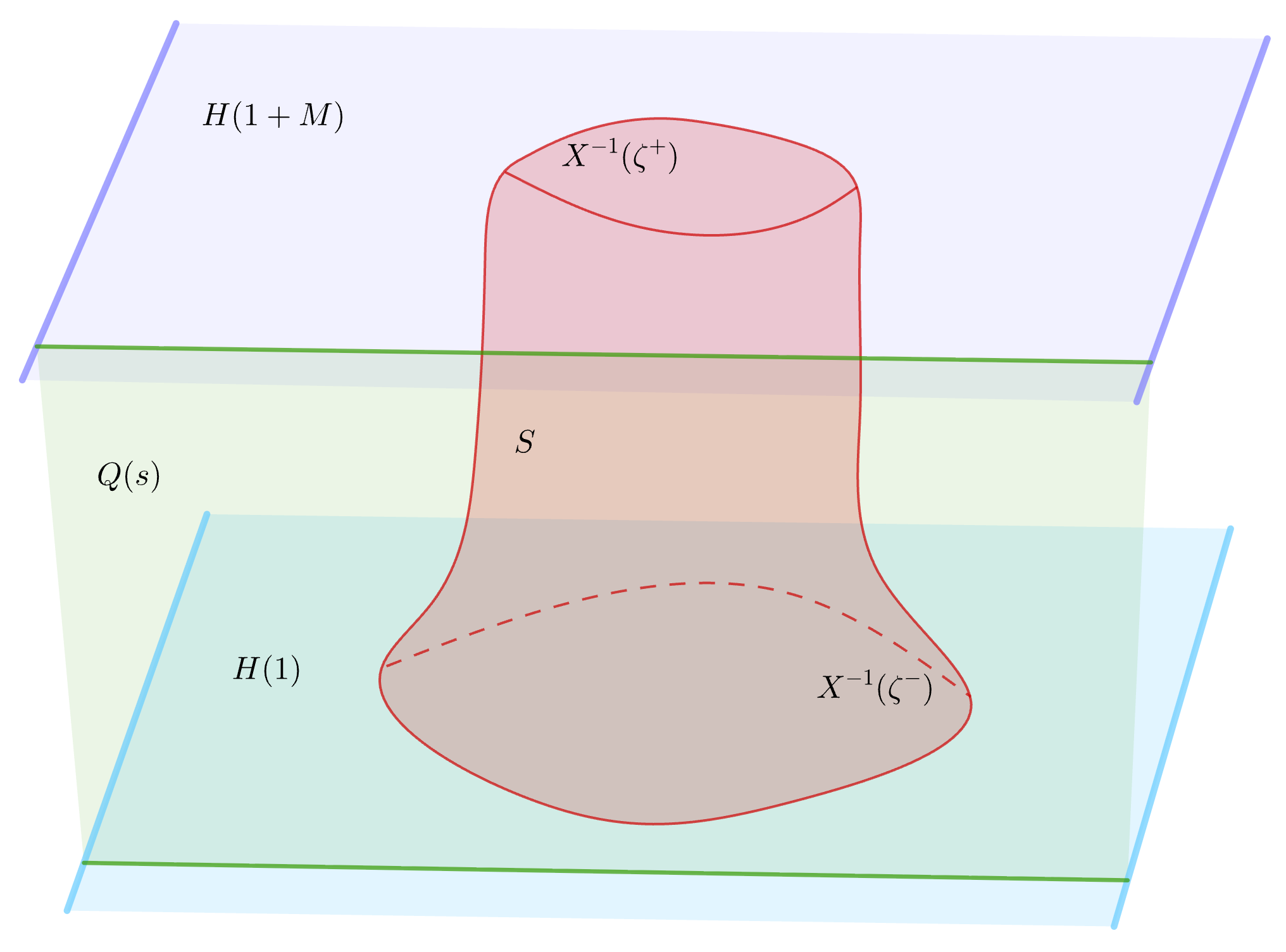}
			\end{subfigure}
			\begin{subfigure}[b]{0.477\linewidth}
				\includegraphics[width=\linewidth]{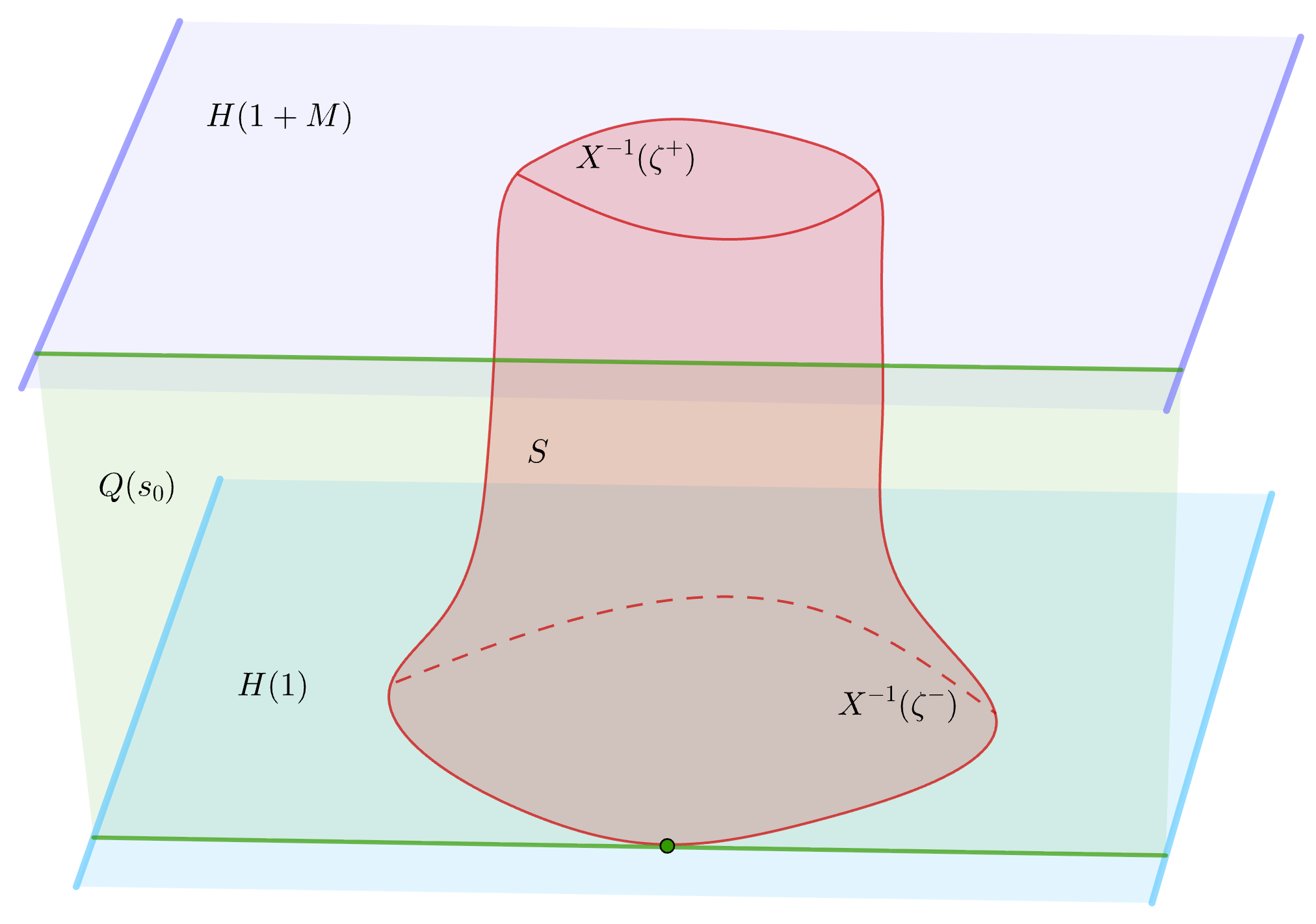}
			\end{subfigure}
			\caption{This figure illustrates how the moving plane method can be applied to $S$ within the domain $ \mathcal{U} $. {\bf Left:} Since $ S $ is compact, there exists a real number $ s \in \mathbb{R} $ such that the plane $ Q(t) $ does not intersect $ S $ for all $t \leq s$; that is, $ Q(t) \cap S = \emptyset $ for all $t \leq s$. {\bf Right:} As we increase $s$, moving $ Q(s) $ towards $ S $, there exists a smallest value $ s_0 > s $ where $ Q(s_0) $ first makes contact with  $S $. This initial point of contact at $ s_0 $ allows us to begin the reflection process, initiating the moving plane method from this position.}
			\label{fig:Liftin2}
		\end{figure}
		
		Then there exists $s_0 \in \r $ such that $Q(s_0)$ touches $S$ for the first time (cf. Figure \ref{fig:Liftin2}, right) and, by applying the standard moving plane method, there exists $\bar{s} \in \mathbb{R}$, $s_0 < \bar s$, such that:
		$$\left(S \setminus Q^-(s)\right) \cap \widetilde{S}_s = \emptyset \text{ and } \widetilde{S}_s \setminus Q(s) \subset X^{-1}(\mathcal{O}) \text{ for all } s \in (s_0, \bar s),$$
		and one of the following occurs:
		\begin{enumerate}
			\item[(a)] $\widetilde{S}_{\bar{s}}$ and $S \setminus Q^-(\bar{s})$ intersect at an interior point;
			\item[(b)] $\widetilde{S}_{\bar{s}}$ and $S \setminus Q^-(\bar{s})$ intersect at a boundary point on the same boundary component;
			\item[(c)] $\widetilde{S}_{\bar{s}}$ is orthogonal to $Q(\bar{s})$ at some point along $\partial \tilde{S}_{\bar{s}}$; 
			\item[(d)] a boundary point of $\widetilde{S}_{\bar{s}}$ intersects an interior point of $(\textup{cl}(S) \setminus Q^-(\bar{s}) ) \cap H(1)$.
		\end{enumerate}
		
		Then, by pushing forward a neighborhood in any of the above cases, the function $w_{\bar \theta} : \widetilde \Omega _{\bar \theta} \to \r $, $\bar \theta = \bar s \, ({\rm mod} \, 2 \pi)$, given by $w _{\bar \theta} (p) = u (\mathcal R_{x, \bar \theta} (p)) -u (p)$ is non-negative and satisfies a second order linear elliptic equation in a neighborhood of any point in the above cases. Thus, we can proceed as in Claim A to rule out all the case and prove radial symmetry. This proves Claim B.
	\end{proof}
Therefore, Claim B implies that $(\Omega ,u)$ is rotationally symmetric and completes the proof of the lemma in case $(i)$.

To prove the lemma in case $(ii)$, we follow exactly the same steps as in the proof of case $(i)$, but considering the radial graph
\[
\Sigma := \{(1+u(p))\,p \in \mathbb{R}^3 : p \in \Omega \subset \mathbb{S}^2\},
\]
where we now define $\partial \Sigma = \zeta^- \cup \zeta^+$ with $\zeta^- = \Gamma_0 \subset  \mathbb{S}^2$ and $\zeta^+ = (1+M)\Gamma_M \subset \mathbb{S}^2(1+M)$. The only difference with respect to case $(i)$ is that, when applying the Alexandrov reflection method to the lifted surface $S = X^{-1}(\Sigma)$ as in the proof of Claim~B, case $(d)$ may occur only if there is an intersection between $\partial \tilde{S}_{\bar s}$ and a different boundary component of $\zeta^-$. Since $|\alpha _j | \leq |\alpha _i|$ for all $j \in \{1, \ldots , m \}$ and $j \neq i$, in this situation we can conclude symmetry by applying the boundary maximum principle in a standard manner.
\end{proof}

Observe that, when we are in case $(ii)$ with $\Gamma_M = \emptyset$, Lemma~\ref{LemSirakov} recovers the theorem of Kumaresan--Prajapat~\cite{KP} under the weaker assumption that a connected component of the complement of the closure of the domain contains two antipodal points. We are now ready to prove:
\\
\\
{\bf Theorem A:} {\it Let $(\Omega, u)$, $\Omega \subset \s^2$ a finite domain, be an analytic solution to \eqref{OEP}. Assume that $\Max(u)$ contains a simple closed curve $\gamma$ that is an isolated component of $\Max(u)$. Then, $(\Omega ,u)$ is either antipodally or rotationally symmetric. In the latter case, $\Omega$ is either a disk or an annulus.}

\begin{proof}[Proof of Theorem A]
First, if $\gamma$ is antipodally symmetric, Lemma \ref{lemmaTheoremA1} implies that $(\Omega,u)$ is antipodally symmetric. Henceforth, we assume that $\Max(u)$ contains a simple closed curve $\gamma$ which is an isolated component of $\Max(u)$ and is not antipodally symmetric.

Since $\gamma$ is an analytic simple closed curve in $\mathbb{S}^2$, it separates the sphere into two topological disks, $\mathbb{S}^2 \setminus \gamma = D^+ \cup D^-$. As $\gamma$ is not antipodally symmetric, there exists $\bar p \in \gamma$ such that $\mathcal{A}(\bar p)=-\bar p \in D^+ \cup D^-$. By a small perturbation, we may choose a nearby point $\bar q$ so that $\bar q$ and $\mathcal{A}(\bar q)$ both lie in the same component of $\mathbb{S}^2\setminus \gamma$; in particular, we can ensure that there exist antipodal points $x,-x$ contained in one of the two disks, say $x,-x \in D^-$.

We now distinguish whether $\gamma$ is contractible in $\Omega$ or not, but we treat both situations in the same way: in each case, we identify a subdomain of $\Omega$ whose complement contains the antipodal pair $x,-x$, and then apply Lemma~\ref{LemSirakov}.

If $\gamma$ is contractible in $\Omega$, then it bounds a topological disk $D\subset\Omega$ with $\partial D=\gamma$, and the remainder $\Omega_\gamma:=\Omega\setminus {\rm cl}(D)$ is connected with $\partial \Omega_\gamma=\gamma\cup \partial\Omega$. Since $x$ and $-x$ lie in the same component of $\mathbb{S}^2\setminus\gamma$, they are contained either in $D$ or in $\mathbb{S}^2\setminus \overline{D}$. In the first case, we apply case $(i)$ of Lemma~\ref{LemSirakov} to the domain $\Omega_{\gamma}$. In the second case, we apply case $(i)$ of Lemma~\ref{LemSirakov} (with $\Gamma_0 = \emptyset$) to the disk $D$. In both instances the hypotheses of Lemma~\ref{LemSirakov} are satisfied, and we conclude that $(\Omega,u)$ is rotationally symmetric.

If $\gamma$ is non-contractible in $\Omega$, then $\Omega\setminus \gamma$ has two connected components, say $\Omega_\gamma^+$ and $\Omega_\gamma^-$. Since $x,-x$ lie in the same component of $\mathbb{S}^2\setminus \gamma$, we may assume (up to relabeling) that $x,-x$ belong to a connected component of $\mathbb{S}^2\setminus {\rm cl}(\Omega_\gamma^+)$. Therefore, case $(i)$ of Lemma~\ref{LemSirakov} applies to $\Omega_\gamma^+$ and yields again that $(\Omega,u)$ is rotationally symmetric.

This finishes the proof of Theorem A.
\end{proof}

If, in addition, $f(x,0)>0$ for all $x>0$, then Theorem A can be refined as follows:
\\
\\
{\bf Corollary A.} {\it Let $(\Omega, u)$, $\Omega \subset \s^2$, be an analytic solution to \eqref{OEP} and assume that $f(x,0)>0$ for all $x>0$. Then, either $\textup{Max}(u)$ contains only finitely many isolated points or the solution $(\Omega ,u)$ is either antipodally or rotationally symmetric. In the latter case, $\Omega$ is an annulus.}

\begin{proof}[Proof of Corollary A]
First, by analyticity and the positivity of $f(\cdot,0)$, the set $\textup{Max}(u)$ consists either of finitely many isolated points or contains at least one closed curve $\gamma\subset\textup{Max}(u)$ (see \cite[Corollary~3.4]{Chr}). Moreover, $\textup{Max}(u)$ cannot contain transversal, tangential, or cuspidal intersections; hence each connected component of $\textup{Max}(u)$ that is not an isolated point is simple, closed, and isolated from the others. Therefore, Theorem~A implies that $(\Omega,u)$ is either antipodally symmetric or rotationally symmetric.

In the latter case, $\Omega$ is a disk or an annulus and $\textup{Max}(u)$ contains an isolated analytic closed curve. Thus, $\Omega$ cannot be a disk in this situation. Indeed, otherwise there would exist a curve $\gamma\subset\textup{Max}(u)$ bounding a contractible domain $\mathcal{D}\subset \Omega\setminus\textup{Max}(u)$. Then $u$ would attain a minimum at some interior point $p \in \mathcal{D}$. But in such a point, it would be $\Delta u (p) = -f(u(p),0) <0$, a contradiction. Hence $\Omega$ must be an annulus, which proves the result.

\end{proof}

\section{Geometric applications}\label{SectGA}

The goal of this section is to translate the analytic information provided by Theorem~A into a geometric rigidity result for constant mean curvature surfaces.

Let $\mathbb{B}^3(r)$ denote the open Euclidean ball of radius $r>0$ centered at the origin,
with boundary $\partial\mathbb{B}^3(r)=\mathbb{S}^2(r)$.
For simplicity, we write $\mathbb{B}^3=\mathbb{B}^3(1)$ and
$\mathbb{S}^2=\mathbb{S}^2(1)$ for the unit ball and the unit sphere, respectively.

Let $\Sigma\subset\mathbb{R}^3$ be a properly immersed surface with boundary, and assume that
each boundary component $\zeta_i\subset\partial\Sigma$ is contained in a sphere
$\mathbb{S}^2(r_i)$ for some $r_i>0$.
Denote by $\nu:\partial\Sigma\to\mathbb{S}^2$ the outer conormal to~$\Sigma$.
We assume that each boundary component intersects its supporting sphere from the inside, that is,
\[
\langle \nu(p),p\rangle \ge 0, \qquad \forall\, p\in\partial\Sigma.
\]

\begin{definition}\label{surfacesWithFreeBoundaries}
Let $\Sigma\subset\mathbb{R}^3$ be a properly immersed surface with boundary.
We say that $\Sigma$ has \emph{capillary boundaries} if each boundary component
$\zeta_i\subset\partial\Sigma$ meets a sphere $\mathbb{S}^2(r_i)$ centered at the origin
at a constant contact angle $\theta_i\in [0,\pi)$, where the radii and angles may vary from one
component to another.

In particular, if $\Sigma\subset\mathbb{B}^3$, $\partial\Sigma\subset\mathbb{S}^2$, and all
boundary components meet $\mathbb{S}^2$ at the same constant angle, then $\Sigma$ is called a \emph{capillary surface} in the unit ball.
If this angle equals $\pi/2$, then $\Sigma$ is said to be a \emph{free boundary surface}
in~$\mathbb{B}^3$.
\end{definition}

Let $I$ and $II$ denote the first and second fundamental forms of $\Sigma$.
The mean curvature $H$ and the extrinsic curvature $K_e$ are defined by
$2H=k_1+k_2$ and $K_e=k_1k_2$, where $k_1$ and $k_2$ are the principal curvatures.
We say that $\Sigma$ is a CMC surface if its mean curvature is constant.
Throughout the paper, the mean curvature is computed with respect to a choice
of unit normal $N$ such that $H>0$.

From now on, let $\Sigma\subset\mathbb{R}^3 \setminus \set{\bf 0}$ be an analytic properly embedded CMC annulus
with capillary boundaries $\partial\Sigma=\zeta^-\cup\zeta^+$, where
$\zeta^\pm\subset\mathbb{S}^2(r_\pm)$ for some $0<r_-<r_+$. Set $0< r_0 := \min\{ |p| : p\in\Sigma \} < r_- $. We recall that a neck of $\Sigma$ is a closed (analytic) curve $\gamma \subset \Sigma$ along which the distance function to the origin on $\Sigma$, denoted by $d_\Sigma$, attains its absolute minimum. Now we are ready to show:
\\
\\
{\bf Theorem B:} {\it Let $\Sigma \subset \mathbb{R}^3 \setminus \set{\bf 0} $ be a properly embedded (analytic) CMC annulus with capillary boundaries, $\partial \Sigma= \zeta ^- \cup \zeta ^+$, such that $\zeta ^- \subset \mathbb{S}^2 (r_-)$ and $\zeta ^+ \subset \mathbb{S}^2 ( r_+)$ for some $0<  r_- \leq r_+$. Set $0< r_0 := {\rm min}_{\Sigma} \, {\bf d}_\Sigma < r_-$ and suppose that 
\begin{enumerate}
\item[(i)] $\Sigma$ contains a neck $\gamma\subset\Sigma\cap\s^2(r_0)$ which is an isolated component of ${\rm Min}({\bf d}_\Sigma)$.

\item[(ii)] Let $\Sigma ^\pm \subset \Sigma \setminus \gamma$ be the connected component such that $\partial \Sigma ^\pm = \zeta ^\pm \cup \gamma$, $i=1,2$. The cone over $\gamma$, 
$$C(\gamma):=\set{ \lambda \cdot p \in \mathbb{R}^3 \, : \, \lambda > 0 \text{ and } p \in \gamma},$$ separates $\Sigma ^-$ and $\Sigma ^+$. 

\item[(iii)]  ${\rm int}(\Sigma ^\pm) \subset \textup{cl}(\b ^3 (r_\pm)) \setminus \b ^3 (r_0)$.	
\end{enumerate}
Then, $\Sigma$ is either antipodally or rotationally symmetric. In the latter case, $\Sigma$ is a part of either a catenoid, if $H=0$, or a Delaunay surface, if $H\neq 0$. Moreover, when $\Sigma$ is antipodally symmetric, then $r_-= r_+$.}

\begin{proof}[Proof of Theorem B]

By condition ${\rm (i)}$, $\Sigma$ contains a simple neck $\gamma$ which is an isolated component of ${\rm Min}({\bf d}_\Sigma)$. By definition of a neck, $\Sigma$ is tangent to the sphere $\mathbb{S}^2(r_0)$ along $\gamma$,
where $r_0$ is the global minimum of the distance function to the origin.
In particular, $\gamma$ is a line of curvature of $\Sigma$ by Joachimsthal's Theorem.

We claim that $\gamma$ is non--contractible in $\Sigma$.
Suppose, by contradiction, that $\gamma$ bounds a topological disk $D\subset\Sigma$.
Since $\gamma$ is a line of curvature, Nitsche's Theorem \cite{JNit} implies that
$\Sigma$ must be totally umbilical.
As $\Sigma$ and $\mathbb{S}^2(r_0)$ are tangent along $\gamma$, this would force
$\Sigma$ to be contained in $\mathbb{S}^2(r_0)$, which is impossible.
Hence, $\gamma$ does not bound a disk in $\Sigma$, and therefore it is non--contractible.

We now distinguish according to whether $\gamma$ is antipodally symmetric.
Assume first that $\gamma$ is antipodally symmetric, that is, $\gamma=-\gamma$.
Let $\mathcal A:\mathbb{R}^3\to\mathbb{R}^3$ denote the antipodal map and set
$\tilde\Sigma:=\mathcal A(\Sigma)$.
Then $\Sigma$ and $\tilde\Sigma$ coincide along $\gamma$ and are both tangent to
$\mathbb{S}^2(r_0)$ along this curve.
In particular, their unit normals agree along $\gamma$.
By the uniqueness of the Bj\"orling problem for analytic surfaces (see \cite{Brander}),
it follows that $\mathcal A(\Sigma)\equiv\Sigma$, and hence $\Sigma$ is antipodally symmetric. Since $\zeta^\pm\subset\mathbb{S}^2(r_\pm)$ and the antipodal map preserves spheres,
we have $\mathcal A(\zeta^\pm)\subset\mathbb{S}^2(r_\pm)$.
As $\Sigma$ is embedded, the antipodal symmetry exchanges the boundary components,
that is, $\mathcal A(\zeta^-)=\zeta^+$.
Consequently, $r_-=r_+$.

Assume next that $\gamma$ is not antipodally symmetric.
Since $\gamma$ is non--contractible in $\Sigma$, it separates $\Sigma$ into two connected components,
\[
\Sigma\setminus\gamma=\Sigma^-\cup\Sigma^+,
\qquad
\partial\Sigma^\pm=\zeta^\pm\cup\gamma .
\]
We apply the Alexandrov reflection method to $\Sigma^+$, in the same spirit as in the proof of Lemma~\ref{LemSirakov}.

By assumption {\it (ii)}, the cone
\[
C(\gamma):=\{\lambda p\in\mathbb{R}^3:\lambda>0,\ p\in\gamma\}
\]
separates $\Sigma^-$ and $\Sigma^+$.
Since $\gamma$ is not antipodally symmetric, there exists a direction $x\in\mathbb{S}^2$
such that the line $\mathcal L(x)$ is disjoint, say, from $\textup{cl}(\Sigma^+)$.

To implement the reflection method, we lift the problem to the universal cover of
$\mathbb{R}^3\setminus \mathcal L(x)$.
Let $X:\mathcal U\to\mathbb{R}^3\setminus \mathcal L(x)$ be the covering map defined in \eqref{lifting},
where $\mathcal U=(0,\pi)\times\mathbb{R}\times\mathbb{R}_+$, and set $S:=X^{-1}(\Sigma^+)$.
Then $S$ is a properly embedded surface in $\mathcal U$.

By assumption {\it (iii)}, ${\rm int}(\Sigma^+)\subset \textup{cl}(\mathbb{B}^3 (r_+))\setminus\mathbb{B}^3(r_0)$,
and therefore $S\subset\{r_0\le\rho\le r_+\}\subset\mathcal U$.
Moreover,
\[
X^{-1}(\gamma)\subset H(r_0):=\{\rho=r_0\},
\qquad
X^{-1}(\zeta^+)\subset H(r_+):=\{\rho=r_+\}.
\]

In $\mathcal U$, consider the family of vertical bands $Q(s):=(0,\pi)\times\{s\}\times\mathbb{R}_+$.
Each $Q(s)$ corresponds, under $X$, to a half--plane in $\mathbb{R}^3$ containing the line $\mathcal L(x)$.
For $s\ll0$, $Q(s)$ is disjoint from $S$, and we start the Alexandrov reflection procedure.

As $s$ increases, $Q(s)$ moves towards $S$ until a first contact occurs at some critical value $\bar s$.
Denote by $\tilde S_{\bar s}$ the reflection of $S\cap Q^-(\bar s)$ across $Q(\bar s)$, where
$Q^-(\bar s)=\{(t,s,\rho)\in\mathcal U:s\le\bar s\}$.
At this first contact, one of the following situations arises:
\begin{enumerate}
\item[(a)] an interior contact point;
\item[(b)] a boundary contact point on the same boundary component;
\item[(c)] an orthogonal contact with $Q(\bar s)$;
\item[(d)] a boundary point of $\tilde S_{\bar s}$ touching an interior point
of $S\setminus Q^-(\bar s)$ along $H(r_0)$.
\end{enumerate}
In each case, the relevant maximum principle (interior, boundary, or Serrin's lemma) forces
$S$ to be symmetric with respect to $Q(\bar s)$. Finally, there exists an open neighborhood $U_x\subset\mathbb{S}^2$ of $x$ such that
$\mathcal L(y)\cap \textup{cl}(\Sigma^+)=\emptyset$ for all $y\in U_x$.
Repeating the reflection argument for all such directions, we obtain an open family of planes of symmetry. Therefore, $\Sigma$ is therefore rotationally symmetric.

This completes the proof of Theorem~B.
\end{proof}

\begin{figure}[!h]
	\centering
	\begin{subfigure}[b]{0.445\linewidth}
		\includegraphics[width=\linewidth]{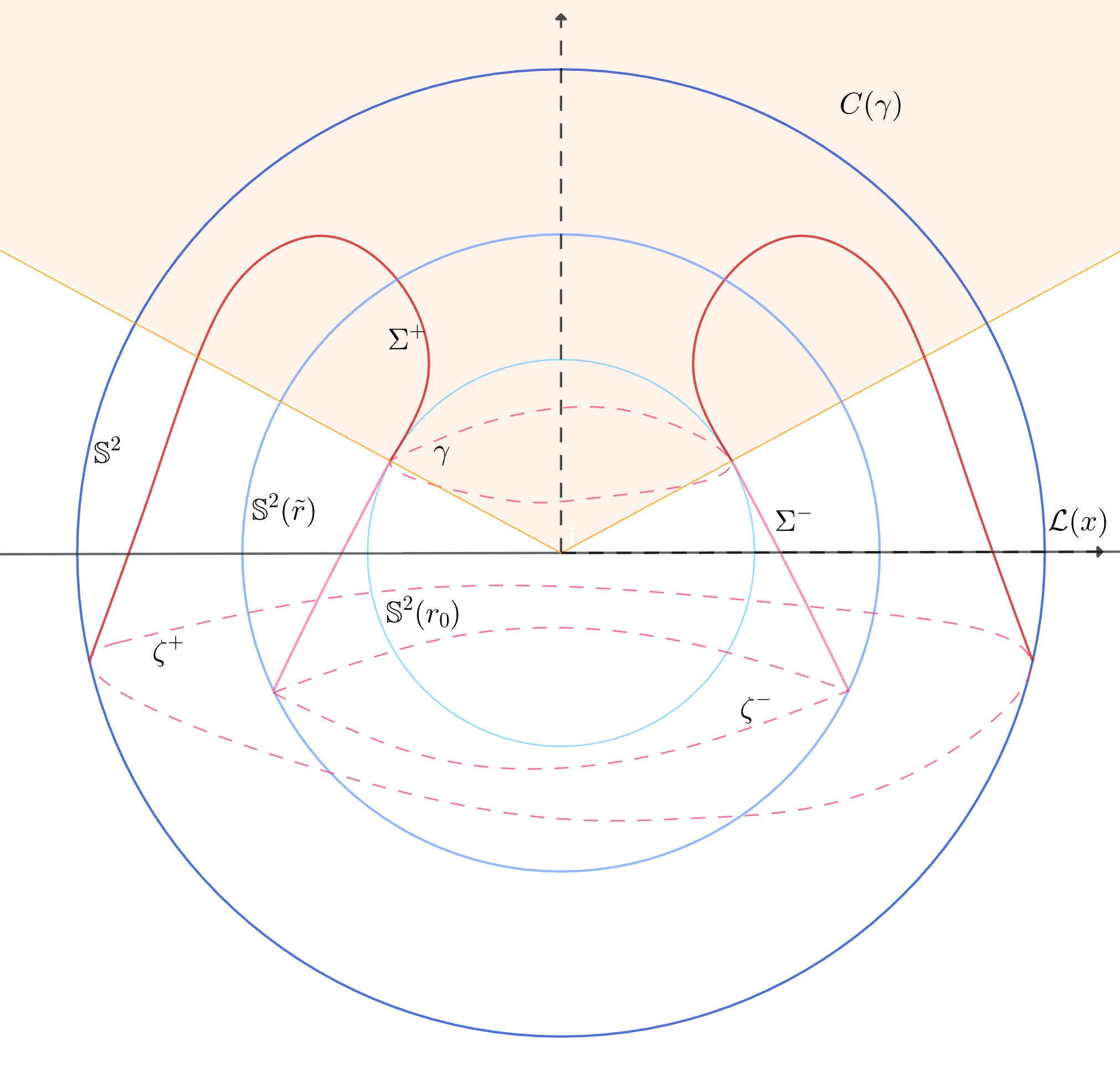}
		\label{fig:AnnulusNonContainedInCone}
	\end{subfigure}
	\begin{subfigure}[b]{0.43\linewidth}
		\includegraphics[width=\linewidth]{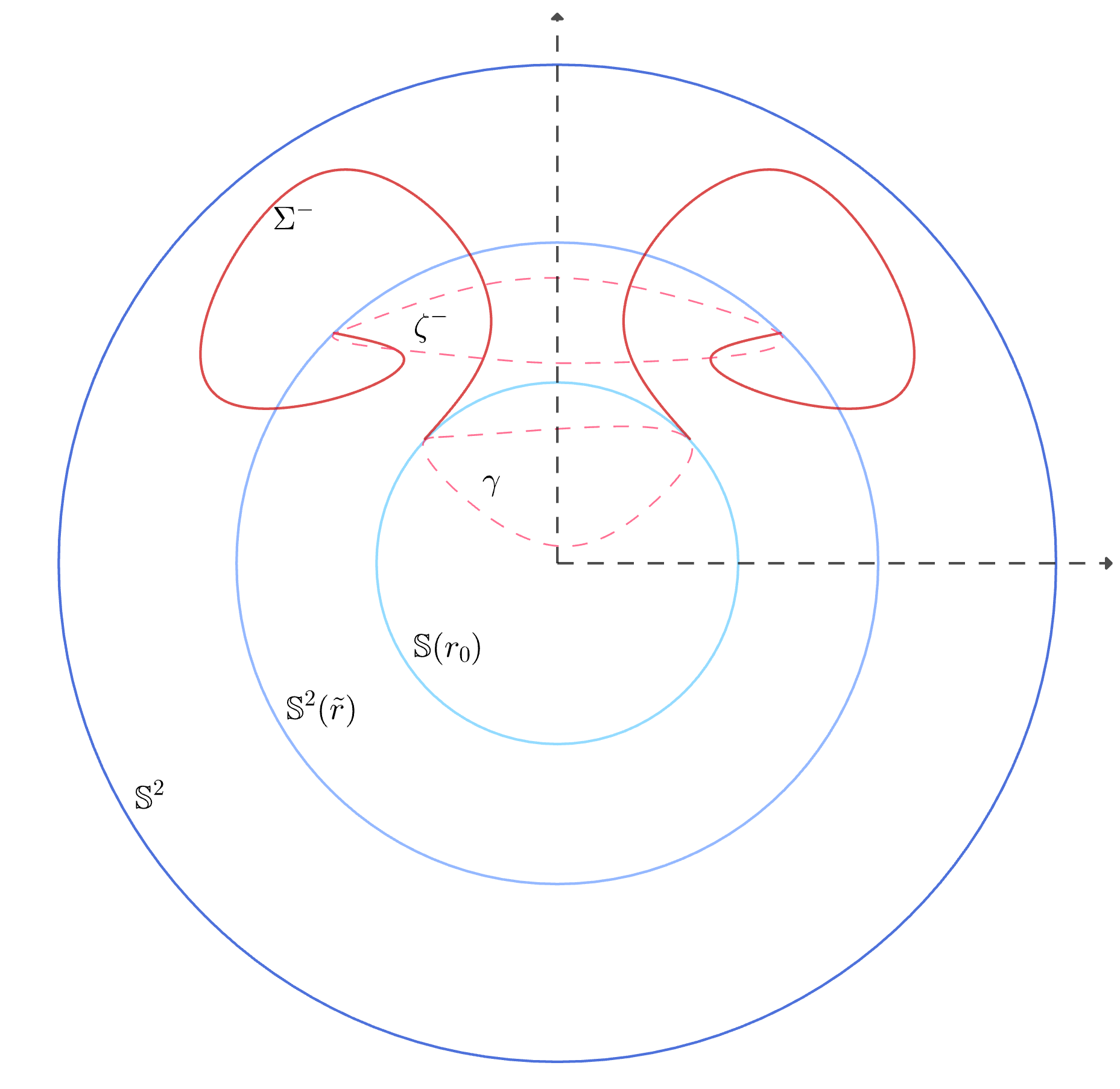}
		\label{fig:NonProperlyOrientedAnnnulus}
	\end{subfigure}
	\caption{This figure presents two vertical cross-sections of surfaces with capillary boundaries, illustrating scenarios where the conditions of Theorem B fail. {\bf Left:} we depict an annulus with capillary boundaries where condition {\it (ii)} is not satisfied; specifically, the cone over $ \gamma $, denoted  $C(\gamma) $, does not separate the connected components of $ \Sigma \setminus \gamma $. In this case, the reflection procedure from Theorem A cannot be initiated because we cannot find a direction $ x \in \mathbb{S}^2$ such that the line $ \mathcal{L}(x) $ is disjoint from one of the connected closed components of $ \Sigma \setminus \gamma $. {\bf Right:} we show a cross-section of the component $\Sigma^-$ of an annulus with capillary boundaries where condition {\it (iii)} is not satisfied. Here, although the method can be applied, an intersection between an interior point and a boundary point may occur at a location where the maximum principle does not apply.} 
	\label{fig:CounterExamplesCapillaryBoundaries}
\end{figure}

\subsection{Remarks on Theorem~B}

We collect here several remarks clarifying the role of assumptions {\it (i)--(iii)}
in Theorem~B and emphasizing the differences between the minimal and the CMC cases.

\paragraph{Remarks on the assumptions of Theorem~B.}

\begin{itemize}
\item
Assumption {\it (i)} cannot be removed in the constant mean curvature setting.
Indeed, there exist properly embedded CMC annuli with capillary boundaries
that do not contain a neck and are not rotationally symmetric
(see, for instance, the examples constructed in \cite{CMF,FHM}).
By contrast, in the minimal case the existence of a neck is expected to follow
from the Critical Catenoid Conjecture and therefore should not be regarded as
an additional restriction.

\item
Assumptions {\it (ii)} and {\it (iii)} encode geometric separation and orientation
conditions that are specific to the CMC case.
In general, they cannot be dropped, as illustrated by the counterexamples
depicted in Figure~\ref{fig:CounterExamplesCapillaryBoundaries}.

\item 
It should be noted that Theorem~B remains valid even when 
$\partial \Sigma^\pm \cap \s^2(r_\pm)$ has more than one connected component.
More generally, the surface $\Sigma$ may have a more complicated topology.
The conclusion of Theorem~B still holds provided that the neck $\gamma$
separates $\Sigma$ into two connected components, $\Sigma \setminus \gamma = \Sigma^+ \cup \Sigma^-$, and each of them is contained in a different connected component of $\mathbb{R}^3 \setminus C(\gamma) = \mathcal C^+ \cup \mathcal C^-$. Since each region $\mathcal C^\pm$ is simply connected, we may lift the
component, say $\mathcal C^+$, that does not intersect the line
$\mathcal L(x)$ to the parameter space, allowing the Alexandrov reflection
argument to proceed exactly as before.
\end{itemize}

\paragraph{Comparison with the minimal and ESW cases.}

\begin{itemize}
\item
In the minimal case ($H=0$), assumption {\it (ii)} is unnecessary when $\Sigma$ has free boundaries.
Indeed, embedded minimal annuli with free boundaries are radial graphs, as shown
in \cite{EMa}.
Moreover, these techniques extend to ESW--surfaces of minimal type, implying
that embedded ESW--annuli with free boundaries are radial graphs as well.

\item
Assumption {\it (iii)} is likewise superfluous in the minimal case and, more
generally, for ESW--surfaces of minimal type, since these surfaces satisfy the
convex hull property (cf.~\cite{AEG}).
\end{itemize}

\paragraph{Relation with previous results.}

\begin{itemize}
\item
Theorem~B extends \cite[Theorem~B]{EMa} to the capillary and CMC setting.
The methods developed in \cite{EMa}, which rely on the free boundary condition
and minimality, do not apply directly in this more general framework.

\item
Nevertheless, \cite[Theorem~B]{EMa} can be recovered from Theorem~B as follows.
If $\Sigma \subset \mathbb{R}^3$ is an embedded minimal annulus with free
boundaries whose support function has infinitely many critical points,
then \cite[Claim~E]{EMa} implies the existence of a neck satisfying
assumption {\rm (i)}. Furthermore, it is proved in \cite[Section 6]{EMa} that $\Sigma$ must have injective Gauss map, so the correspondence given in \cite[Proposition 6.1]{EMa} and \cite[Corollary 3.4]{Chr} implies that $\gamma$ is an analytic curve without singularities.
Moreover, \cite[Claim~D]{EMa} together with the convex hull property yield
assumptions {\rm (ii)} and {\rm (iii)}.
Hence, Theorem~B applies and implies that $\Sigma$ is either rotationally or
antipodally symmetric.
In the latter case, after rescaling, $\Sigma$ becomes an antipodally symmetric
free boundary minimal annulus in the unit ball, and
\cite[Theorem~1.1]{KM} implies that $\Sigma$ is the critical catenoid.
\end{itemize}

\section{Analytic Method: Extended pseudo-radial functions}\label{SectAM}
This section is devoted to the comparison argument leading to the proof of
Theorem~C. Throughout this section we restrict ourselves to the case in which
the nonlinearity $f$ depends only on the function $u$.
Our strategy is to introduce a family of rotationally symmetric solutions, referred to as \emph{model solutions}, and to compare the geometry
of the level sets of a general solution with those of an appropriate model. We follow the philosophy introduced in \cite{ABM,EMa}, but 
adapted here to accommodate a general nonlinearity $f$ for which no global
$P$--function is yet available. This comparison provides quantitative control on the gradient and curvature of solutions and ultimately leads to the rigidity statement of Theorem~C.

The most technical computations underlying the properties of the model solutions
and of the comparison function introduced below are collected in the Appendix,
in order to keep the presentation focused on the ideas.

\subsection{Model Solutions}\label{SectModel}
Recall that $(x,y,z)$ represent Cartesian coordinates in $\r^3$. In this section we will use cylindrical coordinates on $\mathbb{S}^2$, which are defined by
\[
\mathbb{S}^2 = \set{ \left( \sqrt{1-r^2} \cos \theta, \sqrt{1-r^2} \sin \theta, r \right) \, : \, r \in [-1,1], \, \theta \in [0,2\pi)},
\]
where the induced metric on $\mathbb{S}^2$ is given by $ g_{\mathbb{S}^2} = \frac{1}{1-r^2} \, dr^2 + (1 - r^2) \, d\theta^2$.

Consider the overdetermined elliptic problem on $\s^2$:
\begin{eqnarray} \label{OEP2}
    \left\{
    \begin{array}{llll}
        \Delta u + f(u) = 0  &\mathrm{in}~\Omega \subset \s^2,\\
        u > 0  &\mathrm{in}~\Omega, \\
        u = 0  &\mathrm{on}~\partial \Omega,\\
        \langle \nabla u, \eta \rangle = \alpha _i  &\mathrm{on}~\Gamma _i \subset \partial \Omega,
    \end{array}
    \right.
\end{eqnarray}
where \(\Omega \subset \mathbb{S}^2\) is a finite domain, $\eta$ is the outer conormal along $\partial \Omega$, \(\alpha_i < 0\) is a negative constant on each connected component $\Gamma _i \subset \partial \Omega$, $i \in \set{1,\dots ,k}$, and $f \in \mathcal{C}^1 (\r)$ is a non-negative function in $\r_+ := [0, +\infty)$. No additional regularity assumptions are required in this section. 

Radial solutions to \eqref{OEP} exist only when the domain $\Omega$ is a topological disk or an annulus. The disk case was studied in detail in \cite{EM,Sh}, and we adopt several ideas from those works. Since the Neumann condition is automatically satisfied for radial solutions, it is enough to focus on the Dirichlet problem. In cylindrical coordinates, the Laplacian operator can be expressed as
\[
\Delta =  \frac{\partial}{\partial r} \left[(1-r^2) \frac{\partial}{\partial r} \right]+\frac{1}{1-r^2} \frac{\partial ^2}{\partial \theta ^2}.
\]
Thus, if $(\Omega, u)$ is a rotationally symmetric solution to \eqref{OEP2} and we write  $u(r,\theta)=U(r)$ for some function $U$ of class $\mathcal{C}^2$, then $U$ solves the differential equation
\begin{equation}\label{ODE}
	(1-r^2) U'' (r)- 2 r U' (r)+f(U(r))=0.
\end{equation}with the Cauchy conditions 
\begin{equation}\label{CauchyData}
	U(R)=M \quad \textup{and} \quad U' (R)=0,
\end{equation}
for some $R \in [-1,1]$ and $M \in \r_+^* := (0,+\infty)$.

We summarize a few important properties, proven via ODE analysis in Appendix \ref{Appendix_theor_ODE}, of rotationally symmetric solutions:

\begin{theorem}\label{theor_ODE}
Let $f \in \mathcal{C} (\r)$ be a continuous function that is positive in $\r_+^*$. Then, for any $(R,M) \in [-1,1] \times \r_+^*$, there exists a solution $U_{R,M,f}$ to \eqref{ODE}-\eqref{CauchyData}, satisfying the following properties: 
	\begin{enumerate}
		\item $U_{R,M,f}$ is unique. Furthermore, it follows that $U_{-R,M,f} (-r)=U_{R,M,f} (r)$ for each $R$ (being $U_{0,M,f}$ symmetric with respect to $r=0$), so we can consider $R \in [0,1]$.
		\item When $R \in [0,1)$, there exist values $r_\pm := r_\pm (M,R,f)$ satisfying
		$$-1< r_- <R<r_+  < 1$$and such that $U_{R,M,f}(r_+)=U_{R,M,f}(r_+)=0$, $U_{R,M,f}>0$ inside $(r_-,r_+)$, $U_{R,M,f} ' >0$ inside $[r_-,R)$ and $U_{R,M,f} ' <0$ in $(R,r_+]$. 
		\item When $R =1$, there exists a value $r_-(1,M,f)<1$ such that $U_{1,M,f} (r_-(1,M,f)) =0$, $U_{1,M,f} ' >0$ in $[r_-(1,M,f),1)$. Indeed, 
		$$\lim _{R \to 1^-}r_+ (R,M,f)= 1 \text{ and } \lim _{R \to 1^-}r_-(R,M,f)= r_-(1,M,f),$$ which coincides with the value corresponding to the model solution on the geodesic disk (cf. \cite{EM}). 
	\end{enumerate}
	If we suppose further that $f$ is $\mathcal{C}^1$ and:
	\begin{equation}\label{conditionf}
		f(x) \geq x f '(x), \quad \textup{ for all }x>0,
	\end{equation}
	then we have also that:
	\begin{enumerate}
		\item[4.] When $R \in [0,1)$, both functions $r_-$ and $r_+$ depend on the parameters $R$ and $M$ in a $\mathcal{C}^1$ way. Moreover, both functions are non-decreasing in $R$. 
		\item[5.] Given $M \in \mathcal \r_+^*$ fixed, consider the real valued functions
		\begin{equation}\label{functionNormOfGradient}
			g_{M,f}^{\pm} (R) :=(1-r_{\pm}^2) U_{R,M,f}' (r_{\pm})^2 , \quad \forall R \in [0,1).
		\end{equation}
	Then for any $M \in \r_+^*$, the functions $g_{M,f}^- $ and $g_{M,f}^+$ are decreasing and increasing on $[0, 1)$, respectively. 
		\item[6.] The function 
	\begin{equation}\label{hRotacional}
	h_f(M) := \sqrt{1-r_{-}^2} \abs{U_{R,M,f}' (r_{-})}, \qquad M >0
	\end{equation}
is of class $\mathcal{C}^1$ and increasing in $M$.
	\end{enumerate}
	If we suppose further that either $f(x)=1$ or $f \in \mathcal{C}^1(\mathbb{R})$ and:
	\begin{equation}\label{conditionf2}
		f'(x) \geq 2, \quad \textup{ for all }x > 0,
	\end{equation}
	then we also have that:
	\begin{enumerate}
	\item[7.] Given $(R,M )\in [0, 1) \times \mathcal \r_+^*$ fixed, the real valued function 
	$$\mu_{R,M,f} (r) := 2 U''_{R,M,f}(r)^2 - U'_{R,M,f} U''' _{R,M,f}(r),$$ is non-negative on $[r_- , r_+]$.
	\end{enumerate}
\end{theorem}
\begin{remark}
	The existence and properties 1--3 of a solution to \eqref{ODE}--\eqref{CauchyData} stated in the previous result also remain valid if, instead of assuming that the nonlinearity $f$ is positive on $\mathbb{R}_+^*$, we assume that there exists an open interval $\mathcal{I}_f \subset \mathbb{R}_+^*$ with $0 \in \partial \mathcal{I}_f$ such that $f>0$ on $\mathcal{I}_f$ and $M \in \mathcal{I}_f$. 
	
	Moreover, in order to deduce properties 4--7 for a solution of \eqref{ODE}--\eqref{CauchyData}, it is sufficient to assume that conditions \eqref{conditionf} and \eqref{conditionf2} hold on a smaller interval, provided that $M$ belongs to that interval.
\end{remark}

Thus, Theorem \ref{theor_ODE} shows that the problem \eqref{OEP2} admits a $2-$parameter family of rotationally symmetric model solutions, where the parameters $M$ and $R$ that describe this family are the maximum value of the solutions and the height of the parallel in which the solutions attains this value. 

\begin{definition}\label{modelSolutions}
	Using cylindrical coordinates, let us define $$u_{R,M,f} (r, \theta)= U_{R,M,f}(r), \quad \forall \theta \in \s^{1}$$ and $$	\Omega_{R,M,f}:=
	\begin{cases}
		\{\,(\sqrt{1-r^2}\cos \theta , \sqrt{1-r^2} \sin \theta ,r) \in \s^2 :   r_{-}( 1,M,f) <  r \,\}, & R= 1,\\[2pt]
		\set{(\sqrt{1-r^2}\cos \theta , \sqrt{1-r^2} \sin \theta ,r) \in \s^2 ~\colon~ r \in (r_-, r_+)} & R \in[0,1).
	\end{cases}
	$$ where $(R,M) \in [0,1] \times \r_+^*$ and $ U_{R,M,f}$ is a solution to \eqref{ODE}-\eqref{CauchyData}. Then we say that the pair $(\Omega_{R,M,f}, u_{R,M,f})$ is a \textup{model solution}.
\end{definition}

We now fix some notation for the components of any of the above model solutions. When $R \in [0, 1)$, we write the connected components of $\partial \Omega_{R,M,f}$ as 
\begin{equation}\label{boundaryComponent}
\Gamma_{R,M,f}^{\pm} = \set{(\sqrt{1-r^2} \cos \theta , \sqrt{1-r^2} \sin \theta ,r) \in \s^2 ~\colon~ r =r_\pm}, 
\end{equation}
and we also consider the subdomains $\U_{R,M,f}^{\pm} \subset \Omega_{R,M,f} \setminus \textup{Max} (u_{R,M,f})$, where $\textup{cl} (\U_{R,M,f}^{\pm}) \cap \partial \Omega_{R,M,f} = \Gamma_{R,M,f}^{\pm}$. When $R=1$, we simply have $\U_{1,M,f} \subset \Omega_{1,M,f} \setminus \textup{Max} (u_{1,M,f})$ and
$$ \Gamma_{1,M,f} :=\set{(\sqrt{1-r^2} \cos \theta , \sqrt{1-r^2} \sin \theta ,r) \in \s^2 ~\colon~ r =r_-} = \textup{cl} (\U_{1,M,f}) \cap \partial \Omega_{1,M,f}.$$

\paragraph{The Model $\overline{\tau}$-function.} 
For a fixed value $M \in \r_+^*$, we introduce the quantities
\begin{equation}\label{tauModelSolutions}
	\overline{\tau}_{M,f}^\pm (R)
	:= \frac{(1-r_{\pm} (R,M,f)^2)\, U_{R,M,f}' \big(r_{\pm} (R,M,f)\big)^2}
	{(1-r_{-} (1,M,f)^2)\, U_{1,M,f}' \big(r_{-} (1,M,f)\big)^2}
	= \frac{\abs{\nabla u_{R,M,f}}^2 \big|_{\Gamma_{R,M,f}^{\pm}}}{h_f(M)^2},
	\quad R \in [0,1),
\end{equation}
which compare the squared gradient of a ring-type model solution along its zero level sets with the squared gradient of a disk-type model solution along its boundary.
We refer to these functions as the \textit{model $\overline{\tau}$-functions}.

As will be shown in Subsection~\ref{sectionTau}, a crucial step in establishing a one-to-one correspondence between general solutions to the Dirichlet problem associated to \eqref{OEP2} and the corresponding model solutions is to prove that the functions $\overline{\tau}_{M,f}^{\pm}$ are diffeomorphisms.
In the works \cite{ABM,EMa}, this property follows from explicit computations.
By contrast, for a general nonlinearity $f$, such explicit formulas are no longer available, and the proof must rely on the qualitative analysis provided by Theorem~\ref{theor_ODE}.

We define
\[
\tau_f^0 (M) := \overline{\tau}_{M,f}^- (0) = \overline{\tau}_{M,f}^+ (0),
\qquad \forall\, M \in \r_+^*.
\]
The following result summarizes the monotonicity properties of the model $\overline{\tau}$-functions; its proof is deferred to Appendix~\ref{Appendix_lemmaTau}.
\begin{lemma}\label{lemmaTau}
Assume that $f\in \mathcal{C}^1(\mathbb{R})$ is a function that is positive in $\mathbb{R}^*_+$ and satisfies both conditions \eqref{conditionf} and \eqref{conditionf2}. Then, for any $M\in \r_+^*$, the functions
\[
\overline{\tau}_{M,f}^- : [0,1] \longrightarrow [1,\tau_f^0(M)]
\quad \text{and} \quad
\overline{\tau}_{M,f}^+ : [0,1) \longrightarrow [\tau_f^0(M),+\infty),
\]
defined in \eqref{tauModelSolutions}, are decreasing and increasing, respectively.
\end{lemma}

\subsection{The $\overline{\tau}$-function and expected critical height}\label{sectionTau}
Having introduced the family of rotationally symmetric model solutions, we now turn to the problem of associating such a model to a given analytic solution $(\Omega ,u)$ of the Dirichlet problem
\begin{eqnarray} \label{DP}
	\left\{
	\begin{array}{llll}
		\Delta u + f(u) = 0  &\text{in}~\Omega \subset \s ^2,\\
		u > 0  &\text{in}~\Omega, \\
		u = 0  &\text{on}~\partial \Omega ,
	\end{array}
	\right.
\end{eqnarray}
where we assume that $f$ is $\mathcal C^\omega$ (analytic) and satisfies (recall conditions \eqref{conditionf} and \eqref{conditionf2}): 
$$ f(x) \geq x f'(x)  \text{ and } f'(x) \geq 2 \text{ for all } x >0 .$$

\begin{remark}\label{RemCond}
Throughout this section we work with analytic solutions of \eqref{DP}. This choice streamlines the exposition and avoids technicalities. Several steps below remain valid under weaker regularity assumptions on $\Omega$, $u$, and $f$, but we do not pursue this generality here. The computations follow closely \cite[Section~4.1]{EMa}.

Since $f(x)>0$ for all $x>0$ under our assumptions, $\textup{Max}(u)$ is the union of isolated points and, possibly, simple closed curves; moreover, every such curve is an isolated component of $\textup{Max}(u)$ (see \cite[Corollary~3.4]{Chr}). In addition, using the same maximum--principle arguments as in the proof of Corollary A, we deduce that every curve in $\textup{Max}(u)$ is non--contractible (see also \cite[Lemma 3.1]{EMa}).

\end{remark}

More precisely, let $(\Omega,u)$ be an analytic solution to \eqref{DP} and let $\U$ be a connected component of $\Omega\setminus \textup{Max}(u)$. Our aim is to associate to the pair $(\U,u)$, in a canonical way, a corresponding component $(\U_{R,M,f}^{\pm},u_{R,M,f})$ of a model solution, where $M=u_{\textup{max}}$. This association is achieved through an invariant quantity, the $\overline{\tau}$--function, which measures the boundary gradient in a normalized form and, in particular, encodes the relative size of $\U$ within the family of model domains.

\begin{definition}\label{DefinTau}
Let $(\Omega, u)$ be an analytic solution to \eqref{DP}, and let $\Gamma \in \pi_{0} (\partial \Omega)$ be a connected component of the boundary. Then we define
\begin{equation}\label{tau1}
\overline{\tau}(\Gamma):=  \max \set{
\frac{\abs{\nabla u}^2 (p)}{ h_f(u_{\rm max})^2} ~\colon~ p \in \Gamma }.
\end{equation}where $h_f$ is defined in \eqref{hRotacional}.

If $\mathcal{U}$ is a connected component of $\Omega \setminus \textup{Max} (u)$, $\partial \Omega \cap \textup{cl}(\U) \neq \emptyset$, we extend the previous definition as
	\begin{equation*}
		\overline{\tau}(\mathcal{U}):= \max \left\{ \overline{\tau}(\Gamma) ~\colon~ \Gamma \in \pi_{0} ( \partial \Omega \cap \textup{cl}(\U) ) \right\}.
	\end{equation*}
Otherwise, we set $\overline{\tau}(\mathcal{U})=0$.
\end{definition}
\begin{remark}
This definition extends \cite[Definition~3.3]{EMa}, up to taking a square root. Indeed, in that paper we defined
\[
\overline{\tau} (\Gamma) = \max \left\{
\frac{\abs{\nabla u}(p)}{u_{\rm max}} \;:\; p \in \Gamma
\right\}
\]
for a solution $(\Omega, u)$ of \eqref{DP} with $f(x)=2x$. The normalization factor used there is motivated by the $P$--function associated with the problem, namely $P = \abs{\nabla u}^2 + u^2$ (see \cite[Proposition~3.1]{EMa}). In this setting, however, it is straightforward to verify that
\[
h_f(M) = \abs{\nabla u_{1,M,f}} |_{\Gamma_{1,M,f}} = M, 
\quad \forall\, M \in \mathbb{R}_+^*.
\]
Therefore, in this case the quantity defined above coincides with \eqref{tau1}, up to squaring.
\end{remark}

From the definition of the $\overline{\tau}$-function in \eqref{tau1}, it follows immediately that
\[
\overline{\tau} (\Gamma_{R,M,f}^{\pm}) = \overline{\tau}_{M,f}^\pm  (R), 
\qquad \forall\, R \in [0,1), \; M \in \r_+^*,
\]
where $\overline{\tau}_{M,f}^\pm $ are the model $\overline{\tau}$-functions defined in \eqref{tauModelSolutions}, and 
$\Gamma_{R,M,f}^{\pm}$ denote the boundary components of the domain of a model solution (see \eqref{boundaryComponent}). Then Lemma \ref{lemmaTau} yields the following. For any solution $(\Omega,u)$ of \eqref{DP} with $M=u_{\textup{max}}$ and any connected component $\mathcal{U}\in \pi_0(\Omega \setminus \textup{Max}(u))$ such that $\overline{\tau}(\mathcal{U})\ge 1$, there exist a unique $R\in[0,1]$ and a unique model domain $\bar{\mathcal{U}}:=\U^{\pm}_{R,M,f}$ satisfying
\[
\overline{\tau}(\mathcal{U})=\overline{\tau}(\bar{\mathcal{U}}).
\]
This provides a canonical correspondence between general solutions of \eqref{DP} and the family of model solutions through the $\overline{\tau}$--function. The precise formulation is given in the following definition, which is based on \cite[Definition 3.5]{EMa} (see also \cite[Definition 1.4]{ABM}).

\begin{definition}\label{defCriticalHeight}
Let $(\Omega,u)$ be an analytic solution to \eqref{DP}. Set $M:=u_{\textup{max}}$ and fix a connected component $\U\in \pi_0\big(\Omega\setminus \textup{Max}(u)\big)$.
We define the \textup{expected critical height of} $\mathcal U$ as follows:
		\begin{enumerate}
			\item[(i)]  If $\overline{\tau} (\mathcal U) \geq \tau_f^0 (M)$, we set $\bar R(\mathcal U)= (\overline{\tau}_{M,f}^+)^{-1} \left( \overline{\tau} (\mathcal U) \right)$.
		\item[(ii)] If $\overline{\tau} (\mathcal U) \in (1,\tau_f^0 (M))$, we set $ \bar R(\mathcal U)= (\overline{\tau}_{M,f}^-)^{-1} \left( \overline{\tau} (\mathcal U) \right)$.
		
	   \item[(iii)]  If $\overline{\tau} (\mathcal U) \leq 1$, we set $\bar R(\mathcal U)= 1$.
	\end{enumerate}
\end{definition}

Observe that Definition~\ref{defCriticalHeight} differs from \cite[Definition 3.5]{EMa}. In \cite{EMa} we proved that the alternative $\overline{\tau}(\mathcal{U})\le 1$ for some $\mathcal{U}\in \pi_0(\Omega\setminus \textup{Max}(u))$ occurs if and only if $(\Omega,u)$ is the geodesic disk solution (cf.~\cite[Theorem 3.2]{EMa}). For this reason, \cite[Definition 3.5]{EMa} did not include a condition analogous to item ${\it (iii)}$ in Definition~\ref{defCriticalHeight}, since the case $\overline{\tau}(\mathcal{U})\le 1$ was already completely classified.

In the general setting of \eqref{DP}, however, no $P$--function is yet available, and therefore the argument used in the proof of \cite[Theorem 3.2]{EMa} does not apply. Nevertheless, we will prove that $\overline{\tau}(\mathcal{U})>1$ whenever $\textup{cl}(\mathcal{U}) \cap \textup{Max}(u)$ contains infinitely many points (cf.~Proposition~\ref{prop_boundTau}). This will be established in Section~\ref{sectionConsequences} by a geometric argument.

\subsection{Pseudo-radial functions}\label{sectionPseudoRadial}

Now that we have established a procedure to associate a unique model solution to a given general solution of \eqref{DP}, we are in a position to introduce the main tool used to carry out comparison arguments: the \emph{pseudo-radial functions}. In order to introduce this concept we first need to present the following definition, which generalizes \cite[Definition 4.1]{EMa}. 

\begin{definition}\label{associatedSolution}
Let $(\Omega,u)$ be an analytic solution to \eqref{DP}. Set $M:=u_{\textup{max}}$ and fix a connected component $\U\in \pi_0\big(\Omega\setminus \textup{Max}(u)\big)$.

We say that $\big(\U^{\diamond}_{\bar R,M,f},u_{\bar R,M,f}\big)$ is a \textup{comparison pair} associated to $(\U,u)$ if there exist $\bar R\in[0,1]$ and $\diamond\in \set{-,+}$ such that
\[
\overline{\tau}(\U)\le \overline{\tau}\big(\U^{\diamond}_{\bar R,M,f}\big).
\]

If, moreover, $\bar R=\bar R(\U)$ is the expected critical height of $\U$ and
\[
\overline{\tau}\big(\U^{\diamond}_{\bar R,M,f}\big)=\overline{\tau}(\U),
\]
then we say that $\big(\U^{\diamond}_{\bar R,M,f},u_{\bar R,M,f}\big)$ is the \textup{model pair} associated to $(\U,u)$.
\end{definition}

\begin{remark}\label{remarkNotation}
From now on, we omit the dependence of the parameters in Definition~\ref{associatedSolution}, writing $(\bar \U,\bar u)$ for a comparison pair or a model pair. We also denote by $\bar U$ the solution to the O.D.E.~\eqref{ODE} defining $\bar u$. We will also write $(\bar \Omega, \bar u)$ for the model solution of which $(\bar \U,\bar u)$ is a model pair. Finally, we denote by $\bar r_\pm$ the values where $\bar U$ vanishes, that is, $\bar U(\bar r_\pm)=0$.
\end{remark}


Set $M=u_{\textup{max}}$ and choose a comparison pair $(\bar\U,\bar u)$ associated to $(\U,u)$ as in Definition \ref{associatedSolution}. Consider
\[
G:[0,M]\times[\bar r_-,\bar r_+]\to\mathbb{R},\qquad G(u,r)=u-\bar U(r).
\]
Since $\frac{\partial G}{\partial r}(u,r)=-\bar U'(r)$ and $\bar U'(r)=0$ occurs precisely at
$r=\bar R$, the equation $G(u,r)=0$ can be solved for $r$ along the two monotone branches of $\bar U$.
Therefore, by the Implicit Function Theorem, there exist two differentiable functions
\[
\chi_-:[0,M]\to[\bar r_-,\bar R]
\qquad\textup{and}\qquad
\chi_+:[0,M]\to[\bar R,\bar r_+]
\]
such that
\[
G\big(u,\chi_\pm(u)\big)=0 \quad\textup{for all}\quad u\in[0,M].
\]

\begin{figure}[htb]
	\centering
	\includegraphics[width=1\textwidth]{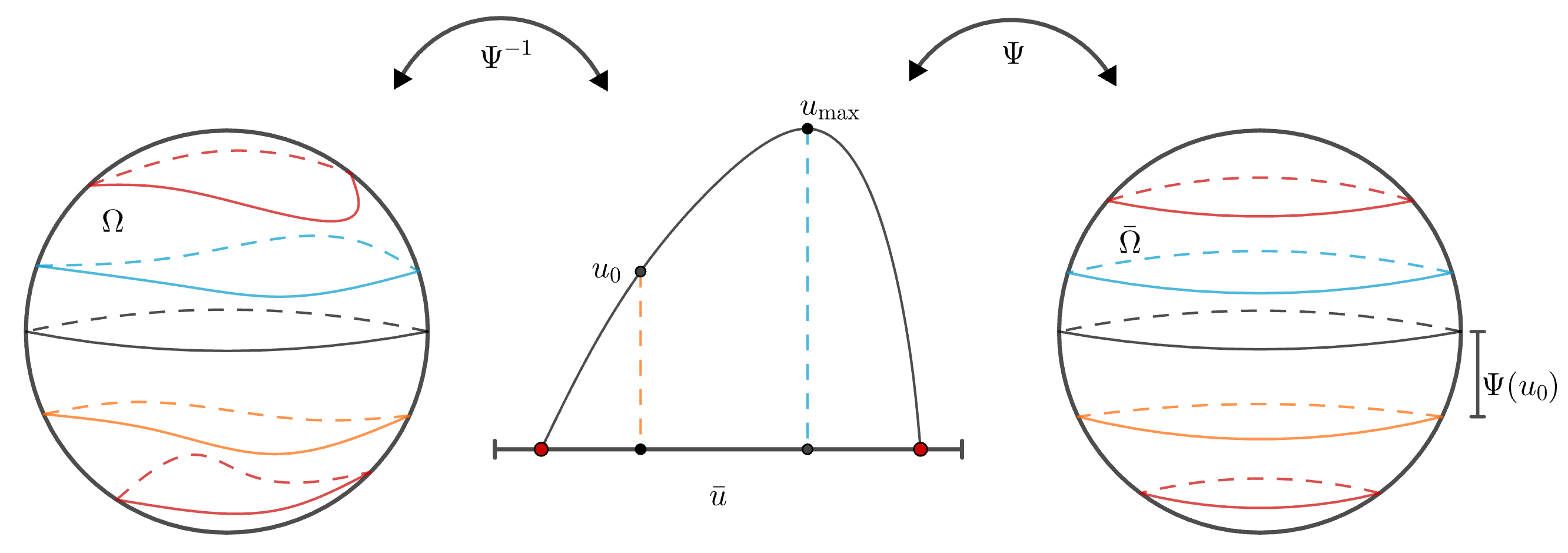}
	\caption{This figure illustrates the concept of pseudo-radial functions $\Psi_{\pm}$, as defined in Definition~\ref{psudoRadialFunctions}. Given a solution $(\Omega, u)$ to problem~\eqref{DP}, and a connected component $\mathcal{U}$ of $\Omega \setminus \textup{Max}(u)$, we associate a comparison pair $(\bar{\mathcal{U}}, \bar{u})$. The pseudo-radial functions $\Psi_{\pm}$ map each point $p \in \mathcal{U}$ to a corresponding height coordinate in the model solution via $\Psi_{\pm}(p) = \chi_{\pm}(u(p))$, depending on whether $\overline{\tau}(\bar{\mathcal{U}})$ is less than or greater than $\tau_f^0(M)$. This mapping establishes a relationship between the original and the model solution, effectively relating the level sets of $u$ in $\mathcal{U}$ to heights with respect to the equator in the model domain $\bar{\mathcal{U}}$. This correspondence allows us to define, in the original solution, quantities that depend on the model solution—for example, the modulus of the gradient of the model solution.}	
	\label{fig:PseudoRadial}
\end{figure}

\begin{definition}\label{psudoRadialFunctions}
Let $(\Omega,u)$ be an analytic solution to \eqref{DP} and fix a connected component $\U\in \pi_0\big(\Omega\setminus \textup{Max}(u)\big)$. Let $(\bar\U,\bar u)$ be a comparison pair associated to $(\U,u)$ (cf. Figure~\ref{fig:PseudoRadial}). Then:
	\begin{itemize}
		\item If $\overline{\tau}(\bar \U) < \tau_f^0 (M)$, we define the pseudo-radial function $\Psi_{-}$ associated to $(\U,u)$ and $(\bar \U,\bar u)$ by
		\begin{equation*}
			\begin{split}
				\Psi_{-}:  \, & \mathcal U \to [ \bar r_{-}, \bar R] \\
				& \, p \,\, \mapsto \Psi_{-} (p):= \chi_{-}\left(u (p)\right).
			\end{split}
		\end{equation*}
		\item If $\overline{\tau}(\bar  \U) \geq \tau_f^0 (M)$, we define the pseudo-radial function $\Psi_{+}$ associated to $(\U,u)$ and $(\bar \U,\bar u)$ by
		\begin{equation*}
			\begin{split}
				\Psi_{+}:  \, & \mathcal U \to [\bar R,\bar r_{+}] \\
				& \, p \,\, \mapsto \Psi_{+} (p):= \chi_{+}\left(u (p)\right).
			\end{split}
		\end{equation*}
	\end{itemize}
\end{definition}

\begin{remark}\label{remarkNotation2}
Following the notation of \cite{EMa}, from now on, we write $\Psi_{\pm}=\Psi$ and $\chi_{\pm}=\chi$, and we treat the two cases simultaneously. Derivatives with respect to $u$ will be denoted by a dot. For later use, we also introduce the auxiliary function
\begin{equation}\label{auxiliarFunction}
	\varphi(r)= \frac{f(\bar U (r))-2r \bar{U}' (r)}{(1-r^2) \bar{U}' (r)^2} = - \frac{\bar U '' (r)}{\bar U '(r)^2}, \quad \forall r \in (\bar r_-, \bar R) \cup (\bar R, \bar r_+),
\end{equation}
and we define
\begin{equation}\label{normGradient}
	\bar W (p)= (1-\Psi ^2 (p)) \bar{U}' (\Psi (p))^2 \quad \text{ for } p \in \mathcal U. 
\end{equation}
\end{remark}

\subsection{Gradient estimates}\label{sectionGradient}

Fix a comparison pair and the associated pseudo--radial functions. We can then compare analytic quantities attached to $u$ with their model counterparts and derive a sharp gradient estimate.

The following lemma is obtained by a lengthy computation, closely paralleling the argument in the proof of \cite[Theorem~3.1]{EMa}. To keep the main text readable, we postpone the proof to Appendix~\ref{Appendix_LemDeltaF}.

\begin{lemma}\label{LemDeltaF}
Let $(\Omega,u)$ be an analytic solution to \eqref{DP} and fix a connected component $\U\in \pi_0\big(\Omega\setminus \textup{Max}(u)\big)$. Let $(\bar\U,\bar u)$ be a comparison pair associated to $(\U,u)$ and let $\bar U$ be the profile function defining $\bar u$ on $[\bar r_-,\bar r_+]$. Then the function
\[
F := \abs{\bar{U}' (\Psi)}^{-1}\big(|\nabla u|^2-\bar{W}\big)
\]
is analytic and satisfies the elliptic inequality
\begin{equation}\label{DeltaF}
\Delta F  - \frac{\mu(\Psi)}{\bar{U}' (\Psi)^2} F \geq 0 \text{ in } \mathcal U,
\end{equation}
with $F\leq 0$ along $\partial \mathcal U$. Here, $ \mu (\Psi) := 2 \bar U'' (\Psi)^2 - \bar U' (\Psi) \bar U''' (\Psi)$ and $\bar W$ is given by \eqref{normGradient}.
\end{lemma}

We now prove the main result of this subsection.

\begin{theorem}\label{theor_Gradient}
Let $(\Omega, u)$, $\Omega \subset \s ^2$ a finite domain, be an analytic solution to \eqref{OEP}, where we assume that either $f\equiv 1$ or $f\in\mathcal{C}^{\omega}$ satisfies $f(x)\geq x f'(x)$ and $f'(x)\geq 2$ for all $x>0$. Fix a connected component $\U\in \pi_0\big(\Omega\setminus \textup{Max}(u)\big)$ and let $(\bar\U,\bar u)$ be a comparison pair associated to $(\U,u)$. 

Then,
\[
\abs{\nabla u}^2(p) \leq \bar{W} (p) \quad \text{ for all } p \in \mathcal U.
\]
Moreover, if equality holds at some point of $\mathcal U$, then $(\Omega,u)$ is congruent to $(\bar \Omega,\bar u)$.
\end{theorem}

\begin{proof}
By the assumptions on $f$ we have $f(x)>0$ for all $x>0$, and therefore $(\Omega,u)$ is, in particular, a solution to \eqref{DP}. Moreover, by Item~7 of Theorem~\ref{theor_ODE} we have $\mu(\Psi)\ge 0$ in \eqref{DeltaF}. 

Consequently, $F$ is a supersolution of a linear elliptic operator with nonnegative zeroth--order coefficient. Since $F\le 0$ on $\partial\U$, the maximum principle yields $F\le 0$ in $\U$, and therefore $|\nabla u|^2\le \bar W$ throughout $\U$.

Assume now that $|\nabla u|^2(p)=\bar W(p)$ at some point $p\in\U$. Then $F(p)=0$, and the strong maximum principle implies $F\equiv 0$ in $\U$, that is, $|\nabla u|^2\equiv \bar W$ on $\U$. In particular, \eqref{laplacianWs2} holds with equality, and therefore
\[
\nabla^2 u = \frac{1}{2}\left(\varphi(\Psi)\abs{\nabla u}^2- f(u)\right) g_{\sp^2} -\varphi(\Psi)\, du \otimes du,
\]
where $\varphi$ is defined as in \eqref{auxiliarFunction}. It follows from the last equation that the level sets of $u$ in $\U$ have constant curvature; hence they are parallel circles in $\s^2$ centered at some point. Consequently, $u$ agrees with its associated model solution on $\U$, and by analytic continuation we conclude that $(\Omega,u)$ and $(\bar\Omega,\bar u)$ are congruents.
\end{proof}

\begin{remark}\label{remarkComparisonFunction}
In \cite{ABM,EMa}, the gradient estimate is formulated in terms of \textup{model pairs}, namely pairs $(\mathcal{U},u)$ and $(\bar{\mathcal{U}},\bar u)$ satisfying $\overline{\tau}(\mathcal{U})=\overline{\tau}(\bar{\mathcal{U}})$. Theorem~\ref{theor_Gradient}, however, shows that this matching condition is not essential for the comparison procedure. In fact, it suffices to work with a \textup{comparison pair} in the weaker sense that
\[
\overline{\tau}(\mathcal{U}) \le \overline{\tau}(\bar{\mathcal{U}}).
\]
\end{remark}

\subsection{Consequences of Theorem \ref{theor_Gradient}}\label{sectionConsequences}

In this subsection we record several geometric consequences of the gradient estimate in Theorem~\ref{theor_Gradient}.
We first bound the geodesic curvature of the level sets of $u$ in terms of the corresponding level circles of a comparison pair; see Propositions~\ref{prop_curvatureZeroSets} and \ref{prop_CurvatureMaxSets}.
These curvature bounds allow us to rule out certain configurations, and in particular imply that if $\Max(u)$ contains a closed curve then $\overline{\tau}(\U)>1$; see Proposition~\ref{prop_boundTau}.
We then derive sharp length estimates for level sets (Propositions~\ref{prop_lengthBounds} and \ref{propConvexity}). Throughout this subsection we work under the hypotheses of Theorem~\ref{theor_Gradient}.

\paragraph{Curvature estimates.}

We compare the curvature of the level sets of the general function $u$ inside $\U$ with that of its associated model solution or, more generally, a comparison pair (cf. \cite[Subsection 4.2]{EMa}). We will use the standard Landau notation as follows: if $g:\mathcal{M}\to\r$ is real analytic function (with $\mathcal{M}=\s^2$ or $\r$) and $m\in\mathbb{N}$, then $O(g^m)$ denotes a real analytic function of $g$ such that $O(g^m)/g^m$ converges to a finite constant as $g\to 0$.

\begin{proposition}\label{prop_curvatureZeroSets}
Let $(\Omega, u)$, $\Omega \subset \s ^2$ a finite domain, be an analytic solution to \eqref{OEP}, where we assume that either $f\equiv 1$ or $f\in\mathcal{C}^{\omega}$ satisfies $f(x)\geq x f'(x)$ and $f'(x)\geq 2$ for all $x>0$. Set $M:=u_{\textup{max}}$. Fix a connected component $\U\in \pi_0\big(\Omega\setminus \textup{Max}(u)\big)$ and let $(\bar\U,\bar u)$ be a comparison pair associated to $(\U,u)$. 

Let $p \in \partial \Omega$ be a point such that
	\begin{equation}\label{condNabla}
		\abs{\nabla u}^2 (p) = \max_{\partial \Omega \cap \textup{cl}(\U)} \abs{\nabla u}^2.
	\end{equation}
If $\kappa (p)$ denotes the (geodesic) curvature of $\partial \Omega$ at $p$ with respect to the inner orientation to $\mathcal{U}$, then
	\begin{equation*}
		\kappa (p) \leq -\frac{\bar{r}_+}{\sqrt{1-\bar{r}_+^2}} \quad \textup{if } \overline{\tau} (\bar{\mathcal{U}}) \geq \tau_f^0(M) \quad \textup{and} \quad \kappa (p) \leq \frac{\bar{r}_{-}}{\sqrt{1-\bar{r}_{-}^2}} \quad \textup{if } \overline{\tau} (\bar{\mathcal{U}}) < \tau_f^0(M).
	\end{equation*}
\end{proposition}

\begin{proof}
	We follow the proof of \cite[Proposition 4.1]{EMa}; that is, we study the Taylor expansions of the functions $ |\nabla u|^2$ and $\bar W$ defined in \eqref{normGradient} in a neighborhood of $\partial \Omega \cap \textup{cl}(\U)$.
	
	Let $p \in \partial \Omega \cap \textup{cl}(\U)$ be a point satisfying \eqref{condNabla} and define $\sigma : [0, 2 \pi) \to \s^2$ to be the unit speed geodesic with initial data $\sigma (0) =p$ and $\sigma ' (0) = \nabla u / \abs{\nabla u} (p)$. Then, since $\Delta u (p) = -f (0)$, the Taylor expansion of the real-valued function $\abs{\nabla u}^2 \circ \sigma$ is given by
\begin{equation}\label{TaylorW}
\abs{\nabla u}^2 (\sigma (s)) = \abs{\nabla u}^2 (p) + 2 \left( \abs{\nabla u} (p) \kappa(p)-f(0) \right) \abs{\nabla u} (p) s + O(s^2),
\end{equation}
On the other hand, using the definition of $\Psi$ and equation \eqref{ODE} it is easy to compute
\[
\nabla \bar W = 2 (\Psi \bar U' (\Psi) - f(u)) \nabla u,
\]
so the corresponding Taylor expansion of $\bar W \circ \sigma$ at $p$ is
\begin{equation}\label{TaylorModelW}
	\bar W (\sigma (s)) 
	= \bar W (p) 
	+ 2 \left( \mp \frac{\bar r_\pm}{\sqrt{1- \bar r_\pm^2}}\, \sqrt{\bar W(p)} - f(0) \right) 
	\sqrt{\bar W(p)}\, s 
	+ O(s^2),
\end{equation}
where the minus sign is chosen when $\overline{\tau} (\bar \U) \geq \tau_f^0(M)$, and the plus sign when $\overline{\tau} (\bar \U) < \tau_f^0(M)$. Since $\abs{\nabla u}^2 (p) = \bar W (p)$ by condition \eqref{condNabla} and $\abs{\nabla u}^2 \leq \bar W$ by Theorem \ref{theor_Gradient}, we conclude the proof of the proposition by subtracting \eqref{TaylorModelW} from \eqref{TaylorW}.
\end{proof}

\begin{remark}\label{remarkSignOfR1}
By Proposition~\ref{prop_curvatureZeroSets}, we can ensure that a boundary component contained in $\textup{cl}(\mathcal{U})\cap\partial\Omega$ will be convex (with respect to the inner orientation to $\mathcal{U}$) only if $\bar r_-<0$.
In general, however, this sign condition fails.

To see this, consider the family of model solutions $\{(\Omega_{R,M,f},u_{R,M,f})\}_{R\in[0,1]}$. For the ball solution ($R=1$), the maximum principle implies that $\Omega_{1,M,f}$ is contained in a closed hemisphere of $\s^2$, say $\Omega_{1,M,f}\subset \s^2_+$.
Assume, for contradiction, that $r_-(1,M,f)<0$.
Let $v\in\mathcal{C}^{\omega}(\s^2)$ be a first eigenfunction of the Laplacian on $\s^2$, positive in $\s^2_+$ and vanishing along $\partial\s^2_+$.
Normalize $v$ so that $v<u$ in $\s^2_+$, and let $\lambda_0>0$ be the smallest value for which the graphs of $u$ and $\lambda_0 v$ touch.
Then $w:=u-\lambda_0 v$ satisfies $w\ge 0$ in $\s^2_+$ and attains a zero minimum at an interior point.

Since $f(u)\ge 2u$, we have
\[
\Delta w=-f(u)+2\lambda_0 v \le -2(u-\lambda_0 v)=-2w\le 0,
\]
and the strong maximum principle yields $u\equiv \lambda_0 v$, contradicting $r_-(1,M,f)<0$.
Therefore $r_-(1,M,f)\ge 0$.
In particular, for $R$ sufficiently close to $1$, the curve $\Gamma_{R,M,f}^-$ lies entirely above the equator whenever $f(x)\neq 2x$.
\end{remark}

Now we consider the case in which $\textup{Max} (u) \cap \textup{cl}(\U)$ contains a curve. In this situation, we can estimate the curvature of the curve using the same argument as in \cite[Proposition 4.2]{EMa}.

\begin{proposition}\label{prop_CurvatureMaxSets}
Let $(\Omega, u)$, $\Omega \subset \s ^2$ a finite domain, be an analytic solution to \eqref{OEP}, where we assume that either $f\equiv 1$ or $f\in\mathcal{C}^{\omega}$ satisfies $f(x)\geq x f'(x)$ and $f'(x)\geq 2$ for all $x>0$. Set $M:=u_{\textup{max}}$. Fix a connected component $\U\in \pi_0\big(\Omega\setminus \textup{Max}(u)\big)$ and let $(\bar\U,\bar u)$ be a comparison pair associated to $(\U,u)$ with critical height $\bar R<1$. 

Assume that $\textup{cl}(\U)\cap \textup{Max}(u)$ contains a curve $\gamma$. If $\kappa(p)$ denotes the (geodesic) curvature of $\gamma$ at $p$ with respect to the inner orientation to $\U$, then
	\begin{equation*}
		\kappa (p) \leq \frac{\bar R}{\sqrt{1-\bar{R}^2}} \quad \textup{if } \overline{\tau} (\bar{\mathcal{U}}) \geq \tau_f^0(M) \quad \textup{and} \quad \kappa (p)\leq -\frac{\bar{R}}{\sqrt{1-\bar{R}^2}} \quad \textup{if } \overline{\tau} (\bar{\mathcal{U}}) < \tau_f^0(M).
	\end{equation*}
\end{proposition}
\begin{proof}
First, observe that $\gamma$ is an analytic simple closed curve (cf. Remark \ref{RemCond}). Hence, by defining $d=\textup{dist}_{\sp^2} (\cdot, \gamma)$ in $\U$, Lemma \ref{lemma_Taylor} in Appendix \ref{AppTaylorMax} provides us with the following Taylor expansions in a neighborhood of $\gamma$:
\[
\begin{split}
	&	\abs{\nabla u}^2 = f(M)^2  \left( 1+ \kappa (p) d \right) d^2+ O(d^4), \\
	& 	\bar{W} = f (M)^2 d^2 \left(1+  \left(\frac{\kappa (p)}{3} \pm \frac{2 \bar R}{3 \sqrt{1-\bar R^2}}\right)d\right)+ O (d^4),
\end{split}
\]
where in the second expansion we take the positive sign if $\overline{\tau} (\bar \U) \geq \tau_f^0(M)$ and the negative one otherwise. Now the result follows as in the proof of Proposition \ref{prop_curvatureZeroSets}.
\end{proof}

As a consequence of the curvature estimates obtained above, we can prove the property about the $\overline{\tau}$-function that was announced at the end of Subsection \ref{sectionTau}. 

\begin{proposition}\label{prop_boundTau}
Let $(\Omega, u)$, $\Omega \subset \s ^2$ a finite domain, be an analytic solution to \eqref{OEP}, where we assume that either $f\equiv 1$ or $f\in\mathcal{C}^{\omega}$ satisfies $f(x)\geq x f'(x)$ and $f'(x)\geq 2$ for all $x>0$. Fix a connected component $\U\in \pi_0\big(\Omega\setminus \textup{Max}(u)\big)$. If $\textup{cl}(\U)\cap \textup{Max}(u)$ contains a closed curve $\gamma$, then $\overline{\tau}(\U)>1$.
\end{proposition}
\begin{proof}
Argue by contradiction and assume that $\overline{\tau}(\U)\le 1$. Set $M:=u_{\textup{max}}$. Then, for every $R\in[0,1)$ the pair $(\U_{R,M,f}^{-},u_{R,M,f})$ is an admissible comparison pair for $(\U,u)$. As discussed above, any closed curve $\gamma\subset \textup{cl}(\U)\cap \textup{Max}(u)$ is an analytic simple closed curve and an isolated component of $\textup{Max}(u)$.

Fix a regular point $p\in\gamma$. Proposition~\ref{prop_CurvatureMaxSets} yields
\[
|\kappa(p)| \ge \frac{R}{\sqrt{1-R^2}} \text{ for each } R\in[0,1).
\]
Since the right-hand side can be made arbitrarily large by choosing $R$ sufficiently close to $1$, this forces $|\kappa(p)|=+\infty$, a contradiction. Therefore $\overline{\tau}(\U)>1$.
\end{proof}

\paragraph{Length estimates.}
To advance our analysis, we derive bounds on the lengths of the level sets of $u$ in terms of a comparison pair.

\begin{proposition}\label{prop_lengthBounds}
Let $(\Omega, u)$, $\Omega \subset \s ^2$ a finite domain, be an analytic solution to \eqref{OEP}, where we assume that either $f\equiv 1$ or $f\in\mathcal{C}^{\omega}$ satisfies $f(x)\geq x f'(x)$ and $f'(x)\geq 2$ for all $x>0$. Set $M:=u_{\textup{max}}$. Fix a connected component $\U\in \pi_0\big(\Omega\setminus \textup{Max}(u)\big)$ and let $(\bar\U,\bar u)$ be a comparison pair associated to $(\U,u)$ with critical height $\bar R<1$. 

Set $\gamma^{\U} := \textup{cl}(\U) \cap \textup{Max}(u)$ and $\Gamma^{\U} := \textup{cl}(\U) \cap \partial \Omega $. Then 
	\begin{equation}\label{lengthBound}
	\abs{\gamma^\U} \leq \sqrt{\frac{1- \bar{R}^2}{1-\bar r_+^2}}\abs{\Gamma^\U} \quad \textup{if } \overline{\tau} (\bar{\mathcal{U}}) \geq \tau_f^0(M) \quad \textup{and} \quad 	\abs{\gamma^\U} \leq \sqrt{\frac{1- \bar{R}^2}{1-\bar r_-^2}}\abs{\Gamma^\U} \quad \textup{if } \overline{\tau} (\bar{\mathcal{U}}) < \tau_f^0(M),
	\end{equation}
	where $\abs{\gamma^\U}$ and $\abs{\Gamma^\U}$ denote the sum of the lengths of each connected component. 
\end{proposition}

\begin{proof}
The argument follows \cite[Proposition~4.3]{EMa}. By Remark~\ref{RemCond}, the sets $\gamma^{\U}$ and $\Gamma^{\U}$ consist of analytic simple closed curves, possibly together with isolated points. If $\gamma^{\U}$ only contains isolated points then the bounds in \eqref{lengthBound} are trivial, so we assume that $\gamma^\U$ contains at least a curve. For $\epsilon>0$, define
\[
\U_{\epsilon} = \set{p \in \U ~\colon~M-u(p)>\epsilon}.
\]
Choose $\epsilon>0$ such that $M-\epsilon$ is a regular value of $u$. Since $u$ is analytic, Sard's Theorem implies that the critical values of $u$ are isolated; hence such choices exist with $\epsilon$ arbitrarily small.

For fixed $\epsilon$, the set $\gamma_{\epsilon}^{\U}:=\partial \U_\epsilon \setminus \Gamma^\U$ is a finite union of analytic curves, and we can apply the divergence theorem on $\U_\epsilon$. Let $\nu$ be the inner unit normal to $\U_\epsilon$. Using that $\bar U$ is the solution to the O.D.E.~\eqref{ODE} defining the comparison function $\bar u$ and that $\Psi$ is the pseudo-radial function from Definition~\ref{psudoRadialFunctions}, we obtain
\[
\int_{\U_\epsilon} \frac{-f(u)}{(1-\Psi^2) \bar{U} ' (\Psi)} \, d A
= -\int_{\U_\epsilon} \frac{f(u) \abs{\nabla u}^2}{(1-\Psi^2)^{2} \bar{U}' (\Psi)^3}  \, dA
+\int_{\partial \U_\epsilon} \frac{\pscalar{\nabla u}{\nu}}{(1-\Psi^2)\bar{U}' (\Psi)} \, d s.
\]
Along $\Gamma^\U$ we have $\nu=\nabla u/\abs{\nabla u}$, whereas along $\gamma_\epsilon^{\U}$ we have $\nu=-\nabla u/\abs{\nabla u}$. Rearranging terms yields
\begin{equation}\label{integralIdentity}
	\int_{\U_\epsilon} \frac{f(u)}{(1-\Psi^2)^{2} \bar{U}' (\Psi)^3} (\abs{\nabla u}^2-\bar{W}) \, dA
	= \int_{\Gamma^\U} \frac{-\abs{\nabla u}}{(1-\Psi^2)\bar{U}' (\Psi)} \, d s
	+\int_{\gamma_\epsilon} \frac{\abs{\nabla u}}{(1-\Psi^2)\bar{U}' (\Psi)} \, d s,
\end{equation}
where $\bar W=(1-\Psi^2)\bar U'(\Psi)^2$ is as in Section~\ref{sectionGradient}.

We now estimate the three terms in \eqref{integralIdentity}. Since $\Psi=\bar r_\pm(=:r_\pm(\bar R,M,f))$ along $\Gamma^\U$, we obtain
\begin{equation}\label{1}
	\int_{\Gamma^\U} \frac{\abs{\nabla u }}{(1-\Psi^2) \bar{U}' (\Psi)} \, ds
	= \mp \frac{1}{\sqrt{1-\bar{r}_\pm^2}} \int_{\Gamma^\U} \frac{\abs{\nabla u}}{\sqrt{\bar W}} \, ds
	\left\{	\begin{matrix}
		\geq - \frac{\abs{\Gamma^\U}}{\sqrt{1-\bar{r}_+^2}} & \text{ if } & \overline{\tau} (\bar{\mathcal{U}}) \geq \tau_f^0 (M) ,\\[2mm]
		\leq \frac{\abs{\Gamma^\U}}{\sqrt{1-\bar{r}_-^2}} & \text{ if } & \overline{\tau} (\bar{\mathcal{U}}) < \tau_f^0 (M),
	\end{matrix}\right.
\end{equation}
where the last inequalities use $\abs{\nabla u}^2 \le \bar W$ from Theorem~\ref{theor_Gradient}.

Next, take a sequence $\epsilon_n\to 0$ such that $M-\epsilon_n$ is a regular value for all $n$. Using the expansions \eqref{RhoExpansionW} and \eqref{RhoExpansionModelW}, we have
\[
\lim_{x \in \U, x \to \textup{Max}(u)} \frac{\abs{\nabla u}^2}{\bar W} =1,
\]
and therefore
\begin{equation}\label{2}
	\lim_{n \to +\infty} \int_{\gamma_{\epsilon_n}}  \frac{\abs{\nabla u}}{(1-\Psi^2)\bar{U}' (\Psi)} \, d s=
	\left\{	\begin{matrix}
		= - \frac{\abs{\gamma^\U}}{\sqrt{1-\bar{R}^2}} & \text{ if } & \overline{\tau} (\bar{\mathcal{U}}) \geq \tau_f^0 (M) ,\\[2mm]
		= \frac{\abs{\gamma^\U}}{\sqrt{1-\bar{R}^2}} & \text{ if } & \overline{\tau} (\bar{\mathcal{U}}) < \tau_f^0 (M).
	\end{matrix}\right.
\end{equation}

Finally, note that $\bar U'(\Psi)^3<0$ if $\overline{\tau}(\bar\U)\ge \tau_f^0(M)$ and $\bar U'(\Psi)^3>0$ if $\overline{\tau}(\bar\U)<\tau_f^0(M)$. Since $\abs{\nabla u}^2-\bar W\le 0$ in $\U$, it follows that
\begin{equation}\label{3}
	\int_{\U} \frac{f(u)}{(1-\Psi^2) \bar{U}' (\Psi)^3} (\abs{\nabla u}^2-\bar{W}) \, dA
	\left\{	\begin{matrix}
		\geq 0& \text{ if } & \overline{\tau} (\bar{\mathcal{U}}) \geq \tau_f^0 (M) ,\\[2mm]
		\leq 0 & \text{ if } & \overline{\tau} (\bar{\mathcal{U}}) < \tau_f^0 (M).
	\end{matrix}\right.
\end{equation}
Letting $\epsilon_n\to 0$ in \eqref{integralIdentity} and combining \eqref{1}, \eqref{2}, and \eqref{3}, we obtain \eqref{lengthBound}.
\end{proof}

The following estimate relies on Proposition \ref{prop_curvatureZeroSets} and on integration along a zero level set, and therefore it is valid only for solutions of the overdetermined problem \eqref{OEP2}. Since its proof follows exactly the same arguments as in \cite[Proposition 5.1]{EMa}, we omit it.

\begin{proposition}\label{propConvexity}
Let $(\Omega, u)$, $\Omega \subset \s ^2$ a finite domain, be an analytic solution to \eqref{OEP}, where we assume that either $f\equiv 1$ or $f\in\mathcal{C}^{\omega}$ satisfies $f(x)\geq x f'(x)$ and $f'(x)\geq 2$ for all $x>0$. Set $M:=u_{\textup{max}}$. Fix a connected component $\U\in \pi_0\big(\Omega\setminus \textup{Max}(u)\big)$ and let $(\bar\U,\bar u)$ be a comparison pair associated to $(\U,u)$. If $\Gamma \in \pi_0 (\partial \Omega \cap \textup{cl} (\U))$, then
	\begin{equation*}
		\abs{\Gamma} \leq 2 \pi \sqrt{1-\bar{r}_+^2} \quad \textup{if } \overline{\tau} (\mathcal{U}) \geq \tau_f^0(M).
	\end{equation*}
If, in addition, $\bar r_-<0$, then
	\[
	\abs{\Gamma} \leq 2 \pi \sqrt{1-\bar{r}_-^2}\quad \textup{if } \overline{\tau} (\mathcal{U}) < \tau_f^0(M).
	\]
\end{proposition}

\subsection{Proof of Theorem C}

Now, we are ready to show:
\\
\\
{\bf Theorem C:} {\it Let $(\Omega, u)$, $\Omega \subset \s ^2$ a finite domain, be an analytic solution to \eqref{OEP}, where we assume that either $f\equiv 1$ or $f\in\mathcal{C}^{\omega}$ satisfies $f(x)\geq x f'(x)$ and $f'(x)\geq 2$ for all $x>0$.

Then, either $\textup{Max}(u)$ contains only finitely many isolated points or the solution $(\Omega ,u)$ is rotationally symmetric. In the latter case, $\Omega$ is an annulus.}

\begin{proof}[Proof of Theorem C]
Assume that $\textup{Max}(u)$ contains infinitely many points. By Remark~\ref{RemCond}, $\Max(u)$ contains a simple closed non--contractible curve $\gamma$ (here we use that $f>0$ on $\r_+^*$). The curve $\gamma$ separates the domain into two connected components,
\[
\Omega\setminus\gamma=\Omega_\gamma^-\cup\Omega_\gamma^+,
\]
and we set $\Gamma_\pm:=\partial\Omega\cap \textup{cl}(\Omega_\gamma^\pm)$.

By Corollary~A, $(\Omega,u)$ is either rotationally symmetric or antipodally symmetric. If $(\Omega,u)$ is rotationally symmetric, then $\Omega$ is an annulus and we are done. Hence, we may assume that $(\Omega,u)$ is antipodally symmetric. In particular, the antipodal map $\mathcal A$ exchanges the two components, namely $\mathcal A(\Omega_\gamma^+)=\Omega_\gamma^-$, and therefore
\[
\overline{\tau}(\Omega_\gamma^+)=\overline{\tau}(\Omega_\gamma^-).
\]
Fix a comparison pair $(\bar\U,\bar u)$ that is admissible for both $(\Omega_\gamma^+,u)$ and $(\Omega_\gamma^-,u)$. Applying Proposition~\ref{prop_CurvatureMaxSets} to $\gamma$ from either side yields that its geodesic curvature vanishes identically. Hence $\gamma$ is an equator. By the Cauchy--Kovalevskaya Theorem, it follows that $(\Omega,u)\equiv (\Omega_{0,M,f},u_{0,M,f})$.
\end{proof}


We collect here several remarks clarifying the role of the assumptions on the non-linearity $f$ in Theorem~C.

\begin{itemize}
\item Both assumptions \eqref{conditionf} and \eqref{conditionf2} ensure the existence of the two--parameter family of ring--shaped, rotationally symmetric solutions $(\Omega_{R,M,f},u_{R,M,f})$ to \eqref{OEP}, together with a precise control of their boundary derivatives; see Section~\ref{SectModel}. 
This is the input needed to associate to a general solution $u$ a comparison pair in terms of the maximum value $u_{\rm max}=M$ and the boundary Neumann data. 
We encode this correspondence in the $\overline{\tau}$-\textit{function} (Definition~\ref{DefinTau}), in analogy with \cite{ABM} (and \cite{EMa} on $\s^2$), where $\overline{\tau}$ is introduced as the \textit{Normalized Wall Shear Stress}. 
In \cite{ABM,EMa}, the analysis is driven by a suitable $P$--function, which also yields strong topological restrictions. 
For a general nonlinearity $f$, such a $P$--function is not currently available, so here $\overline{\tau}$ is defined geometrically, using the model solutions on geodesic disks. 
The lack of a $P$--function means that $\overline{\tau}$ alone does not control the topology of $\Omega$ as in \cite{ABM,EMa}; this is compensated in Section~\ref{sectionConsequences} by the curvature estimates derived from Theorem~\ref{theor_Gradient}.

\item On the other hand, assumption \eqref{conditionf2} is also essential in the comparison argument leading to Theorem~\ref{theor_Gradient}. 
More precisely, when comparing a model pair $(\U^\pm_{R,M,f},u_{R,M,f})$ with a general solution $(\Omega,u)$, we obtain gradient and curvature estimates for $u$ (Lemma~\ref{LemDeltaF} and Theorem~\ref{theor_Gradient}) whose validity is tied to the sign of the function $\mu_{R,M,f}$ (item~7 in Theorem~\ref{theor_ODE}). 
Condition \eqref{conditionf2} is precisely the assumption guaranteeing $\mu_{R,M,f}\geq 0$ for the relevant model solutions. 
In the settings of \cite{ABM,EMa} the model solutions are explicit, and $\mu_{R,M,f}$ can be computed exactly.
\end{itemize}

\appendix

\section{Technical results and auxiliary estimates}\label{AppendixTech}

This appendix collects several technical results and auxiliary arguments used throughout Section~\ref{SectAM}. They are gathered here to streamline the presentation of the main proofs and to emphasize the structural role they play in the comparison arguments leading to Theorem~C. 

\subsection{Proof of Theorem \ref{theor_ODE}}\label{Appendix_theor_ODE}

Let $R \in [-1,1]$ and let $M \in \r_+^*$. We denote by $U_{R,M,f}$ the solution to the ODE with Cauchy data \eqref{ODE}-\eqref{CauchyData}. To prove this theorem, it will be useful sometimes to consider the change in the height parameter $r= \cos (s)$, from the cylindrical coordinate $r$ to the geographical coordinate $s$. Hence, in the coordinate $s$, the equation \eqref{ODE} reads as
\begin{equation}\label{ODERadial}
	V'' (s) + \cot (s) V' (s) + f(V(s)) =0,
\end{equation}
and then, by imposing the Cauchy data 
\begin{equation}\label{Cauchy2}
	V(S)=M, \quad V' (S)=0,
\end{equation}
for some $S \in [0, \pi], M \in \r_+^*$, we get that a solution $V_{S,M,f}$ to \eqref{ODERadial}-\eqref{Cauchy2} satisfies the relation $V_{S,M,f} (s) = U_{\cos (S),M,f} (\cos (s))$.

Thus, item 1 in the case $R=1$ follows from the existence and uniqueness of a solution to \eqref{ODERadial}-\eqref{Cauchy2} when $S=0$, which is proved in \cite[Lemma 3.5]{EM}. Reproducing the proof of this lemma in the case $S= \pi$, one obtains item~1 when $R=-1$.  When $R \in (-1,1)$, classical theory of differential equations implies the existence of a unique solution $U_{R,M,f}$ of $\eqref{ODE}$ satisfying \eqref{CauchyData} defined over an open interval containing $R$. The symmetry $U_{R,M,f} (-r) = U_{-R,M,f} (r)$ now follows from the fact that $U_{R,M,f} (-r)$ solves \eqref{ODE} and $U_{R,M,f}' (-R)=0, U_{R,M,f} (-R)=M$.

Item $2$ follows the ideas of \cite[Lemma 3.6]{EM} using that we can obtain an implicit representation of the solutions as
\begin{equation}\label{solutions}
U_{R,M,f}(r)
= M-\int_{R}^{r}\frac{1}{1-y^2}\left(\int_{R}^{y} f\bigl(U_{R,M,f}(x)\bigr)\,dx\right)dy.
\end{equation}
Differentiating \eqref{solutions} with respect to $r$ yields
\begin{equation*}
U'_{R,M,f}(r)
= -\frac{1}{1-r^2}\int_{R}^{r} f\bigl(U_{R,M,f}(x)\bigr)\,dx.
\end{equation*}
In particular, we deduce that $U'_{R,M,f} (r)>0$ if $r<R$ and $U'_{R,M,f} (r)<0$ if $r>R$. 

Item 3 follows from \cite[Lemma 3.6]{EM}, taking into account that $V_{0,M,f} (s) = U_{1,M,f} (\cos (s))$.

The first part of item $4$ follows from the Implicit Function Theorem. To prove the monotonicity properties of $r _ \pm$ we use the arguments employed in \cite{EM}. Fix $M \in \r_+^*$ and define the function
\begin{equation}\label{functionHR}
	Z_R (r)= \frac{\partial }{\partial R} \left(U_{R,M,f}\right) (r), \quad r \in [r_-, r_+].
\end{equation}

By a direct computation, $Z_R$ is a solution to the differential equation
\begin{equation*}
	(1-r^2)Z''(r)- 2 r Z' (r)+ f' (U_{R,M,f} (r)) Z (r)=0, 
\end{equation*}
and, if we write the rotationally symmetric function
\begin{equation*}
	z_R (r, \theta )= Z_R (r), \quad  (\sqrt{1-r^2} \cos \theta , \sqrt{1-r^2} \sin \theta ,r) \in \Omega_{R,M,f},
\end{equation*}
using cylindrical coordinates, $z_R$ lies in the kernel of the linearized operator, that is, $z_R$ is a solution to $\Delta z_R + f' (u_{R,M,f}) z_R =0$ in $\Omega_{R,M,f}$. 

\begin{quote}
{\bf Claim A:} {\it
For any $(R,M)\in [0,1)\times \mathcal \r_+^*$, the function $Z_R$ is negative on $(r_-(R,M,f),R)$ and positive on $(R,r_+(R,M,f))$.}
\end{quote}
\begin{proof}[Proof of Claim A]
We follow the proof of \cite[Lemma 3.8]{EM}. Let 
\begin{equation}\label{subannulus}
A_{r'}^{r''}= \set{(\sqrt{1-r^2} \cos \theta , \sqrt{1-r^2} \sin \theta ,r) \in \Omega_{R,M,f} ~\colon~ r \in (r',r'')}
\end{equation}for each $r_- \leq r' < r'' \leq r_+$, we view $A_{r'}^{r''}$ as a subannulus of $\Omega_{R,M,f}$.
	
Computing the derivative of \eqref{solutions} with respect to $R$, we get the implicit relation
\[
Z_R (r)=f(M) \int_{R}^r \frac{1}{1-x^2} \, dx - \int_{R}^r \frac{1}{1-y^2} \left( \int_R^y f'(U_{R,M,f} (x)) Z_R (x) \, dx \right) \, dy,
\]	
Evaluating at $r=R$ yields $Z_R(R)=0$. More explicitly, differentiating the identities
$U_{R,M,f}(R)=M$ and $U'_{R,M,f}(R)=0$ with respect to $R$ gives
\[
Z_R(R)+U'_{R,M,f}(R)=0,
\qquad
Z_R'(R)+U''_{R,M,f}(R)=0.
\]
Since $U'_{R,M,f}(R)=0$ and $(1-R^2)U''_{R,M,f}(R)+f(M)=0$ by \eqref{ODE}, we obtain
$Z_R(R)=0$ and $Z_R'(R)=\frac{f(M)}{1-R^2}>0$. It then follows that there exists $r'<R<r''$ such that $Z_R$ is negative in $(r',R)$ and positive in $(R,r'')$. We first prove that, in fact, $Z_R >0$ in $(R, r_+)$. We argue by contradiction. If this were not the case, then there exists a first $r_0 \in (R, r_+)$ such that $Z_R (r_0)=0$. Define the quantity $\delta= \min_{\textup{cl}(A_R^{r_0})} \frac{u_{R,M,f}}{z_R}>0$, and consider the function $w=u_{R,M,f}- \delta z_{R}$. Then $w$ is a nonnegative function inside $A_R^{r_{0}}$, and since $z_R$ is a solution to $\Delta z_R + f' (u_{R,M,f}) z_R =0$, using \eqref{conditionf} yields
	\begin{equation*}
	\Delta w + f' (u_{R,M,f}) w = -f(u_{R,M,f})+ f'(u_{R,M,f})u_{R,M,f} \leq 0.
	\end{equation*}
But $w$ attains zero at some point inside $A_R^{r_{0}}$, so the strong maximum principle implies that $w \equiv 0$ in $A_R^{r_0}$. This is a contradiction, since $U_{R,M,f}(r_0)-\delta Z_R (r_0)=U_{R,M,f}(r_0)>0$. Finally, the sign of $Z_R$ on $(r_-,R)$ is proved in the same way, working on a subannulus to the left of $R$ and applying the above argument to $-z_R$. This completes the proof of Claim~A.
\end{proof}

Now, since $U_{R,M,f} (r_\pm) =0$ for all $R \in [0,1)$, by the Implicit Function Theorem it follows that 
	\[
	\frac{\partial r_\pm}{ \partial R} (R,M) = -\frac{Z_R}{U_{R,M,f}'} (r_\pm) \geq 0, \quad \forall R \in [0,1), 
	\]which proves item $4$.

To prove item $5$, it is better to use geographical coordinates in $\s^2$. Define the functions $s_{\pm}:= s_{\pm} (S,M,f) = \arccos (r_{\mp} (\cos (S),M,f))$, for all $S \in [0, \pi)$. Then, recalling \eqref{functionNormOfGradient}, we have that
\[
g_{M,f}^{\mp} (\cos (S)) = V_{S,M,f} ' (s_{\pm} (S,M,f))^2, \quad \forall S \in (0, \pi / 2],
\]
Defining $\tilde{Z}_S = \frac{\partial}{\partial S} (V_{S,M,f})$, we note that
\begin{equation}\label{eq:derivativeVS}
	\begin{split}
		\frac{\partial }{\partial S} \left( V_{S,M,f} ' (s_{\pm})^2 \right) &= 2 \frac{\partial V_{S,M,f} '}{\partial S} V_{S,M,f}' (s_{\pm}) + 2 V_{S,M,f} ' (s_{\pm}) V_{S,M,f} '' (s_\pm) \frac{\partial s_{\pm}}{ \partial S} \\
		&=-2 \left( \tilde{Z}_S (s_{\pm}) V_{S,M,f} '' (s_{\pm})- \tilde{Z}_S ' (s_{\pm}) V_{S,M,f} ' (s_{\pm}) \right),
	\end{split}
\end{equation}
where we have used that $\partial s_{\pm} / \partial S = -(\tilde{Z}_S / V_{S,M,f} ') (s_{\pm})$.

Thus, item $5$ will follow as a corollary of the following:
\begin{quote}
	{\bf Claim B:} \it The function
	$$G_S (s) := V_{S,M,f}''(s)\tilde{Z}_S (s)-V_{S,M,f}'(s)\tilde{Z}_S '(s), \quad s \in [s_-,s_+],$$is positive in $[s_-,S)$ and negative in $(S,s_+]$ for any $(S,M) \in (0, \pi / 2] \times \mathcal I _f$.
\end{quote}
\begin{proof}[Proof of Claim B]

Now let us fix some notation about Killing vector fields in $\s^2$. For each $p \in \s^2$, $v \in T_p \s^2$ with $\abs{v}^2 =1$, let $\gamma_{p,v}: \r \to \s^2$ be the unique geodesic parametrized by arc-length with $\gamma_{p,v} (0)=p$ and $\gamma_{p,v}' (0)=v$. Then, recalling that $\textbf{n}=(0,0,1)$ and writing $Y=(1,0 , 0) \in T_{\textbf{n}} \s^2$, the Killing vector field with flow of isometries that fix the trace of the geodesic $\gamma_{\textbf{n},Y}$ is given by
\begin{equation}\label{killingField}
	\mathcal{Y} (q)= -\pscalar{q}{Y} \textbf{n}+ \pscalar{\textbf{n}}{q} Y, \quad \forall q \in \s^2.
\end{equation}
We follow the arguments contained in the proof of \cite[Lemma 3.10]{EM}. Thus, we take $S<s_0 \leq s_+$ and we will show that $G_S (s_0)<0$.

Define $\tilde{z}_S (s,\theta) = \tilde{Z}_S (s)$ on $\s^2$. Note that $\tilde{Z}_S (s) = -\sin (S) Z_{\cos (S)} (\cos (s))$, where $Z_R$ is defined in \eqref{functionHR}, so Claim A implies that $\tilde{Z}_S$ is negative in $(s_- , S)$ and positive in $(S, s_+)$. Hence, if $\tilde{Z}_S (s_0)=0$ then it must be $\tilde{Z}_S ' (s_0)<0$ by the maximum principle (since $\Delta \tilde{z}_S + f' (u_{\cos (S),M,f}) \tilde{z}_S =0$), so it follows that $G_S (s_0)=-V_{S,M,f}' (s_0) \tilde{Z}_S ' (s_0)<0$. Then, suppose that $\tilde{Z}_S (s_0)>0$ and consider the function $v$ defined by
	\[
	v=\frac{V_{S,M,f}' (s_0)}{\tilde{Z}_S (s_0)} \tilde{z}_{S} - \meta{\mathcal{Y}}{\nabla u_{\cos (S),M,f}}, \quad \textup{in } A_{\cos (S)}^{\cos (s_0)},
	\]
	where $\mathcal{Y}$ is the Killing vector field defined in \eqref{killingField} and $A_{\cos (S)}^{\cos (s_0)}$ is defined as in \eqref{subannulus}. In geographical coordinates, we have that 
	\[
	v(s,\theta)= \frac{V_{S,M,f}' (s_0)}{\tilde{Z}_S (s_0)} \tilde{Z}_{S}(s)-V_{S,M,f} ' (s) \cos \theta,
	\]
	and since $V_{S,M,f}' <0$ in $(S,s_+]$, we conclude that $\max_{\partial A_{\cos (S)}^{\cos (s_0)}} v = 0$, and this maximum is reached at 
	\[
	(\sin (s) \cos \theta , \sin (s) \sin \theta ,\cos (s))=\left(\sin (s_0) , 0 , \cos (s_0) \right)=: q_0. 
	\]
	Now we prove that $v \leq 0$ in $A_{\cos (S)}^{\cos (s_0)}$. To do so, consider the open set $$\mathcal{V}= \set{q \in A_{\cos (S)}^{\cos (s_0)} ~\colon~v(q)>0},$$ and let us suppose that $\mathcal{V} \neq \emptyset$. Then the value $\delta=\min_{\textup{cl}(\mathcal{V})} \frac{\tilde{z}_S}{v}$ is reached inside it, so the function $w=\tilde{z}_S-\delta v \geq 0$ in $\mathcal{V}$ and $w(p)=0$ for some point $p \in \mathcal{V}$. But since $\mathcal{Y}$ is a Killing field, it follows that $\Delta w + f'(u_{\cos (S),M,f}) w=0$ over $\mathcal{V}$, so by the maximum principle it must be $w \equiv 0$ in $\mathcal{V}$. But this is impossible since $w >0$ along $\partial \mathcal{V}$, so we get a contradiction, and it must be $\mathcal{V} = \emptyset$.\\
	We conclude that it must be $v \leq 0$ in $A_{\cos (S)}^{\cos (s_0)}$, so, as we have that $\Delta v + f' (u_{\cos (S),M,f}) v =0$ in $A_{\cos (S)}^{\cos (s_0)}$, by the maximum principle it must be $v \equiv 0$ inside $A_{\cos (S)}^{\cos (s_0)}$ or $\frac{\partial v}{\partial \nu} (q_0) <0$, where $\nu$ is the inner conormal to $\partial A_{\cos (S)}^{\cos (s_0)}$. However, the first case is ruled out because $v \neq 0$ along $\partial A_{\cos (S)}^{\cos (s_0)}$, so we have that
	\[
	\frac{\partial v}{\partial \nu} (q_0) = V_{S,M,f}'' (s_0)-\frac{V_{S,M,f} ' (s_0)}{\tilde{Z}_S (s_0)} \tilde{Z}_S ' (s_0) <0,
	\]
	and from here it follows that it must be $G_S (s_0) <0$. To conclude, we note that to prove that $G_S>0$ in $[s_-,S)$ we can apply the same argument as above but with $-v$ instead of $v$. This proves Claim B.
\end{proof}
Recall that $R=\cos S$ and
$s_{\pm}(S,M,f)=\arccos\bigl(r_{\mp}(\cos S,M,f)\bigr)$, so that
\[
g_{M,f}^{\mp}(\cos S)=V_{S,M,f}'\bigl(s_{\pm}(S,M,f)\bigr)^2.
\]
Differentiating with respect to $S$ gives \eqref{eq:derivativeVS}, namely
\[
\frac{\partial}{\partial S}\Bigl(V_{S,M,f}'\bigl(s_{\pm}\bigr)^2\Bigr)=-2\,G_S\bigl(s_{\pm}\bigr).
\]
Since $R=\cos S$, we have $dS/dR=-1/\sqrt{1-R^2}$. Therefore, after rewriting
$S=\arccos R$ and $s_{\pm}=\arccos(r_{\mp}(R,M,f))$, we obtain the equivalent form
\[
(g_{M,f}^\pm)'(R)
= \frac{2}{\sqrt{1-R^2}}\,
G_{\arccos(R)}\bigl(\arccos(r_\pm(R,M,f))\bigr),
\qquad \forall\, R\in[0,1).
\]
In particular, the sign information in Claim~B yields the monotonicity asserted in Item~5.
For item $6$, we note that
\[
h_f(M) = -V_{0,M,f}' (s_+ (0,M,f)), \quad \forall M \in \r_+^*,
\]
so we conclude that the above function is increasing after \cite[Subsection 3.2.3]{EM}.

Let us now show item $7$. When $f(x) =1$, we can explicitly solve the first derivative of a solution to \eqref{ODE}-\eqref{CauchyData} as 
$$U_{R,M,f}' (r) = - \frac{r-R}{1-r^2}, \quad r \in [r_- , r_+] \text{ and } (R,M) \in [0,1) \times \r_+^*.$$

A straightforward computation shows that $ \mu_{R,M,f} (r) = 2\frac{1-R^2}{(1-r^2)^3} >0$.

Let us focus on the case $f' (x)\geq 2$. Since $U_{R,M,f}$ is decreasing on $(R,r_+]$ and increasing on $[r_-,R)$, isolating $U_{R,M,f}''$ in \eqref{ODE} we obtain that $U_{R,M,f}$ is concave on $[R,r_+]$. Furthermore, if $r_-<0$, then $U_{R,M,f}$ is concave on $[r_-,0]$. Note that since $f$ is $\mathcal{C}^1$ we have that $U_{R,M,f}$ is of class $\mathcal{C}^3$ and \eqref{ODE} implies that $U_{R,M,f}$ solves the third-order equation
\begin{equation}\label{ThirdDerivative}
	(1-r^2)U''' _{R,M,f}(r)-4rU'' _{R,M,f}(r)+(f'(U_{R,M,f}(r))-2) U'_{R,M,f} (r)=0, 
\end{equation}which yields
$$ \mu_{R,M,f} (r) = \frac{1}{1-r^2}\left( (f'(U_{R,M,f}(r))-2) U'_{R,M,f} (r)^2 -2 f(U_{R,M,f} (r))U''_{R,M,f} (r)\right)$$which is non-negative on $[R,r_+]$ and $[r_-,0]$ if $r_- <0 $ by the assumption $f ' (x)\geq 2 $ and the concavity of $U_{R,M,f}$. 

To show that $\mu_{R,M,f}$ is non-negative on $[{\rm max}\set{0,r_-} ,R]$, it is enough to prove that $U'''_{R,M} \leq 0$ on this interval. Integrate \eqref{ThirdDerivative} between $r_0 :={\rm max}\set{0,r_-} \geq 0$ and $r \in (r_0 ,R)$ to obtain the following identity for the second derivative of $U_{R,M,f}$:
	\[
	U_{R,M,f}'' (r) = -f(U_{R,M,f} (r_0))-\frac{1}{(1-r^2)^{2}} \int_{r_0}^r (1-x^2) \left(f' (U_{R,M,f} (x))-2\right) U_{R,M,f}' (x) \, dx , \quad \forall r \in [0,R].
	\]
In order to check that $U_{R,M,f}''' \leq 0$ on $[0,R]$ we shall prove that $U_{R,M,f}''$ is non-increasing on this interval. Take $r_0 \leq a < b \leq R$, and write $\mathcal H(x)=(1-x^2) \left(f' (U_{R,M,f} (x))-2\right)U_{R,M,f}' (x)$ for $x \in [r_0,R]$. Note that, assuming \eqref{conditionf2}, it follows that $\mathcal H \geq 0$ on $[r_0,R]$. Hence, we have that
	\[
	\begin{split}
			U_{R,M,f}'' (a)-U_{R,M,f}'' (b) &= -\frac{1}{(1-a^2)^2}  \int_{r_0}^a \mathcal H(x) \, dx + \frac{1}{(1-b^2)^2}  \int_{r_0}^b \mathcal H(x) \, dx \\
			&= \left( \frac{1}{(1-b^2)^2}-\frac{1}{(1-a^2)^2} \right) \int_{r_0}^a \mathcal H(x) \, dx + \frac{1}{(1-b^2)^2}  \int_{a}^b \mathcal H(x) \, dx  \geq 0.
	\end{split}
	\]
	Then we conclude that $U_{R,M,f}'' (a) \geq U_{R,M,f}'' (b)$, and this proves Item~7.

\subsection{Proof of Lemma \ref{lemmaTau}}\label{Appendix_lemmaTau}

By item 5 of Theorem \ref{theor_ODE} the functions $\overline{\tau}_{M,f}^\pm$ have the claimed monotonicity properties. So, to finish, we shall study the image of each function. 

Observe that $g_{M,f}^- (1)=h_f(M)^2$ by item 6 of Theorem \ref{theor_ODE} and $g_{M,f}^-  (0) = g_{M,f}^+ (0)$ because of the symmetry of $U_{0,M,f} (r)$, item 1 of Theorem \ref{theor_ODE}. Then, from \eqref{tauModelSolutions}, we obtain that $\overline{\tau} _{M,f}^-([0,1)) = (1, \tau _f^0 (M)]$ and $\overline{\tau} _{M,f}^+ ([0,1)) = [\tau _0 (M) , \tau_f^\infty (M))$, for some $\tau_f^\infty (M) \in (\tau_f^0 (M), + \infty]$. It remains to show that $\tau_f^\infty (M) = + \infty$ for all $M \in \r_+^*$, which is the content of the next claim:
\begin{quote}
	{\bf Claim C:} {\it We have that
		$$\displaystyle\lim_{R \to 1^-} \abs{\nabla u_{R,M,f}} (r_+(R,M,f),\theta ) =+\infty, \quad \forall \theta \in \s^{1}.$$}
\end{quote}
\begin{proof}[Proof of Claim C]
	In the case $f(x)=1$ the solution is explicit, and the above limit can be checked easily. Thus, we restrict to the case in which $f' \geq 2$. We argue by contradiction. Assume therefore that $c_f(M) := h_f(M)\,\tau_f^{\infty}(M) < +\infty.$

	Note that conditions \eqref{conditionf} and \eqref{conditionf2} imply that $f(x) \geq 2x$ for all $x \in \r_+^*$. Let us define $\tilde f(x)=2x$ and denote by $(\U_{R,M,\tilde{f}}^+,u_{R,M,\tilde{f}})$
	the upper subdomains of the corresponding model solutions associated with the eigenvalue problem. We write $ U_{R,M,\tilde f}$ for the solution of \eqref{ODE}--\eqref{CauchyData} defining $u_{R,M,\tilde f}$, and we denote by $r_+(R,M,\tilde f)$ the first zero of $ U_{R,M,\tilde f}$ after $R$. By \cite[Lemma~3.2]{EMa}, there exists $R_0<1$ such that
	\[
	\sqrt{1- r_+(R,M,\tilde f)^2}\,
	\bigl| U_{R,M,\tilde f}'\bigl( r_+(R,M,\tilde f)\bigr)\bigr|
	> c_f(M),
	\qquad \forall\, R\in (R_0,1).
	\]
	Fix $R_1\geq  r_+(R,M,\tilde f)$. For a fixed $M\in \r_+^*$ and for each
	$R\in [R_0,1)$, define the function
	\[
	\omega_R :=  u_{R,M,\tilde f}-u_{R_1,M,f}
	\quad \text{in } \mathcal O_R := \U_{R_1,M,f}^+ \cap \U_{R,M, \tilde f}^+ .
	\]
	For $R$ sufficiently close to $R_0$, we have $\omega_R<0$ in $\mathcal O_R$. Increasing $R$ slightly, we reach a first value for which $\omega_R \leq 0$ in $\mathcal O_R$ and there exists a point $x\in \mathrm{cl}(\mathcal O_R)$ such that
	$\omega_R(x)=0$.
	
	On the other hand, in $\mathcal O_R$ we compute
	\[
	\Delta \omega_R
	= -2\, u_{R,M,\tilde f} + f(u_{R_1,M,f})
	\geq -2 \omega_R,
	\]
	where the inequality follows from $f(x)\geq 2x$. Therefore, by the strong interior or boundary maximum principles, we conclude that $	\omega_R \equiv 0$ in $\mathcal{O}_R$,	which is impossible.
	
	This contradiction proves Claim C.
\end{proof}

Hence, we complete the proof of Lemma~\ref{lemmaTau}.

\subsection{Proof of Lemma \ref{LemDeltaF}}\label{Appendix_LemDeltaF}

We follow the computations outlined in \cite[Section 4.1]{EMa}. We begin by computing $\Delta (\abs{\nabla u}^2-\bar{W})$. Using Bochner's identity, we get that
\[
\Delta \abs{\nabla u}^2 = 2 \left( 1-f' (u) \right) \abs{\nabla u}^2 + 2 \abs{\nabla^2 u}^2.
\]
On the other hand, we have that
\[
\Delta \bar W = \frac{\partial \bar W }{\partial u} \Delta u + \frac{\partial^2 \bar W}{\partial u^2} \abs{\nabla u}^2.
\]
Then, using the definition of the function $\chi_{\pm} = \chi$ (see Subsection \ref{sectionPseudoRadial}), it is easy to get that
\begin{equation}\label{firstDerivativeModelW}
	\frac{\partial}{\partial u}\left( (1-\chi (u) ^2) \bar{U}' (\chi (u))^2 \right)= 2 \Psi \bar{U}' (\chi (u)) -2 f(u)
\end{equation}
and
\begin{equation}\label{secondDerivativeModelW}
	\frac{\partial^2}{\partial u^2}\left( (1-\chi (u) ^2) \bar{U}' (\chi (u))^2 \right)= 2 \left( 1+ \frac{2 \chi (u)^2}{1-\chi (u)^2}-\frac{f (u) \chi (u)}{(1-\chi (u)^2) \bar{U}' (\chi (u))} \right)-2 f' (u),
\end{equation}
hence, since $\Psi (p) = \chi \circ u (p)$, it follows from \eqref{firstDerivativeModelW} and \eqref{secondDerivativeModelW} that
\[
\Delta \bar{W} =2 \left( \Psi \bar{U}' (\Psi) -2 f(u)\right)+ 2 \left( 1+ \frac{2 \Psi^2}{1-\Psi^2}-\frac{f (u) \Psi}{(1-\Psi^2) \bar{U}' (\Psi)}-f' (u) \right)\abs{\nabla u}^2.
\]
We conclude that
\begin{equation}\label{laplacianWs}
\Delta (\abs{\nabla u}^2-\bar{W}) =  2 \abs{\nabla^2 u}^2 - 2 \left( \Psi  \bar{U}' (\Psi) -2 f(u)\right)  -\frac{2}{1-\Psi^2} \left(2 \Psi^2 - \frac{f(u)\Psi}{\bar{U}' (\Psi)} \right)\abs{\nabla u}^2
\end{equation}
Next, we estimate the norm of the Hessian of $u$. To do so, we note that, for the model solution, the following relation is satisfied outside the set of critical points:
\[
\nabla^2 \bar{u} (r,\theta) = \frac{1}{2}\left( \varphi(r) \abs{\nabla u}^2 -f(\bar u) \right) g_{\sp^2}-\varphi(r) d\bar u \otimes d\bar u,\,  \text{ for }(r,\theta) \in \bar{\Omega} \setminus \textup{Max}(\bar u).
\]
Hence, we can obtain a bound for $\abs{\nabla^2 u}^2$ as it is done in \cite{ABM,EMa}, computing
\begin{equation}\label{normHessianRelation}
	\begin{split}
		0 & \leq \abs{\nabla^2 u + \varphi(\Psi) du \otimes du +  \frac{1}{2}\left( f(u)- \varphi(\Psi) \abs{\nabla u}^2 \right) g_{\sp^2}}^2 \\
		&= \abs{\nabla^2 u}^2+ 2\varphi (\Psi) \pscalar{\nabla^2 u}{du \otimes du}+\frac{1}{2} \varphi (\Psi)^2 \abs{\nabla u}^4-\frac{1}{2} f(u)^2 \\
		& \hspace{0.5cm} +\varphi (\Psi) f(u) \abs{\nabla u}^2.
	\end{split}
\end{equation}
Now observe that
\begin{equation}\label{gradientModelW}
	\begin{split}
		\nabla \bar W = 2 \left( \Psi \bar{U}' (\Psi)-f(u) \right) \nabla u 
		= - \left( \varphi (\Psi) \bar{W} + f(u) \right) \nabla u.
	\end{split}
\end{equation}
Hence, recalling that $\meta{\nabla \abs{\nabla u}^2}{\nabla u}=2 \pscalar{\nabla^2 u}{du \otimes du}$ and \eqref{gradientModelW}, we get
\[
\pscalar{\nabla^2 u}{du \otimes du}= \frac{1}{2} \pscalar{\nabla (\abs{\nabla u}^2-\bar{W})}{\nabla u}-\frac{1}{2} \left(  \varphi (\Psi) \bar{W} (\Psi) + f(u) \right) \abs{\nabla u}^2.
\]
Substituting this in \eqref{normHessianRelation}, we finally obtain
\begin{equation}\label{boundNormHessian}
\abs{\nabla^2 u}^2 \geq  - \varphi (\Psi)\pscalar{\nabla (\abs{\nabla u}^2-\bar{W})}{\nabla u}-\frac{1}{2} \varphi (\Psi)^2 \abs{\nabla u}^4 + 2 \frac{1}{2} \varphi (\Psi)^2 \bar{W} \abs{\nabla u}^2 +\frac{1}{2} f(u)^2.
\end{equation}
Then, noting that
\[
2 \left( \Psi  \bar{U}' (\Psi) -2 f(u)\right) = f(u) \left(  \varphi (\Psi) \bar{W} + f(u) \right)
\]
and
\[
-\frac{2 }{1-\Psi^2} \left(2 \Psi^2 - \frac{f(u)\Psi}{\bar{U}' (\Psi)} \right)\abs{\nabla u}^2 =  \left( f(u)-\varphi (\Psi) \bar{W} \right) \varphi (\Psi),
\]
we substitute \eqref{boundNormHessian} in \eqref{laplacianWs} and, after rearranging the terms, we get the differential inequality
\begin{equation}\label{laplacianWs2}
\Delta (\abs{\nabla u}^2-\bar{W}) \geq - 2 \varphi (\Psi) \pscalar{\nabla (\abs{\nabla u}^2-\bar{W})}{\nabla u}  +  \varphi (\Psi) \left( f(u)-\varphi (\Psi) \abs{\nabla u}^2 \right) (\abs{\nabla u}^2-\bar{W}).
\end{equation}
Consider now $F_\beta = \beta (\Psi) \cdot (\abs{\nabla u}^2-\bar{W})$ for a smooth function $\beta >0$. We follow the computations of \cite{EMa}. Observe that
\begin{equation}\label{laplacianF}
	\Delta F_{\beta} = \beta \cdot \Delta (\abs{\nabla u}^2-\bar{W}) + 2 \pscalar{\nabla \beta}{\nabla (\abs{\nabla u}^2-\bar{W})} +   (\abs{\nabla u}^2-\bar W) \Delta \beta,
\end{equation}
so we compute each term of the previous sum and then we substitute. It can be checked that
\begin{equation}\label{secondSum}
	2 \pscalar{\nabla \beta}{\nabla (\abs{\nabla u}^2-\bar{W})} = \frac{2}{\bar{U}' (\Psi)} \frac{\beta ' (\Psi)}{\beta} \pscalar{\nabla F_\beta}{\nabla u}-\frac{2}{\bar{U}' (\Psi)^2} \left(\frac{\beta ' (\Psi)}{\beta (\Psi)}\right)^2  \abs{\nabla u}^2 F_\beta
\end{equation}
and
\begin{equation}\label{thirdSum}
	(\abs{\nabla u}^2-\bar{W}) \Delta \beta = -\frac{f(u)}{\bar{U}' (\Psi)} \frac{\beta ' (\Psi)}{\beta (\Psi)} F_\beta + \frac{1}{\bar{U}' (\Psi)} \left( \varphi (\Psi) \frac{\beta ' (\Psi)}{\beta (\Psi)}+  \frac{1}{\bar{U}' (\Psi)} \frac{\beta '' (\Psi)}{\beta (\Psi)} \right) \abs{\nabla u}^2 F_\beta.
\end{equation}
Now, observe that
\[
\pscalar{\beta \nabla (\abs{\nabla u}^2-\bar{W})}{\nabla u}= \pscalar{\nabla F_\beta}{\nabla u}-\frac{1}{\bar{U}' (\Psi)}\frac{\beta ' (\Psi)}{\beta (\Psi)} F_\beta \abs{\nabla u}^2,
\]
so, taking this into account, after substituting \eqref{secondSum} and \eqref{thirdSum} into \eqref{laplacianF}, applying the bound \eqref{laplacianWs2} and rearranging the terms, we get
\begin{equation*}
	\begin{split}
		\Delta F_\beta \geq & \frac{2}{\bar{U}' (\Psi)} \left( \frac{\beta ' (\Psi)}{\beta (\Psi)}-\varphi (\Psi) \bar{U}' (\Psi) \right) \pscalar{\nabla F_\beta}{\nabla u} \\
		&-\frac{f(u)}{\bar{U}' (\Psi)} \left( \frac{\beta ' (\Psi)}{\beta (\Psi)}-  \varphi (\Psi) \bar{U}' (\Psi) \right) F_\beta \\
		&+\frac{\abs{\nabla u}^2}{\bar{U}' (\Psi)^2} \left( \left(  \frac{\beta ' (\Psi)}{\beta (\Psi)} \right)' - \left( \frac{\beta ' (\Psi)}{\beta (\Psi)}\right)^2 + 3 \varphi (\Psi) \bar{U}' (\Psi)  \frac{\beta ' (\Psi)}{\beta (\Psi)}- \varphi (\Psi)^2 \bar{U}' (\Psi)^2 \right) F_\beta.
	\end{split}
\end{equation*}
Now we substitute 
	\begin{equation*}
	\beta (\Psi) := \frac{1}{\abs{\bar{U}' (\Psi)}},
\end{equation*}
in the above expression, which is clearly an analytic function inside $\U$. Note that, in this case, we have that
	\[
\frac{\beta ' (\Psi)}{\beta (\Psi)}=  \varphi (\Psi) \bar{U}' (\Psi),
\]
and
\[
- \left( \frac{\beta ' (\Psi)}{\beta (\Psi)}\right)^2 + 3 \varphi (\Psi) \bar{U}' (\Psi)  \frac{\beta ' (\Psi)}{\beta (\Psi)}- \varphi (\Psi)^2 \bar{U}' (\Psi)^2 = \left(\frac{\beta ' (\Psi)}{\beta (\Psi)}\right)^2.
\]
Thus, we conclude that the function $F$ satisfies the elliptic inequality
\[
\Delta F - \left( \frac{2 \bar U'' (\Psi)^2- \bar U' (\Psi) \bar U ''' (\Psi)}{\bar U ' (\Psi)^2} \right) F \geq 0.
\]
To complete the proof, we evaluate
$$F = \abs{\bar{U}' (\Psi)}^{-1} (\abs{\nabla u}^2-\bar{W})= \frac{\abs{\nabla u}^2}{\abs{\bar{U}' (\Psi)}}-(1-\Psi^2) \abs{\bar{U}' (\Psi)} \text{ along }\partial \U.$$ 

Observe that, since $(\bar \U,\bar{u})$ is a comparison pair, it is clear that $\abs{\nabla u}^2(p) \leq \bar{W} (p)$ for each $p \in  \textup{cl} (\U) \cap \partial \Omega$ and hence $F \leq 0$ along $\textup{cl} (\U) \cap \partial \Omega$. On the other hand, Lemma \ref{lemma_limit} implies that
$$\displaystyle\lim_{p \in \U, p \to \textup{Max}(u)} \frac{\abs{\nabla u}^2}{\abs{\bar{U}' (\Psi)}} = \displaystyle\lim_{p \in \U, p \to \textup{Max}(u)} \frac{\sqrt{1-\Psi^2} \abs{\nabla u}^2}{ \sqrt{2 f(M)} \sqrt{M-u}}=0,$$where the last identity follows from \cite[Corollary 2.3]{Chr}. Thus, we conclude that 
$$\displaystyle\lim_{p \in \U, p \to \textup{Max}(u)} F (p) =0 .$$
Therefore, it follows that $F \leq 0$ along $\partial \U$, which proves the lemma.

\subsection{A Taylor expansion for $\abs{\nabla u}^2$ and $\bar W$ on $\textup{Max}(u)$}\label{AppTaylorMax}

In this subsection, let $(\Omega, u)$, $\Omega \subset \s ^2$ a finite domain, be an analytic solution to \eqref{OEP}, where we assume that either $f\equiv 1$ or $f\in\mathcal{C}^{\omega}$ satisfies $f(x)\geq x f'(x)$ and $f'(x)\geq 2$ for all $x>0$. Fix a connected component $\U\in \pi_0\big(\Omega\setminus \textup{Max}(u)\big)$ and let $(\bar\U,\bar u)$ be a comparison pair associated to $(\U,u)$.

Given the functions $\abs{\nabla u}^2$ and $\bar W$ defined in Subsection \ref{sectionGradient}, we obtain here a useful expansion for these functions near the set of maximum points $\textup{Max} (u)$. We will need first the following limit:
\begin{lemma}\label{lemma_limit}
	If $p \in \textup{Max}(u)$, $M= u_{\textup{max}}$, it holds:
	\[
	\displaystyle\lim_{x \in \mathcal{U},x \to p} \frac{\bar{W}}{M-u} = 2 f(M).
	\]	
\end{lemma}
\begin{proof}
This is a direct computation using L'H\^{o}pital's rule. Recall that $\bar U$ denotes the solution to
\eqref{ODE}--\eqref{CauchyData} defining the comparison pair in $\U$, and that
$u=\bar U\circ \Psi$ and $\bar W=(1-\Psi^2)\bar U'(\Psi)^2$ in $\U$. Hence,
\[
\begin{split}
\lim_{x\in \U,\,x\to p}\frac{\bar W}{M-u}
&=\lim_{\Psi\to \bar R}\frac{(1-\Psi^2)\bar U'(\Psi)^2}{M-\bar U(\Psi)} \\
&=\lim_{\Psi\to \bar R}\frac{-2\Psi \bar U'(\Psi)^2+2(1-\Psi^2)\bar U'(\Psi)\bar U''(\Psi)}{-\bar U'(\Psi)} \\
&=\lim_{\Psi\to \bar R}\Bigl(2\Psi \bar U'(\Psi)-2(1-\Psi^2)\bar U''(\Psi)\Bigr).
\end{split}
\]
Using \eqref{ODE} for $\bar U$, namely
$(1-\Psi^2)\bar U''(\Psi)=2\Psi \bar U'(\Psi)-f(\bar U(\Psi))$, we obtain
\[
\lim_{x\in \U,\,x\to p}\frac{\bar W}{M-u}
=\lim_{\Psi\to \bar R}\Bigl(-2\Psi \bar U'(\Psi)+2 f(\bar U(\Psi))\Bigr)
=2 f(M),
\]
since $\bar U(\bar R)=M$ and $\bar U'(\bar R)=0$.
\end{proof}
Since $u$ is real analytic in $\Omega$, the set $\textup{Max}(u)$ consists of finitely many isolated points and real-analytic curves. Suppose that $\textup{Max}(u)$ contains a curve $\gamma$, and fix a point $p \in \gamma$. If $\gamma$ is not closed, assume moreover that $p$ is not an endpoint of $\gamma$. Then there exists an open neighborhood $\mathcal V$ of $p$ such that $\mathcal V \setminus \gamma = \mathcal V_+ \cup \mathcal V_-$, where $\mathcal V_+$ and $\mathcal V_-$ are open, connected sets, with $\mathcal V_+ \subset \mathcal U$.

We define the signed distance function $d : \mathcal V \to \mathbb R$ by
\[
d(x) :=
\begin{cases}
	\phantom{-}\mathrm{dist}_{\sp^2}(x,\gamma), & \text{if } x \in \mathcal V_+ \cup \gamma, \\[2mm]
	-\,\mathrm{dist}_{\sp^2}(x,\gamma), & \text{if } x \in \mathcal V_- .
\end{cases}
\]
Note that, possibly after shrinking the neighborhood $\mathcal V$ of $p$, we may assume that the signed distance function $d$ is real analytic on $\mathcal V$; see \cite{Kr2}.
\begin{lemma}\label{lemma_Taylor}
There is a neighborhood of $p$ contained in $\mathcal{V}$ where the following expansions hold:
\[
\begin{split}
&	\abs{\nabla u}^2 = f(M)^2  \left( 1+ \kappa (p) d \right) d^2+ O(d^4), \\
& 	\bar{W} = f (M)^2 d^2 \left(1+  \left(\frac{\kappa (p)}{3} \pm \frac{2 \bar R}{3 \sqrt{1-\bar R^2}}\right)d\right)+ O (d^4),
\end{split}
\]
where in the second expansion we take the positive sign if $\overline{\tau} (\bar \U) \geq \tau_f^0 (M)$ and the negative one if $\overline{\tau} (\bar \U) < \tau_f^0 (M)$.  
\end{lemma}
\begin{proof}
	We follow \cite[Proposition 4.2]{EMa}.  Let $t = \Psi - \bar{R}$, with $\Psi$ denoting the pseudo-radial function defined in Definition \ref{psudoRadialFunctions} and $\bar R$ being the critical height of $(\bar \U, \bar u)$. We also recall that $\bar U$ denotes the solution to \eqref{ODE}-\eqref{CauchyData} such that $u = \bar U \circ \Psi$ in $\U$. Then, using the equation \eqref{ODE}, we compute the Taylor expansion of $M-\bar U \circ \Psi$ with respect to the variable $\Psi$, obtaining
	\begin{equation}\label{Taylor1}
		M-u = \frac{f (M)}{2(1-\bar R^2)} t^2+\frac{2 \bar R f(M)}{3 (1-\bar R^2)^2} t^3+ O(t^4).
	\end{equation}
	Next, recalling that $\bar W = (1-\Psi^2) \bar{U}' (\Psi)^2$, we compute the following Taylor expansion in a neighborhood of $t=0$,
	\begin{equation}\label{Taylor2}
		\frac{\bar{W}}{M-u} = 2 f(M) + \frac{4 \bar R f(M)}{3 (1-\bar R^2)}t + O (t^2).
	\end{equation}
	Finally, using \eqref{Taylor1} we get that 
	\[
	t= \pm \sqrt{\frac{2 (1-\bar R^2)}{f(M)}} \sqrt{M-u}+ O(M-u),
	\]
	where we take the positive sign if $\overline{\tau} (\bar \U) \geq \tau_f^0 (M)$ and the negative sign if $\overline{\tau} (\bar \U) < \tau_f^0 (M)$. Substituting this into \eqref{Taylor2} we conclude that
	\begin{equation}\label{Taylor3}
		\bar{W} = 2 f(M)(M-u) \pm \frac{4\bar R \sqrt{2 f (M)}}{3 \sqrt{1-\bar R^2}} (M-u)^{3/2}+ O ((M-u)^2).
	\end{equation}
	Now,  \cite[Theorem 3.1]{Chr} implies that in a neighborhood of $p$ contained in $\mathcal{V}$ the following expansion hold:
	\[
	u = M- \frac{f(M)}{2} d^2 - \frac{ f(M)}{6} \kappa(p) d^3 + O(d^4).
	\]
	Then, we get that 
	\begin{equation}\label{RhoExpansionU}
		\nabla u = - f (M) \left(1+\frac{\kappa (p)}{2}d + O(d^2)\right)d \cdot \nabla d,
	\end{equation}
	so, taking into account that $\abs{\nabla d} =1$, we conclude that
	\begin{equation}\label{RhoExpansionW}
		\abs{\nabla u}^2 = f(M)^2  \left( 1+ \kappa (p) d \right) d^2+ O(d^4).
	\end{equation}
	In addition, substituting \eqref{RhoExpansionU} into \eqref{Taylor3} yield
	\begin{equation}\label{RhoExpansionModelW}
		\bar{W} = f (M)^2 d^2 \left(1+  \left(\frac{\kappa (p)}{3} \pm \frac{2 \bar R}{3 \sqrt{1-\bar R^2}}\right)d\right)+ O (d^4),
	\end{equation}
	where we recall that the positive sign is taken if $\overline{\tau} (\bar \U) \geq \tau_f^0 (M)$ and the negative is taken otherwise. This finishes the proof of the lemma. 
\end{proof}

\section*{Acknowledgments}

The authors would like to thank H. Rosenberg for his insightful comments and valuable suggestions, which have significantly improved the original version of this paper.

The first author is partially supported by the {\it Maria de Maeztu} Excellence Unit IMAG, reference CEX2020-001105-M, funded by MCINN/AEI/10.13039/ 501100011033/CEX2020-001105-M and  Spanish MIC Grant PID2024-160586NB-I00.

The second author is partially supported by the {\it Maria de Maeztu} Excellence Unit IMAG, reference CEX2020-001105-M, funded by MCINN/AEI/10.13039/ 501100011033/CEX2020-001105-M, and Spanish MIC Grant PID2023-150727NB.I00.

\end{document}